\documentclass[epsf]{article}
\usepackage{amsmath,amsthm,amsfonts,latexsym,amscd}
\usepackage[mathscr]{eucal}
\usepackage{eufrak}
\usepackage{graphicx}
\title{Bumpy Riemannian Metrics and Closed Parametrized Minimal Surfaces in Riemannian Manifolds\footnote{This is a revision of the article, {\em Bumpy Riemannian metrics and closed parametrized minimal surfaces in Riemannian manifolds}, Trans. AMS {\bf 358} (2006), 5193-5256, with erratum in Trans. AMS {\bf 359} (2007), 5117-5123.}}
\author{John Douglas Moore\\Department of 
Mathematics\\University of California\\
Santa Barbara, CA, USA 93106\\e-mail: moore@math.ucsb.edu}
\date{}

\begin{document}

\maketitle

\begin{abstract}
This article is concerned with conformal harmonic maps $f:\Sigma \rightarrow M$, where $\Sigma $ is a closed Riemann surface and $M$ is a compact Riemannian manifold of dimension at least four.  We show that when the ambient manifold $M$ is given a generic metric, all prime closed parametrized minimal surfaces are free of branch points, and are as Morse nondegenerate as allowed by the group of complex automorphisms of $\Sigma $.
\end{abstract}

\section{Introduction}\label{S:intro}

This article is devoted to providing part of the foundation needed for a partial Morse theory for parametrized minimal surfaces in Riemannian manifolds which should parallel the well-known Morse theory of smooth closed geodesics in a compact Riemannian manifold $M$.  Recall that the theory of closed geodesics is concerned with the action integral
$$J : L^2_1(S^1,M) \rightarrow {\mathbb R} \qquad \hbox{defined by} \qquad J(\gamma ) = \frac{1}{2} \int _{S^1} \langle \gamma '(t), \gamma '(t)\rangle dt,$$
where $\langle \cdot , \cdot \rangle $ denotes a Riemannian metric on $M$ and $t$ is the usual coordinate on $S^1$, which is regarded as the unit interval $[0,1]$ with endpoints identified.  Moreover, $L^2_1(S^1,M)$ denotes the Hilbert manifold of $L^2_1$-maps of $S^1$ into $M$, maps which are $L^2$ and have $L^2$ first derivatives, as described for example in \cite{K} or \cite{H}.  The function $J$ is invariant under a continuous right action of the circle group $S^1$ on $L^2_1(S^1,M)$,
$$\phi : L^2_1(S^1,M) \times S^1 \longrightarrow L^2_1(S^1,M), \qquad \phi (\gamma , s)(t) = \gamma (t + s),$$
and therefore whenever any point is critical for $J$, so is the entire $S^1$-orbit.  Thus nonconstant critical points for $J$ are never Morse nondegenerate in the usual sense.  However, the Bumpy Metric Theorem of Abraham \cite{A} states that for generic choice of Riemannian metric, all nonconstant smooth closed geodesics lie on nondegenerate critical submanifolds of dimension one, and have the property that their only Jacobi fields are those generated by the $S^1$-action.

The purpose of this article is to provide a similar bumpy metric theorem for the energy function on maps from compact closed Riemann surfaces into a compact Riemannian manifold $M$ of dimension at least four.  Thus suppose that $\Sigma $ is a compact Riemann surface without boundary of genus $g$, that $\hbox{Map}(\Sigma ,M)$ is a completion of the space of smooth maps from $\Sigma $ to $M$ with respect to a suitable Sobolev norm, and that ${\mathcal T}$ is the Teichm\"uller space of marked conformal structures on $\Sigma $.  Then the {\em energy\/} $E : \hbox{Map}(\Sigma ,M) \times {\mathcal T} \rightarrow {\mathbb R}$ is defined by
$$E(f, \omega ) = \frac{1}{2}\int _{\Sigma}\left[ \frac{\partial f}{\partial x} \cdot \frac{\partial f}{\partial x} + \frac{\partial f}{\partial y} \cdot \frac{\partial f}{\partial y} \right] dx dy = \frac{1}{2} \int_\Sigma |df|^2 dA, $$
where $(x,y)$ are local conformal coordinates for the conformal structure $\omega \in {\mathcal T}$ (the integrand being independent of the choice) and $|df|$ and $dA$ represent the norm of $df$ and the area element with respect to any metric on $\Sigma $ within the conformal structure $\omega $.  For a fixed choice of $\omega \in {\mathcal T}$, we define
$$E_\omega : \hbox{Map}(\Sigma ,M) \longrightarrow {\mathbb R} \qquad \hbox{by} \qquad E_\omega (f) = E(f,\omega ).$$
Critical points for $E_\omega $ are called harmonic maps (or $\omega $-harmonic maps when we want to emphasize the specific conformal structure), while critical points for $E$ are (weakly) conformal harmonic maps, or equivalently minimal surfaces, possibly containing branch points.  In Osserman's survey \cite{O}, in which the ambient manifold $M$ is Euclidean space, critical points of $E$ are called {\em generalized\/} minimal surfaces when branch points are allowed, but we will use the simpler term parametrized minimal surfaces for such objects.

In some cases one can prove existence of critical points for these functions by means of a perturbed energy introduced by Sacks and Uhlenbeck \cite{SU1}, \cite{SU2}.  Thus for $\omega \in {\mathcal T}$ and $\alpha > 1$, we define a function $E_{\alpha ,\omega} : \hbox{Map}(\Sigma ,M) \rightarrow {\mathbb R}$ by
\begin{equation}E_{\alpha ,\omega} (f) = \frac{1}{2} \int \int _{\Sigma } (1 + |df|^2)^\alpha dA, \label{E:alphaomega}\end{equation}
where the Riemannian metric on $\Sigma $ used in the integrand is chosen to be a metric within the conformal equivalence class defined by $\omega $.  Although the ability to choose a metric on $\Sigma $ provides potential flexibility, we will normalize the metric on $\Sigma $ to be the unique constant curvature metric of total area one.  We call $E_{\alpha ,\omega}$ the {\em $(\alpha ,\omega )$-energy\/}.  We recall some important properties of the function $E_{\alpha ,\omega }$ when $ \hbox{Map}(\Sigma ,M)$ is the completion with respect to the Sobolev $L^p_1$-norm, where $p = 2\alpha $: The function $E_{\alpha ,\omega}$ has continuous derivatives up to order two, it satisfies condition C of Palais and Smale and its critical points are as smooth as the metric on $M$.  For fixed choice of $\alpha > 1$, one can use techniques of global analysis \cite{U2} to prove existence of critical points for $E_{\alpha ,\omega}$ corresponding to minimax constraints.  As $\alpha \rightarrow 1$, one can often show that a \lq\lq minimax sequence" will approach a harmonic map in the limit, together with a finite number of \lq\lq bubbles," which are harmonic maps of two-spheres.

There are several difficulties in constructing a Morse theory for the energy function $E$ via this approach.  First, bubbling interferes with the most straightforward version of the Morse inequalities (but also suggests a procedure for analyzing to what extent the Morse inequalities fail).  Second, minimax sequences for various constraints might tend to the boundary of Teichm\"uller space, a phenomenon already observed in the Plateau problem for mappings of Riemann surfaces with boundary into Euclidean space ${\mathbb R}^3$---thus a minimizing sequence might have a handle pinch off in the limit, for example.  Third, a given minimal surface has infinitely many different markings, giving rise to an orbit of critical points in $\hbox{Map}(\Sigma ,M) \times {\mathcal T}$ under an action of the mapping class group $\Gamma $.  (We describe this action at the beginning of \S\ref{S:diffeogroup}.)  This action leaves the energy invariant, so the energy actually descends to a function on the quotient
$$E :  (\hbox{Map}(\Sigma ,M) \times {\mathcal T} )/\Gamma \rightarrow {\mathbb R},$$
where the domain projects to the moduli space ${\mathcal M} = {\mathcal T}/\Gamma $ with typical fiber $\hbox{Map}(\Sigma , M)$.  One really wants a partial Morse theory to estimate the number of critical points of given index in the quotient, and although the base space ${\mathcal M}$ is diffeomorphic to a ball when $g=1$, it is only an orbifold for large genus. 

Finally, just as in the theory of smooth closed geodesics, multiple covers of a given minimal surface should not count as geometrically distinct, even though they will appear as distinct critical points for $E$.  Indeed, this difficulty is more pronounced for maps from surfaces because the covers might be branched:  Suppose, for example, that $h : \Sigma _2 \rightarrow M$ is a smooth harmonic map and that $g : \Sigma _1 \rightarrow \Sigma _2$ is a nontrivial conformal branched cover.  (If $\Sigma _2 = S^2$, the Riemann sphere, we can regard $g$ as a meromorphic function on $\Sigma _1$.)  Then the composition $f = h \circ g : \Sigma _1 \rightarrow M$ is called a {\em branched cover\/} of the harmonic map $h : \Sigma _2 \rightarrow M$.  Any compact simply connected Riemannian manifold must contain minimal two-spheres \cite{SU1}, so branched covers are often present.

In spite of these difficulties, it is sometimes possible to control bubbling and use critical point theory for the $\alpha $-energy to prove existence of harmonic maps or minimal surfaces.  One application is the theorem of Sacks and Uhlenbeck \cite{SU1} that any compact simply connected Riemannian manifold contains at least one minimal two-sphere.  Another application involving equivariant Morse theory for minimal two-spheres is presented in \cite{M}, in which we use bounds on the energy to prevent bubbling, which yields existence of many minimal two-spheres of low energy in suitably pinched $n$-spheres, which satisfy appropriate hypotheses.  For further results of this nature, we need a bumpy metric theorem for parametrized minimal surfaces as a foundation.

In view of transversality techniques of Smale \cite{Sm2} and Uhlenbeck \cite{U1}, it is not surprising that such a bumpy metric theorem should hold when properly formulated.  Indeed, White proves a bumpy metric theorem for imbedded minimal submanifolds of arbitrary dimension in a Riemannian ambient manifold \cite{W2}.  However, for the theory of parametrized minimal surfaces, it is crucial that we allow for the possibility of branch points.  Branch points always occur on the nontrivial branched covers of parametrized surfaces of lower genus mentioned before; such branched covers can have arbitrarily large spaces of Jacobi fields that cannot be removed by generic perturbations of the metric on $M$.  Removing branch points by perturbation on the remaining surfaces is the key issue which needs to be resolved for establishing a bumpy metric theorem in the parametrized context which is strong enough for our projected applications.  (A related transversality theorem for Gromov's theory of pseudoholomorphic curves is presented in the book of McDuff and Salamon \cite{McS}.)  

Another complication to the parametrized theory is that when the genus of $\Sigma $ is small, the energy function $E$ is invariant under a nontrivial connected Lie group $G$ of symmetries.  In the case where $\Sigma $ is the torus, the two-dimensional Lie group $G = S^1 \times S^1$ acts on $\hbox{Map}(\Sigma ,M)$ preserving the energy, which implies that nonconstant minimal tori lie on orbits of critical points of dimension two.  In the case where $\Sigma $ is the Riemann sphere $S^2$, the energy is invariant under the action of the six-dimensional Lie group $G = PSL(2,{\mathbb C})$ via linear fractional transformations on the domain.  This action preserves the energy and hence nonconstant minimal two-spheres lie on orbits of critical points having dimension six.  Since $PSL(2,{\mathbb C})$ is noncompact, it is often convenient (in accordance with the procedure used in \cite{M}) to replace $\hbox{Map}(S^2,M)$ by the subspace of maps with \lq\lq center of mass zero,"
\begin{equation} \hbox{Map}_0(S^2 ,M) = \left\{ f \in \hbox{Map}(S^2 ,M) : \int _{S^2} X |df|^2 dA = 0 \right\},\label{eq:centerofmass}\end{equation}
where $X : S^2 \rightarrow {\mathbb R}^3$ is the standard inclusion onto the unit sphere in ${\mathbb R}^3$.  This allows us to reduce the symmetry group of the energy function from $PSL(2,{\mathbb C})$ to its maximal compact subgroup $SO(3)$.  (This is particularly useful for the approach to existence via $\alpha $-energy, because the $\alpha $-energy is only invariant under the smaller group $SO(3)$.)  Because of the group actions, the most one could hope for in the cases of minimal spheres and tori (as well as nonorientable projective planes and Klein bottles) is that in analogy with the closed geodesic problem, the critical points for energy lie on nondegenerate critical submanifolds for generic choice of metric on the ambient manifold.  Let us recall Bott's definition: 

\vskip .1 in
\noindent
{\bf Definition.} let $F : {\mathcal M} \rightarrow {\mathbb R}$ be a $C^2$ function on a smooth manifold ${\mathcal M}$.  A {\em nondegenerate critical submanifold\/} of ${\mathcal M}$ is a finite-dimensional submanifold $N \subset {\mathcal M}$ such that every $f \in N$ is a critical point for $F$ and
\begin{equation} f \in N \quad \Rightarrow \quad T_fN = \{ X \in T_f{\mathcal M} : d^2F(f)(X,Y) = 0 \hbox{ for all } Y \in T_f{\mathcal M} \}. \label{E:noncritsub} \end{equation}
Here $d^2F(f)$ is the Hessian of $F$ at the critical point, and elements $X \in T_f{\mathcal M}$ which satisfy the condition on the right-hand side of (\ref{E:noncritsub}) are called {\em Jacobi fields\/}.

\vskip .1 in
\noindent
{\bf Definition.}  We say that a parametrized minimal surface $f: \Sigma \rightarrow M$ is {\em prime\/} if it is nonconstant and is not a nontrivial cover (possibly branched) of another parametrized minimal surface $f_0: \Sigma _0\rightarrow M$ of lower energy, where we explicitly allow the surface $\Sigma _0$ to be nonorientable.

\vskip .1 in
\noindent 
{\bf Definition.}  By a {\em generic choice of Riemannian metric\/}, we mean a metric belonging to a countable intersection of open dense subsets of the spaces of $L^2_k$ Riemannian metrics on $M$, where $k$ ranges over the positive integers.

\vskip .1 in
\noindent 
With these definitions in place, we can now state our main theorem for critical points of the energy function $E$:

\vskip .1 in
\noindent
{\bf Main Theorem.} {\sl Suppose that $M$ is a compact connected smooth manifold of dimension at least four.  Then
\begin{enumerate}
\item For a generic choice of Riemannian metric on $M$, every prime minimal two-sphere $f:S^2 \rightarrow M$ is free of branch points and lies on a nondegenerate critical submanifold of dimension six which is an orbit for the G-action, where $G = PSL(2,{\mathbb C})$.
\item For a generic choice of Riemannian metric on $M$, every prime minimal two-torus $f:T^2 \rightarrow M$ is free of branch points and lies on a nondegenerate critical submanifold of dimension two which is an orbit for the G-action, where $G = S^1 \times S^1$.
\item For a generic choice of Riemannian metric on $M$, every prime oriented minimal surface $f: \Sigma \rightarrow M$ of genus at least two is free of branch points and is Morse nondegenerate in the usual sense.
\end{enumerate}
Moreover, in each of the cases, the prime minimal immersions $f : \Sigma \rightarrow M$ have transversal crossings.}

\vskip .1 in
\noindent
In particular, if $M$ has a generic metric, prime oriented closed parametrized minimal surfaces within $M$ are immersions.  For simplicity, we will focus on oriented minimal surfaces throughout most of the article, although the argument can be modified to give a bumpy metric theorem for nonorientable minimal surfaces by arguments we sketch in \S\ref{S:nonorientable}.

The freedom from branch points asserted in the Main Theorem provides an analog of a well-known theorem of B\"ohme and Tromba (\cite{BoTr1} and \cite {BoTr2}) that parametrized minimal surfaces bounded by a generic smooth closed Jordan curve in ${\mathbb R}^n$ have no branch points when $n \geq 4$.  Moreover, the Main Theorem shows that the nondegeneracy hypothesis of Theorems 1 and 2 of \cite{M} hold for generic metrics.

Here is a brief outline of the remainder of the article.  We first review the parametric theory of closed harmonic and minimal surfaces in Riemannian manifolds in \S\ref{S:prelim}, and explain the relationship with the theory of Riemann surfaces.  In \S\ref{S:harmonic} we prove a bumpy metric theorem for $\omega $-harmonic immersions via a proof schema which will be used as a model for the proof of the Main Theorem.  \S\ref{S:two-variableenergy}, we describe second variation for a two-variable energy in which the metric on the range is allowed to vary.  Although this reduces to the usual second variation formula for minimal surfaces when restricted to normal variations, the more general formula uncovers additional tangential Jacobi fields when applied to minimal surfaces with branch points.  The number of such Jacobi fields is calculated by means of the Riemann-Roch Theorem.  We then establish the expected the easy bumpy metric theorem for immersed minimal surfaces.

The key idea behind the Main Theorem is that some of the tangential Jacobi fields generated by the branch points can be perturbed away by suitable perturbations in the metric.  The reason one expects such perturbations is described at the beginning of \S\ref{S:eliminatesimple}, and we then go on to construct explicit tangential variations that would perturb away all tangential Jacobi fields not generated by the action of $G$ for minimal surfaces with only simple branch points and with self-intersection set consisting only of double points.  This gives a bumpy metric theorem for simple branch points as we describe in \S\ref{S:simple}.

To proceed further, we need to examine the self-intersection set of arbitrary minimal immersion in \S\ref{S:selfintersection}, modifying slightly well-known results on self-intersections of minimal surfaces in three-dimensional manifolds.  These extend the easier fact that prime parametrized minimal surfaces are automatically somewhere injective, which follows from a theorem of Aronszajn.  Abraham's original sketch for the bumpy metric theorem for geodesics \cite{A0} relied on a transversal density theorem presented in \cite{AR}.  We use that theorem in \S\ref{S:genericselfintersection} to prove that for generic metric on $M$, any parametrized minimal surface $f : \Sigma \rightarrow M$ is an immersion with transversal crossing on the complement of its branch locus.  But to establish the hypotheses of the transversality theory, we need that prime minimal surfaces lie on countably many submanifolds parametrized by metrics in \S\ref{S:submanifoldsoftildeM}.  A modification of this argument then finished the proof of the Main Theorem in \S\ref{S:proof}.  Even though we do not know how to perturb away all the tangential Jacobi fields, the perturbations we can construct allow us to show that for any branch type with nontrivial branch locus, the minimal surfaces of this branch type lie in a countable sequence of manifolds which project to the space of metrics via a Fredholm map with Fredholm index $\leq 2$.  The Main Theorem then follows from the Sard-Smale Theorem.

In \S\ref{S:nonorientable} we describe how the argument for the Main~Theorem can be modified to treat prime nonorientable minimal surfaces, while a brief final section (\S\ref{S:families}) describes how the Main Theorem can be extended to generic one-parameter families of metrics on $M$.

Throughout this article, we will assume some familiarity with the techniques of global analysis on infinite-dimensional manifolds, presented in \cite{A0} and with more details in \cite{Pal1.5}.  More leisurely introductions are found in \cite{L} and \cite{AMR}.  One of the core theorems to keep in mind is nicely stated in a survey article by Eells \cite{Ee} (see page~780):  If $S$ and $M$ are finite-dimensional smooth manifolds and $C^k(S ,M)$ denotes the completion of the space of smooth maps from $S$ to $M$ with the $C^k$ norm, a smooth infinite-dimensional manifold, then
\begin{equation} \Phi : C^{k+s}(M,N) \times C^k(S,M) \rightarrow C^k(S,N), \qquad \Phi (g,f) = g \circ f, \label{E:lossofderivatives} \end{equation}
is a $C^s$ map.  In conjunction with the Sobolev Lemmas, this can be extended to yield smoothness of composition in Sobolev completions, and specializes to yield various versions of the so-called $\alpha $- and $\omega $-Lemmas of \cite{A0}.  One version of the $\omega $-Lemma for $L^p_k$ spaces states that if $\Sigma $ is a smooth surface and $M_1$ and $M_2$ are smooth manifolds of arbitrary dimension, a $C^\infty$ map $g : M_1 \rightarrow M_2$ induces a $C^\infty$ map
$$\omega _ g : L^p_k(\Sigma ,M_1) \longrightarrow L^p_k(\Sigma ,M_2), \qquad \omega _g(f) = g \circ f,$$
where $L^p_k(\Sigma ,M)$ denotes the completion of $C^\infty $ maps with respect to the $L^p_k$ norm, for any $p$ and $k$ satisfying $k - 2/p > 0$.  On the other hand, if we take $S$ to be a one-point space, (\ref{E:lossofderivatives}) only shows that the evaluation map
$$\hbox{ev}: C^{s}(M,N) \times M \longrightarrow N, \qquad \hbox{ev}(f,p) = f(p)$$
is $C^s$.  The loss of derivatives implicit in (\ref{E:lossofderivatives}) often forces us to work in a Sobolev completion $L^2_k(\Sigma ,M)$ in which the number $k$ of generalized derivatives is high.  For example, in arguments involving branch points, the needed value of $k$ is sometimes determined by the largest possible branching order, which in turn depends on an upper bound $E_0$ for the energy, as we will see in the next section.  This is harmless for projected applications, since regularity theorems imply that critical points for $E$ and $E_\alpha $ are always $C^\infty$.  It would be an easy exercise to rephrase the argument in terms of Fr\'echet manifolds as described in \cite{Ham}.

The author thanks Vincent Borrelli of the Universit\'e Claude Bernard of Lyon, France for helpful discussions on the theory presented here during a visit to Lyon in December of 2004.  He also thanks the referee for numerous questions and suggestions that helped him make many major improvements to the article.

\section{Preliminaries}
\label{S:prelim}

We recall a few basic concepts from the theory of harmonic and minimal surfaces in Riemannian manifolds (some of which are presented in more detail in \cite{MM}).  We focus first on the $\omega $-energy $E_\omega $ for a fixed choice of Riemannian metric $g = \langle \cdot , \cdot \rangle$ on $M$.  The first derivative of $E_\omega $ is given by the formula
\begin{equation} dE_\omega (f)(X) = \int _\Sigma \langle F_\omega (f,g), X \rangle dA, \label{E:firstvariationfore} \end{equation}
where $F_\omega( \cdot ,g) = 0$ is the Euler-Lagrange equation for the variational problem.  If $f$ is a critical point for $E_\omega $, we can differentiate once again to obtain the Hessian,
\begin{equation}d^2E_\omega (f)(X,Y) = \int _\Sigma \langle D_1F_\omega (f,g)(X), Y \rangle dA = \int _\Sigma \langle L_{(f,g)}(X), Y \rangle dA, \label{E:jacob} \end{equation}
where $D_1F_\omega $ denotes the derivative with respect to the variable $f \in \hbox{Map}(\Sigma ,M)$ and $L_{(f,g)}$ is the Jacobi operator, a formally self-adjoint second-order elliptic differential operator which acts on sections of the pullback $f^*TM$ of the tangent bundle to $M$.  Such sections are the elements of the tangent space to $\hbox{Map}(\Sigma ,M)$ at $f$.  Elements $X$ of this tangent space $T_f\hbox{Map}(\Sigma ,M)$ which satisfy the {\em Jacobi equation\/} $L_{(f,g)}(X) = 0$ are called {\em Jacobi fields\/} at $f$.  Note that since $L_{(f,g)} = D_1F_\omega (f,g)$, the Jacobi equation is just the linearization of the Euler-Lagrange equation.

In terms of a local complex coordinate $z = x + iy$ on $\Sigma $, with $\lambda ^2$ denoting the conformal factor (so that the area element is given by $dA = \lambda ^2 dxdy$), the Euler-Lagrange equation for an $\omega $-harmonic map $f : \Sigma \rightarrow M$ is simply
\begin{equation} - \frac{1}{4\lambda ^2}\frac{D}{\partial \bar z}\left(\frac{\partial f}{\partial z}\right) = 0, \qquad \hbox{where} \qquad \frac{\partial f}{\partial z} = f_*\left(\frac{\partial }{\partial z}\right) \label{E:harmonic}\end{equation}
is regarded as a section of the complex bundle ${\bf E} = f^*TM \otimes {\mathbb C}$ over $\Sigma $, and $D$ denotes the pullback of the Levi-Civita connection of $M$ to ${\bf E}$.  Thus the Euler-Lagrange equation just says that the section $\partial f/\partial z$ is holomorphic, when ${\bf E}$ is given the holomorphic structure (guaranteed to exist by a theorem of Koszul and Malgrange) such that if $V$ is an element of $\Gamma ({\bf E})$, the space of smooth sections of ${\bf E}$, then
$$\hbox{$V$ is holomorphic} \qquad \Leftrightarrow \qquad \frac{D V}{\partial \bar z} = 0.$$
Although the locally defined holomorphic section
$$\frac{\partial f}{\partial z} \quad \hbox{may have isolated zeros, the induced map} \quad \left[\frac{\partial f}{\partial z} \right] : \Sigma \rightarrow {\mathbb P}({\bf E})$$
is independent of choice of local coordinate $z$ and is a globally defined section on $\Sigma $, where ${\mathbb P}({\bf E})$ denotes the bundle of complex projective lines in fibers of ${\bf E}$.  This globally defined family of lines defines a holomorphic line bundle ${\bf L}$ within ${\bf E}$.  A point $p \in \Sigma$ is said to be a {\em branch point\/}  if $(\partial f/\partial z)(p) = 0$.  In this case, if $z$ is a holomorphic coordinate centered at $p$, we can write
$$\frac{\partial f}{\partial z} = z^\nu g(z),$$
where $g$ is a holomorphic section of ${\bf E}$ such that $g(p) \neq 0$.  The integer $\nu $ is called the {\em branching order\/} of $f$ at $p$.

Observe that the vector-valued differential
$$\partial f = \frac{\partial f}{\partial z}dz \quad \hbox{is a globally defined holomorphic section of} \quad {\bf L} \otimes {\bf K},$$
where ${\bf K} = T^*_h\Sigma$ is the holomorphic cotangent bundle of $\Sigma $, also known as the canonical bundle, and this section has zeros precisely at the branch points of $f$.  In particular, if $f$ has no branch points, this section trivializes ${\bf L} \otimes {\bf K}$, so in this case we can conclude that ${\bf L}$ must be isomorphic to the holomorphic tangent bundle ${\bf K}^{-1}$ to $\Sigma $.  In the general case, we let $\nu _p(f)$ denote the branching order of the harmonic map $f$ at $p$ and call the finite sum
$$D(f) = \sum \{ \nu _p(f) p : \nu_p(f) > 0 \}$$
the {\em divisor\/} of the harmonic map $f$.  By standard arguments in Riemann surface theory (as presented in \cite{Gu}, \S 7), we can determine the isomorphism class of the holomorphic bundle ${\bf L}$, the result being that
$${\bf L} \otimes {\bf K} \cong \zeta_{p_1}^{\nu _{p_1}} \otimes \cdots \otimes \zeta_{p_n}^{\nu _{p_n}}, \quad \hbox{and hence} \quad {\bf L} \cong {\bf K}^{-1} \otimes \zeta_{p_1}^{\nu _{p_1}} \otimes \cdots \otimes \zeta_{p_n}^{\nu _{p_n}},$$
where $p_1, \ldots ,p_n$ are the branch points of $f$ and $\zeta _p$ is the holomorphic point bundle at $p$.  The Kronecker product of the first Chern class of ${\bf L}$ and the fundamental class of $\Sigma $ is given by the formula
\begin{equation}\langle c_1({\bf L}),[\Sigma ]\rangle = \langle c_1({\bf K}^{-1}),[\Sigma ]\rangle + \sum \nu _p = (2-2g) + \sum \nu _p, \label{E:branch}\end{equation}
where $\sum \nu _p$ is the total branching order of $f$.

In the case of the sphere $S^2 = {\mathbb C} \cup \{ \infty \}$ we can take the standard complex coordinate $z$ on ${\mathbb C}$ and check that
$$\frac{\partial f}{\partial z} \quad \hbox {extends to a holomorphic vector field which vanishes at $\infty$.}$$
Thus the holomorphic function
$$\left\langle \frac{\partial f}{\partial z}, \frac{\partial f}{\partial z} \right\rangle \quad \hbox {is globally defined and vanishes at $\infty$,}$$
and must therefore be identically zero.  By dividing into real and imaginary parts, we conclude that
$$\left\langle \frac{\partial f}{\partial x},\frac{\partial f}{\partial x}\right\rangle = \left\langle \frac{\partial f}{\partial y}, \frac{\partial f}{\partial y}\right\rangle, \qquad \left\langle \frac{\partial f }{\partial x},\frac{\partial f}{\partial y}\right\rangle = 0,$$
which is just the condition that the harmonic two-sphere $f$ be conformal.  (Some authors would say weakly conformal when $f$ has branch points.) 

In the case of the torus, the coordinate vector field
$$\frac{\partial f}{\partial z} \qquad \hbox{is globally defined and} \qquad \left\langle \frac{\partial f}{\partial z}, \frac{\partial f}{\partial z} \right\rangle$$
is a globally defined holomorphic function, which must therefore be constant by the maximum modulus principle.  The constant is zero precisely when $f$ is conformal.  Since the constant is nonzero for a nonconformal harmonic torus, such a torus cannot have branch points.

To deal with Riemann surfaces of genus $g \geq 2$ we utilize the Hopf differential
$$\Omega _f = \langle \partial f, \partial f \rangle = \left\langle \frac{\partial f}{\partial z}, \frac{\partial f}{\partial z} \right\rangle dz^2,$$
a quadratic differential which is holomorphic when $f$ is $\omega $-harmonic and vanishes when $f$ is also conformal.  If $\Sigma $ has genus $g \geq 2$, all nonzero holomorphic quadratic differentials on $\Sigma $ must have the same number of zeros counting multiplicity, namely
$$\langle c_1({\bf K}^2) ,[\Sigma ]\rangle = \langle c_1(T_h^*\Sigma \otimes T_h^*\Sigma ) ,[\Sigma ]\rangle = 4g-4.$$
On the other hand, if a nonconformal $\omega $-harmonic map $f:\Sigma \rightarrow M$ has a branch point of branching order $k$ at $p$, then the Hopf differential $\Omega _f$ must have a zero of order at least $2k$ at $p$, and hence
$$\nu (f) = (\hbox{total branching order of $f$}) \leq \frac{1}{2} (\hbox{number of zeroes of $\Omega _f$}) = 2g-2.$$
Thus a nonconformal harmonic surface of genus $g \geq 2$ cannot have more than $2g-2$ branch points, counting multiplicity.  Moreover,
$$\langle c_1({\bf L}) ,[\Sigma ]\rangle = 2 - 2g + \nu (f) \leq 0,$$
so ${\bf L}$ has no holomorphic sections, unless it is the trivial line bundle and $\nu (f) = 2g-2$.  The case where ${\bf L}$ is the trivial bundle can occur if $f = h \circ g$, where $g: \Sigma \rightarrow T^2$ is a branched cover and $h:T^2 \rightarrow M$ is a nonconformal harmonic torus.

Returning to the case of general genus, we note that if $f$ is not only harmonic but also conformal, the rank of $df_p$ can never be one, and $f$ is an immersion except at branch points.  This contrasts with the fact that nonconformal harmonic maps can have points at which $df_p$ has rank one; this happens for example when $f$ is a torus parametrization of a geodesic.  If $f$ is an immersion except for branch points, the real and imaginary parts of sections of the line bundle ${\bf L}$ generate a two-dimensional subbundle $(f^*TM)^\top$ of the pullback tangent bundle $f^*TM$.  This subbundle has an induced metric and orientation which determine an almost complex structure
$$J : \Gamma ((f^*TM)^\top) \rightarrow \Gamma ((f^*TM)^\top),$$
where $\Gamma ((f^*TM)^\top)$ denotes the space of smooth sections of $(f^*TM)^\top$, such that
$${\bf L} \cong \{ v \in (f^*TM)^\top\otimes {\mathbb C} : J(v) = iv \}, \quad {\bf L^*} \cong \{ v \in (f^*TM)^\top\otimes {\mathbb C} : J(v) = -iv \}.$$
More explicitly, the map $\tau : (f^*TM)^\top \rightarrow {\bf L}$ defined by
\begin{equation} \tau \left( M \frac{\partial f}{\partial x} + N \frac{\partial f}{\partial y}\right) = 2(M + i N) \frac{\partial f}{\partial z}\label{E:isomorphism}\end{equation}
is a complex linear isomorphism, the inverse of which is just the projection of a section of ${\bf L}$ to its real part, an isomorphism we will frequently use.

The bundle $(f^*TM)^\top$ has an orthogonal complement $(f^*TM)^\bot$ and we can define the second fundamental form
$$A(f) : \Gamma (f^*TM)^\top \times \Gamma (f^*TM)^\top \rightarrow \Gamma (f^*TM)^\bot \quad \hbox{by} \quad A(f)(X,Y) = (D_XY)^\bot,$$
with $D$ denoting once again the Levi-Civita connection.  It follows immediately from the equation for harmonic maps (\ref{E:harmonic}) that
\begin{equation} A(f)\left(\frac{\partial f}{\partial x},\frac{\partial f}{\partial x}\right) + A(f)\left(\frac{\partial f}{\partial y}, \frac{\partial f} {\partial y}\right) = 0. \label{E:secondfundamentalform}\end{equation}
If we normalize the Riemannian metric on the compact manifold $M$ so that the sectional curvatures of $M$ are $\leq 1$, the Gauss equations together with (\ref{E:secondfundamentalform}) imply that the curvature $K$ of $(f^*TM)^\top$ is also $\leq 1$, and hence that
\begin{equation} \langle c_1({\bf L}),[\Sigma ]\rangle = \frac{1}{2 \pi} \int _\Sigma KdA \leq \frac{1}{2 \pi} \int _\Sigma dA = \frac{\hbox{Area of $f$}}{2 \pi} \leq \frac{\hbox{Energy of $f$}}{2 \pi}. \label{E:branchpointest}\end{equation}
It follows from (\ref{E:branchpointest}) that minimal spheres of sufficiently small energy cannot have branch points, and the total branching order of a parametrized minimal surface of genus $g$ grows at most linearly with the energy.  Thus given a bound $E_0$ on the energy, we can bound the total branching order $\nu (f)$. 

Let $\hbox{Met}(M)$ denote the space of Riemannian metrics on $M$, completed with respect to a suitable Sobolev norm.  Given an element $\omega $ of the Teichm\"uller space ${\mathcal T}$ for the closed oriented surface $\Sigma $ of genus $g$, we let
\begin{multline*}{\mathcal S}_\omega = \{ (f,g) \in \hbox{Map}(\Sigma ,M) \times \hbox{Met}(M) : \hbox{$df_p$ is injective} \\ \hbox{for some $p \in \Sigma$ and $f$ is $\omega $-harmonic for the metric $g$} \}. \end{multline*}
We also let
\begin{multline}{\mathcal S} = \{ (f,\omega ,g) \in \hbox{Map}(\Sigma ,M) \times {\mathcal T} \times \hbox{Met}(M) : \\ \hbox{$f$ is both weakly conformal and $\omega $-harmonic for the metric $g$} \}. \label{E:mathcalS} \end{multline}
In analogy with the presentation in B\"ohme and Tromba \cite{BoTr1}, \cite{BoTr2}, we can construct stratifications of ${\mathcal S}_\omega$ and ${\mathcal S}$, the strata being
$${\mathcal S}_\omega^\nu = \{ (f,g) \in {\mathcal S}_\omega : \hbox{$f$ has total branching order $\nu $} \},$$
$${\mathcal S}^\nu = \{ (f,\omega ,g) \in {\mathcal S}: \hbox{$f$ has total branching order $\nu $} \},$$
where $\nu$ is a nonnegative integer.

Holomorphic sections of the line bundle ${\bf L}$ can be identified with meromorphic sections $X$ of the holomorphic tangent bundle $T_h\Sigma $ which have the property that
\begin{equation}(X) = (\hbox{divisor of $X$}) \geq - \sum \nu _p p = - (\hbox{divisor of $f$}).\label{eq:divisor}\end{equation}
The space ${\mathcal O}({\bf L})$ of such sections is easily determined when $\Sigma $ has genus either zero or one.  If $\Sigma $ has genus zero, it follows from the Riemann-Roch theorem that
$$(f,g) \in {\mathcal S}_\omega^\nu \qquad \Rightarrow \qquad \hbox{dim}_{\mathbb C} {\mathcal O}({\bf L}) = \nu + 3,$$
there being only one conformal structure $\omega $ in this case.  If $\Sigma $ has genus one, the lowest degree stratum ${\mathcal S}^0$ only occurs when ${\bf L}$ is trivial and $\hbox{dim}_{\mathbb C} {\mathcal O}({\bf L}) = 1$, as forced by the action of $S^1 \times S^1$.  Moreover, if $\nu \geq 1$, it follows from the Riemann-Roch theorem that
$$(f,g) \in {\mathcal S}_\omega^\nu \quad \hbox{or} \quad (f,\omega ,g) \in {\mathcal S}^\nu \qquad \Rightarrow \qquad \hbox{dim}_{\mathbb C} {\mathcal O}({\bf L}) = \nu.$$

The Jacobi operator, defined in terms of the second variation by (\ref{E:jacob}), plays a central role, so we review the complex version of this operator as presented in \cite{MM}.  It follows from the usual second variation formula for harmonic maps that the Jacobi operator has the explicit expression
\begin{equation} L_{(f,g)}( \cdot ) = - \frac{1}{\lambda ^2}\left[\frac{D}{\partial x} \circ \frac{D}{\partial x} + \frac{D}{\partial y} \circ \frac{D}{\partial y} + R \left( \cdot, \frac{\partial f} {\partial x} \right)\frac{\partial f} {\partial x} + R \left( \cdot, \frac{\partial f} {\partial y} \right)\frac{\partial f} {\partial y}\right], \label{E:Jacobiharmonicmap}\end{equation}
where $R$ is the Riemann-Christoffel curvature tensor of $M$.  Straightforward calculations show that on the one hand,
\begin{multline*} 4\frac{D}{\partial z} \circ \frac{D}{\partial \bar z} = \frac{D}{\partial x} \circ \frac{D}{\partial x} + \frac{D}{\partial y} \circ \frac{D}{\partial y} + \sqrt{-1} \left(\frac{D}{\partial x} \circ \frac{D}{\partial y} - \frac{D}{\partial y} \circ \frac{D}{\partial x}\right) \\ = \frac{D}{\partial x} \circ \frac{D}{\partial x} + \frac{D}{\partial y} \circ \frac{D}{\partial  y} + \sqrt{-1} R \left( \frac{\partial f}{\partial x}, \frac{\partial f}{\partial y}\right),\end{multline*}
while on the other hand, if $R$ is extended to be complex linear,
\begin{multline*} 4 R\left(\cdot , \frac{\partial f}{\partial z} \right)\frac{\partial f} {\partial \bar z} = R \left( \cdot, \frac{\partial f} {\partial x} \right)\frac{\partial f} {\partial x} + R \left( \cdot, \frac{\partial f} {\partial y} \right)\frac{\partial f} {\partial y} \\ + \sqrt{-1} \left[ R \left( \cdot , \frac{\partial f}{\partial x}\right) \frac{\partial f}{\partial y} - R \left( \cdot , \frac{\partial f} {\partial y}\right) \frac{\partial f}{\partial x} \right] \\ = R \left( \cdot, \frac{\partial f} {\partial x} \right)\frac{\partial f} {\partial x} + R \left( \cdot, \frac{\partial f} {\partial y} \right)\frac{\partial f} {\partial y} - \sqrt{-1} R \left( \frac{\partial f}{\partial x}, \frac{\partial f}{\partial y}\right),\end{multline*}
with the last step following from the Bianchi symmetry.  We conclude that the Jacobi equation, or the linearization of the Euler-Lagrange equation (\ref{E:harmonic}) at an $\omega$-harmonic map $f$, is
\begin{equation} L_{(f,g)}(Z) = 0, \quad \hbox{where} \quad L_{(f,g)} = - \frac{4}{\lambda ^2} \left[ \frac{D}{\partial z} \circ \frac{D}{\partial \bar z} + R\left(\cdot , \frac{\partial f}{\partial z} \right)\frac{\partial f} {\partial \bar z} \right]. \label{E:Jacobiinholoform}\end{equation}

Integration by parts shows that the Hessian of the $\omega $-energy at an $\omega $-harmonic map $f$ extends to the symmetric complex bilinear form
$$d^2E_\omega (f) : \Gamma ({\bf E}) \times \Gamma ({\bf E}) \longrightarrow {\mathbb C}, \quad \hbox{where} \quad {\bf E} = f^*TM \otimes {\mathbb C},$$
which satisfies the formula
\begin{multline} d^2E_\omega (f)(Z,\overline Z) = 4 \int _\Sigma \left[ \left| \frac{DZ}{\partial \bar z} \right|^2 - \left\langle R \left(Z, \frac{\partial f}{\partial z}\right)\frac{\partial f}{\partial \bar z},\overline Z \right\rangle \right] dx dy \\ = 4 \int _\Sigma \left[ \left| \frac{DZ}{\partial \bar z} \right|^2 - \left\langle {\mathcal R}\left( Z \wedge \frac{\partial f}{\partial z}\right), \overline Z \wedge \frac{\partial f}{\partial \bar z} \right\rangle \right] dx dy, \label{E:secondvariationforharmonic} \end{multline}
$\overline Z$ denoting the complex conjugate of $Z$.  Here $\langle \cdot , \cdot \rangle$ denotes the complex bilinear extension of the Riemannian metric on $M$ and ${\mathcal R} : \Lambda ^2{\bf E} \rightarrow \Lambda^2{\bf E}$ is the complex linear extension of the curvature operator.

It follows directly from these formulae that a holomorphic section of ${\bf L}$ is automatically a Jacobi field (since the curvature term vanishes for sections of ${\bf L}$).  It is easily checked that a section $X$ of $(f^*TM)^\top$ is a Jacobi field if and only if the corresponding section $\tau (X)$ of ${\bf L}$ is also a Jacobi field.  Thus for an $\omega $-harmonic map $f$, the real dimension of the space of Jacobi fields which are sections of $(f^*TM)^\top$ is twice the complex dimension of ${\mathcal O}({\bf L})$, the space of holomorphic sections of ${\bf L}$.

\section{Bumpy metrics  for harmonic surfaces}
\label{S:harmonic}

To illustrate the approach used throughout the article, we first prove a bumpy metric theorem for harmonic maps which satisfy additional hypotheses.

Following \cite{McS}, we say that a harmonic map $f:\Sigma \rightarrow M$ is {\em somewhere injective\/} if there is a point $p \in \Sigma $ such that $f^{-1}(f(p)) = \{ p \}$.  The hypothesis somewhere injective eliminates degenerate harmonic maps, such as torus parametrizations of geodesics, which would need to be treated in a complete theory of harmonic surfaces.

A map which is somewhere injective satisfies the condition $f^{-1}(f(p)) = \{ p \}$, for $p$ lying in an open dense set of points.  To prove this, we can use a theorem of Sampson (\cite{Samp}, Theorem~1) which states that two harmonic maps $f,f' : D \rightarrow M$ from a connected domain $D$ in a Riemann surface which agree on an open set must be identical as maps.  As in Yang-Mills Theory (Theorem 6.38 of \cite{FU}) and in the theory of $J$-holomorphic curves \cite{McS}, the proof is an application of Aronszajn's unique continuation theorem \cite{Ar} to the equation
$$\frac{\partial ^2 u_k}{\partial x_1^2} + \frac{\partial ^2 u_k}{\partial x_2^2} + \sum _{i,j} \Gamma ^k_{ij} \left(\frac{\partial u_i}{\partial x_1} \frac{\partial u_j}{\partial x_1} + \frac{\partial u_i}{\partial x_2} \frac{\partial u_j}{\partial x_2}\right) = 0$$
for harmonic maps, in which we use the abbreviation $u_k$ for $u_k \circ f$.

\vskip .1 in
\noindent
{\bf Proposition~3.1.} {\sl Suppose that $M$ is a compact connected manifold of dimension at least three and $\omega \in {\mathcal T}_g$, the Teichm\"uller space for the Riemann surface of genus $g$.  For a generic choice of Riemannian metric on $M$, every compact somewhere injective $\omega $-harmonic surface $f : \Sigma \rightarrow M$ with no branch points is as nondegenerate as allowed by its connected group $G$ of symmetries.  If $G$ is trivial, it is Morse nondegenerate for $E_\omega $ in the usual sense, while if $G = S^1\times S^1$ or $PSL(2,{\mathbb C})$, it lies on a nondegenerate critical submanifold for $E_\omega$, which is an orbit for the $G$-action.}

\vskip .1 in
\noindent
The techniques for the proof of Proposition~3.1 are similar to those employed in \cite{U1} and \cite{W2}.  We will utilize the spaces
\begin{multline} \hbox{Map}'(\Sigma, M) = \{ f \in \hbox{Map}(\Sigma , M) \hbox{ such that} \\ \hbox{$f$ is a somewhere injective immersion} \},\label{E:defofmapprime}\end{multline}
and
$$\hbox{Met}(M) = \{ \hbox{Riemannian metrics on $M$}\},$$
which we assume have been completed with respect to suitable norms, a Sobolev $L^2_k$ norm for the maps and an $L^2_{\ell}$ norm for the metrics, when $k$ is sufficiently large and $\ell \geq k-1$ is frequently much larger.  When it is necessary to be specific about the value of $k$ and $\ell$, we will write
$$(L^2_k)'(\Sigma ,M), \quad L^2_k(\Sigma ,M) \quad \hbox{and} \quad \hbox{Met}^2_{\ell}(M),$$
instead of
$$\hbox{Map}'(\Sigma ,M), \quad \hbox{Map}(\Sigma ,M) \quad \hbox{and} \quad \hbox{Met}(M).$$

\begin{figure}[tbh]
\centering
\includegraphics[width=.9\textwidth]{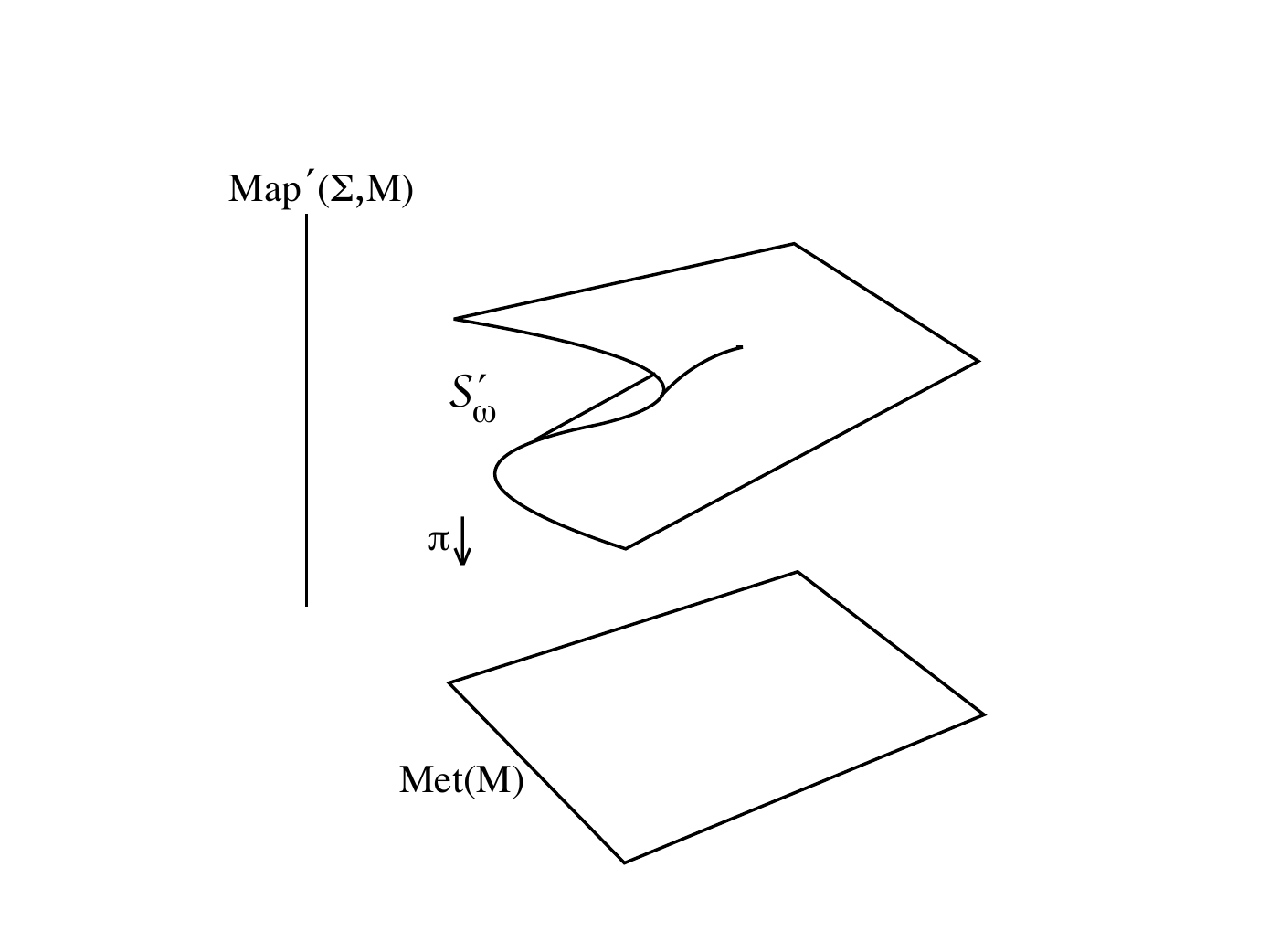}
\caption{The overall strategy of the proof is to show that ${\mathcal S}'_\omega$ is a submanifold and the projection $\pi : {\mathcal S}'_\omega \rightarrow \hbox{Met}(M)$ is Fredholm of Fredholm index zero.}
\label{fig:fredholm4}
\end{figure}

The overall strategy of the proof is show that
\begin{multline} {\mathcal S}'_\omega = \{ (f,g) \in (L^2_k)'(\Sigma ,M) \times \hbox{Met}^2_{\ell }(M) \hbox{ such that} \\
\hbox{$f$ is an $\omega $-harmonic immersion for $g$ } \}, \label{E:somegaprime}\end{multline}
where $\ell \geq k-1$, is a submanifold of the product $\hbox{Map}'(\Sigma ,M) \times \hbox{Met}(M)$, and that the projection $\pi : {\mathcal S} \rightarrow \hbox{Met}(M)$ is a Fredholm map wiith Fredholm index being the dimension of the group $G$ of symmetries.  We then apply the Sard-Smale Theorem.  To carry this out, we consider the parametrized Euler-Lagrange map
\begin{equation} F = F_\omega : (L^2_k)'(\Sigma ,M) \times \hbox{Met}^2_{\ell }(M) \longrightarrow L^2_{k-2}(\Sigma ,TM) \label{E:peulerlagrangemap}\end{equation}
defined in terms of local conformal coordinates $(x^1,x^2)$ on $\Sigma $, by
\begin{equation} F(f,g) = - \frac{1}{\lambda ^2} \left[ \frac{D^g}{\partial x_1} \left(\frac{\partial f}{\partial x_1}\right) + \frac{D^g}{\partial x_2} \left(\frac{\partial f}{\partial x_2}\right) \right],\label{E:before}\end{equation}
where $\lambda ^2 (dx_1^2 + dx_2^2)$ is the canonical Riemannian metric on $\Sigma $ and $D^g$ is the covariant derivative for the Levi-Civita connection on $TM$ determined by the Riemannian metric $g$ on $M$.  Note that $\pi \circ F(f,g) = f$, where
$$\pi : L^2_{k-2}(\Sigma ,TM) \rightarrow L^2_{k-2}(\Sigma ,M)$$
is the obvious projection, a map which is smooth when $k \geq 4$ by the $\omega $-Lemma.  If we can regard
$${\mathcal E} = \{ (f,g, X) \in (L^2_k)'(\Sigma ,M) \times \hbox{Met}^2_{\ell }(M) \times L^2_{k-2}(\Sigma ,TM) : \pi (X) = f \}$$
as (the total space of) a smooth vector bundle over $(L^2_k)'(\Sigma ,M)$, then $F$ defines a smooth section
$$F : (L^2_k)'(\Sigma ,M) \times \hbox{Met}^2_{\ell}(M) \longrightarrow {\mathcal E}.$$
Our goal is to show that $F$ is as transversal as possible to the image ${\mathcal Z}$ of the zero section of ${\mathcal E}$.

The transversality condition is standard when $g \geq 2$ and in that case we can take $\ell = k-1$.  On the other hand, when $\Sigma $ has a positive-dimensional group $G$ of symmetries, there is an advantage to demanding a higher level of regularity of the metrics, because it then follows from regularity theory that
$$F(f,g) = 0 \quad \Rightarrow \quad f \in (L^2_m)'(\Sigma ,M),$$
where $m$ is significantly larger than $k-2$; then (\ref{E:lossofderivatives}) provides as many derivatives as we want for the $G$-action
\begin{equation} (L^2_m)'(\Sigma ,M) \times G \longrightarrow L^2_{k-2}(\Sigma ,M) \label{E:smoothnessofG}\end{equation}
when we choose $\ell$ to be sufficiently large.  Hence for each zero $(f,g)$ of $F$, we have a corresponding family of Jacobi fields
$${\mathcal G}(f,g) \subseteq {\mathcal E}_{(f,g)},$$
the dimension of ${\mathcal G}(f,g)$ being equal to the dimension of $G$.  Straightforward differentiation shows that
$$\int _\Sigma \langle F(f,g) , X \rangle dA = dE_\omega (f,g)(X) = 0, \qquad \hbox{when} \qquad X \in {\mathcal G}(f,g).$$
Thus the transversality condition we seek in this case is
\begin{equation} (\hbox{image of $DF_{(f,g)}$}) + {\mathcal G}(f,g) + T_{(f,g)}{\mathcal Z} = T_{(f,g,0)}{\mathcal E}. \label{E:transversality}\end{equation}
It is then relatively easy to use the $G$-action and the implicit function theorem to show that ${\mathcal S}'_\omega = F^{-1}({\mathcal Z})$ is a submanifold.

\vskip .1 in
\noindent
{\bf Digression~3.2.}  Here are more details regarding one approach to proving that ${\mathcal S}'_\omega$ is a submanifold In the case where $\Sigma = S^2$ or $T^2$, and hence the dimension of $G$ is positive.  The idea is to break the $G$-symmetry.  Suppose first that $\Sigma = S^2$.  If $N$ is a compact codimension two submanifold of $M$ with boundary $\partial N$, we let
\begin{multline}{\mathcal U}(N) = \{ f \in (L^2_k)'(S^2, M) : \hbox{ $f$ does not intersect $\partial N$ and has nonempty }\\ \hbox{ transversal intersection with the interior of $N$ } \}, \label{E:un1} \end{multline}
which is an open subset of $\hbox{Map}'(S^2, M)$.  Given three disjoint compact codimension two submanifolds with boundary, say $Q$, $R$ and $S$, we let
\begin{equation}{\mathcal U}(Q,R,S) = {\mathcal U}(Q) \cap {\mathcal U}(R) \cap {\mathcal U}(S), \label{E:un2} \end{equation}
also an open subset of $\hbox{Map}'(S^2, M)$.  We cover $\hbox{Map}'(S^2, M)$ with a countable collection of sets ${\mathcal U}(Q_j,R_j,S_j)$ defined by a sequence $j \mapsto (Q_j,R_j,S_j)$ of triples of such codimension two submanifolds.  We choose three points $q$, $r$ and $s$ in $S^2$  and let
\begin{equation}{\mathcal F}_j(S^2 ,M) = \{ f \in {\mathcal U}(Q_j,R_j,S_j) : f(q) \in Q_j, f(r) \in R_j, f(s) \in S_j  \},\label{E:defoffiS2}\end{equation}
noting that ${\mathcal F}_j(S^2 ,M)$ meets each $PSL(2,{\mathbb C})$-orbit in ${\mathcal U}(Q_j,R_j,S_j)$ in a finite number of points.  It follows from the Sobolev imbedding theorem and smoothness of the evaluation map on the space of $C^{k-2}$ maps (for example by Proposition~2.4.17 of \cite{AMR}) that the evaluation map
$$\hbox{ev} : L^2_k(S^2, M) \times S^2 \longrightarrow M, \qquad \hbox{defined by} \qquad \hbox{ev}(f,p) = f(p),$$
is $C^{k-2}$.  Thus when $k$ is large, we can regard ${\mathcal F}_j(S^2,M)$ as a submanifold of $\hbox{Map}(S^2, M)$ of codimension six with tangent space
\begin{multline*}T_f{\mathcal F}_j(S^2,M) = \{ \hbox{sections $X$ of $f^*TM$} : \\ X(q) \in T_{f(q)}Q_j, X(r) \in T_{f(r)}R_j, X(s) \in T_{f(s)}S_j \}.\end{multline*}
In the case where $\Sigma $ is the torus, we need fix only one point to break the symmetry.  We define ${\mathcal U}(N)$ by (\ref{E:un1}) with $S^2$ replaced by $T^2$ and choose a sequence $j \mapsto N_j$ of smooth compact codimension two submanifolds of $M$ such that ${\mathcal U}(N_j)$ cover $\hbox{Map}(T^2, M)$.  We choose a base point $q \in T^2$ and let
\begin{equation} {\mathcal F}_j(T^2 ,M)= \{ f \in {\mathcal U}(N_j) : f(q) \in N_j  \}. \label{E:defoffiT2}\end{equation}
If we restrict $F$ to either ${\mathcal F}_j(S^2,M)$ or ${\mathcal F}_j(T^2,M)$, (\ref{E:transversality}) gets replaced by the simpler condition
$$(\hbox{image of $DF_{(f,g)}$}) + T_{(f,g)}{\mathcal Z} = T_{(f,g,0)}{\mathcal E},$$
which makes it easier to directly apply the implicit function theorem.  It is then straightforward to use the $G$-action to construct local coordinates about any element in ${\mathcal S}'_\omega$, and to verify (\ref{E:transversality}).  Although not absolutely needed for maps with no branch points, this idea of breaking the $G$-symmetry will turn out to be useful in later steps of the argument for the Main Theorem.

\vskip .1 in
\noindent
Note that at the image of a zero $(f,g)$ for $F$, the tangent space to the total space ${\mathcal E}$ can be divided into a direct sum
$$T_{(f,g,0)}{\mathcal E} = H \oplus V,$$
where $H$ is the horizontal space (tangent to the space of zero sections) and $V$ is the vertical space (tangent to the fiber).  Let $\pi _V$ denote the projection onto $V$ along $H$.  Our strategy for proving transversality is to show that ${\mathcal V} = {\mathcal G}$, where
\begin{multline}{\mathcal V} = \{\hbox{ elements $X$ of ${\mathcal E}_{(f,g)}$} : \hbox{$X$ is $L^2$-perpendicular to }\\
\hbox{ the image of $\pi _V \circ DF_{(f,g)}$} \}. \label{E:mathcalV}\end{multline}
It follows from (\ref{E:jacob}) that differentiation with respect to the map $f$ yields
\begin{equation} \pi _V \circ (D_1F)_{(f,g)} (X) = L(X), \label{E:jacobioperator1}\end{equation}
where $L = L_{(f,g)}$ is the Jacobi operator at the critical point $(f,g)$.  On the other hand, differentiation with respect to the second variable $g \in \hbox{Met}(M)$,
$$\pi_V \circ (D_2F)_{(f,g)} (\dot g),$$
gives the effect of deformations in the metric on $M$.

\vskip .1 in
\noindent
We note that the following lemma does not require any hypothesis on the branch locus.

\vskip .1 in
\noindent
{\bf Lemma~3.3.} {\sl If $f$ is a somewhere injective critical point for $E_\omega $, the only sections of $f^*TM$ which are $L^2$-perpendicular to the image of every variation in the metric are real and imaginary parts of holomorphic sections of ${\bf L}$.}

\vskip .1 in
\noindent
Proof:  To prove this, we need to calculate
$$\pi_V \circ (D_2F)_{(f,g)} (\dot g),$$
at an element $(f,g) \in {\mathcal S}$ for various choices of the deformation $\dot g \in T_g\hbox{Met}(M)$.  Note that since $F(f,g)$ is the zero section,
\begin{equation} \langle \pi_V \circ (D_2F), X \rangle = D_2 \left\langle F(f,g), X \right\rangle, \label{E:secondderivwithmetric}\end{equation}
where $\langle \cdot , \cdot \rangle$ is another name for the metric $g$ on the ambient space $M$. 

Suppose that $f$ is an $\omega $-harmonic map for $g$.  The set $W \subset \Sigma $, consisting of points $p$ such that $p$ is not a branch point and $f^{-1}(f(p)) = p$, is open and dense by the somewhere injective hypothesis.   Choose an open neighborhood $U$ of a point $p \in W$ such that $U \subset W$ and $f$ imbeds $U$ onto $f(\Sigma )\cap V$ for some open set $V \subset M$.  Arrange, moreover, that $V$ is the domain of local coordinates $(u_1 ,\ldots , u_n)$ such that $u_i(f(p)) = 0$ and
\begin{enumerate}
\item $f(U)$ is described by the equations $u_3 = \cdots = u_n = 0$,
\item $u_a \circ f = x_a$ on $f(U)$, for $a = 1,2$, where $x_1 + i x_2$ is a conformal parameter on $U$, and
\item the Riemannian metric $g$ on the ambient space takes the form $\sum g_{ij} du_idu_j$, such that when restricted to $f(\Sigma ) \cap V$, $g_{ir} = \delta _{ir}$, for $1 \leq i \leq n$ and $3 \leq r \leq n$.
\end{enumerate}
Such coordinates can be constructed using the exponential map restricted to the normal bundle of the surface $f(\Sigma )\cap V$ in $M$.  Note that in the special case in which $f$ is not only $\omega $-harmonic, but also conformal,
\begin{equation} g_{ab} = \sigma ^2 \eta_{ab}, \quad \eta_{ab} = \lambda ^2 \delta_{ab}, \quad \hbox{for $1 \leq a,b \leq 2$,} \label{E:conformalfactor}\end{equation}
where $\sigma ^2$ is a positive function and $(\eta_{ab})$ is a fixed metric on $\Sigma $ within the conformal equivalence class, usually normalized to have constant curvature and total area one.

Recall that in terms of local coordinates, the equation for harmonic maps can be written as
$$\frac{1}{\lambda ^2} \left(\frac{\partial ^2 u_k}{\partial x_1^2} + \frac{\partial ^2 u_k}{\partial x_2^2}\right) + \sum _{i,j} \frac{1}{\lambda ^2} \Gamma ^k_{ij} \left(\frac{\partial u_i}{\partial x_1} \frac{\partial u_j}{\partial x_1} + \frac{\partial u_i}{\partial x_2} \frac{\partial u_j}{\partial x_2}\right) = 0,$$
the canonical constant curvature Riemannian metric on $\Sigma $ being $\lambda ^2 (dx_1^2 + dx_2^2)$.  The expression on the left-hand side of this equation is often called the {\em tension\/} of the map $f$, a harmonic map being one which has zero tension.  In terms of the special local coordinates, the equation for harmonic maps yields
$$\Gamma ^k_{11} + \Gamma ^k_{22} = 0.$$

A perturbation in the metric $\dot g \in T_g\hbox{Met}(M)$ with compact support in $V$ can be written in the form $\dot g = \sum \dot g_{ij} dx_idx_j$, where the $\dot g_{ij}$'s are smooth functions on $V$.  Under such a perturbation, the only part of the tension that changes is the Christoffel symbol
$$\Gamma ^k_{ij} = \sum g^{kl} \Gamma_{l,ij}, \qquad \hbox{where} \qquad \Gamma_{l,ij} = \frac{1}{2}\left(\frac{\partial g_{il}}{\partial u_j} + \frac{\partial g_{jl}}{\partial u_i} - \frac{\partial g_{ij}}{\partial u_l} \right).$$
If $\dot \Gamma _{k,ij}$ denotes the derivative of $\Gamma _{k,ij}$ in the direction of the perturbation,
$$\dot \Gamma _{k,ij} = \frac{1}{2}\left( \frac{\partial \dot g_{ik}}{\partial u_j} +  \frac{\partial \dot g_{jk}}{\partial u_i} - \frac{\partial \dot g_{ij}}{\partial u_k} \right).$$
From this formula, we can calculate that
\begin{multline*} D_2 \left\langle F(f,g), X \right\rangle (\dot g) = - \frac{1}{\lambda ^2}\sum _{i,j,k=1}^n \dot \Gamma _{k,ij} \left(\frac{\partial u_i} {\partial x_1} \frac{\partial u_j}{\partial x_1} + \frac{\partial u_i}{\partial x_2} \frac{\partial u_j}{\partial x_2}\right) f^k\\
= - \frac{1}{\lambda ^2} \sum _{k=1}^n (\dot \Gamma _{k,11} + \dot \Gamma _{k,22}) f^k, \qquad \hbox{when} \qquad X = \sum _{i=1}^n f^i \frac{\partial }{\partial u_i}\end{multline*}
is a fixed vector field along $f$.  Thus the definition (\ref{E:mathcalV}) of ${\mathcal V}$ implies that
$$X \in {\mathcal V} \quad \Rightarrow \quad \int _\Sigma \sum _{k=1}^n \sum _{a = 1}^2 f^k \dot \Gamma _{k,aa} dx_1dx_2 = 0,$$
for all perturbations in the metric.

If we set
$$ \dot g_{11}(u_1, \ldots ,u_n) = \dot g_{22}(u_1, \ldots ,u_n) = \sum_{r=3}^n u_r \phi_r (u_1,u_2),$$
and let $\dot g_{ij} = 0$ for the other choices of indices $i,j$, a straightforward calculation shows that the fiber projection of the partial derivative of $F$ with respect to $g$ is given by the expression
$$\pi _V \circ (D_2F)_{(f,g)} (\dot g) = \frac{1}{2\lambda ^2} \sum _{r=3}^n\sum _{a=1}^2 \frac{\partial \dot g_{aa}}{\partial u_r} \frac{\partial}{\partial u_r} = \frac{1}{\lambda ^2} \sum _{r=3}^n \frac{\partial \dot g_{11}}{\partial u_r} \frac{\partial}{\partial u_r}.$$ 
Thus we see that any vector field of the form
$$\sum _{r=3}^n \phi_r (u_1, u_2) \frac{\partial}{\partial u_r},$$
where $\phi _r$ is a smooth function on $U$ with compact support, lies in the image of $\pi _V \circ (D_2F)_{(f,g)}$, and hence elements of ${\mathcal V} \otimes {\mathbb C}$ must restrict to sections of ${\bf L} \oplus {\bf \bar L}$ over $U$.  Since a dense open subset of points of $\Sigma $ can be covered by sets of the form $U$, we see that elements of ${\mathcal V}$ must be sections of ${\bf L} \oplus {\bf \bar L}$ over all of $\Sigma$.

To determine the conditions satisfied by the tangential components of elements of ${\mathcal V}$, we set
\begin{equation}\dot g_{11} = - \dot g_{22} = M , \quad \dot g_{12} = N , \quad \dot g_{ij} = 0 \quad \hbox{for other choices of indices $i,j$,}\label{E:tangentialvar}\end{equation}
where $M$ and $N$ have compact support within $V$, and a short calculation shows that
$$\dot \Gamma _{1,11} + \dot \Gamma _{1,22} = \frac{1}{2}\frac{\partial \dot g_{11}}{\partial u_1} + \frac{\partial \dot g_{12}}{\partial u_2} - \frac{1}{2}\frac{\partial \dot g_{22}}{\partial u_1} = \frac{\partial M}{\partial u_1} + \frac{\partial N}{\partial u_2},$$
$$\dot \Gamma _{2,11} + \dot \Gamma _{2,22} = \frac{\partial \dot g_{12}}{\partial u_1} - \frac{1}{2}\frac{\partial \dot g_{11}}{\partial u_2}  + \frac{1}{2}\frac{\partial \dot g_{22}}{\partial u_2} = \frac{\partial N}{\partial u_1} - \frac{\partial M}{\partial u_2}.$$
Hence in this case,
\begin{multline} X = \sum _{a=1}^2 f^a \frac{\partial}{\partial u_a} \Rightarrow \\ D_2\langle F(f,g),X \rangle (\dot g) = - \frac{1}{\lambda ^2} \left[\left(\frac{\partial M}{\partial u_1} + \frac{\partial N}{\partial u_2}\right) f^1+ \left(\frac{\partial N}{\partial u_1} - \frac{\partial M} {\partial u_2}\right) f^2 \right] \label{E:derinmetricdir} \end{multline}
and
\begin{multline} \int _\Sigma \left\langle \pi _V \circ (D_2F)_{(f,g)} (\dot g), X \right\rangle dA \\ = - \int _\Sigma \left[\left(\frac{\partial M}{\partial u_1} + \frac{\partial N}{\partial u_2}\right) f^1 + \left(\frac{\partial N} {\partial u_1} - \frac{\partial M} {\partial u_2}\right) f^2 \right] dx_1dx_2. \label{E:varinmetric1} \end{multline}
In the current situation $u_a = x_a$, for $1 \leq a \leq 2$, so we conclude that $f^1 (\partial/\partial x_1) + f^2 (\partial/\partial x_2)$ is perpendicular to the range of the map $\pi _V \circ (D_2F)_{(f,g)}$ if and only if
$$ \int_\Sigma \left[ f^1 \left(\frac{\partial M}{\partial x_1} + \frac{\partial N}{\partial x_2}\right) + f^2 \left(\frac{\partial N}{\partial x_1} - \frac{\partial M}{\partial x_2}\right) \right] dx_1dx_2 = 0, $$
or equivalently, if and only if
$$\int_\Sigma \left[ M \left(\frac{\partial f^1}{\partial x_1} - \frac{\partial f^2}{\partial x_2}\right) + N \left(\frac{\partial f^1}{\partial x_2} + \frac{\partial f^2}{\partial x_1}\right) \right] dx_1dx_2 = 0,$$
for all choices of functions $M$ and $N$, which one verifies is exactly the condition that $f^1 (\partial/\partial x_1) +f^2 (\partial/\partial x_2)$ be the real or imaginary part of a holomorphic section of ${\bf L}$.  It follows that the $L^2$ orthogonal complement of the range of $\pi _V \circ (D_2F)_{(f,g)}$ is contained in the space of real and imaginary parts of holomorphic sections of ${\bf L}$, and Lemma~3.3 is proven.

\vskip .1 in
\noindent
{\bf Remark~3.4.}  The tangential variation $\dot g$ used in the above proof can be regarded as defined by a quadratic differential $\phi dz^2 = (1/2)(M - i N) dz^2$.  If
$$Z = \tau (X) =  (f^1 + i f^2) \left( \frac{\partial }{\partial x_1} - i \frac{\partial }{\partial x_2} \right) = 2(f^1 + i f^2) \frac{\partial }{\partial z},$$
so that the real part of $Z$ is the vector field $X$ considered above, then (\ref{E:derinmetricdir}) yields
\begin{multline} D_2 \langle F(f,g), Z \rangle (\dot g) = - \frac{1}{\lambda ^2} \left( (f^1 + i f^2)\frac{\partial}{\partial \bar z} (M - iN) \right) = - \frac{2}{\lambda ^2} \frac{\partial \phi}{\partial \bar z}dz(Z), \label{E:holocomparison} \end{multline}
where we have extended $\langle F(f,g), \cdot \rangle$ to be complex linear.  An integration by parts yields
\begin{multline} \int _\Sigma \langle \pi _V \circ (D_2F)_{(f,g)} (\dot g),Z \rangle dA = - \int _\Sigma \frac{2}{\lambda ^2} \left( \frac{\partial \phi}{\partial \bar z }dz (Z)\right) dA \\
=  \int _\Sigma \frac{2}{\lambda ^2} \phi dz \left( \frac{DZ}{\partial \bar z} \right) dx_1dx_2.  \label{E:comparison} \end{multline}
If $Z$ is a section of ${\bf L}$. say $Z = \zeta (\partial f/\partial z)$, where $\zeta $ is a complex-valued function and $f$ is conformal, say $f^*\langle \cdot , \cdot \rangle = \sigma ^2(dx_1^2 + dx_2^2)$, we can write
$$\left\langle \frac{\partial f}{\partial \bar z} , \frac{\partial f}{\partial z}\right\rangle = \frac{\sigma ^2}{2} \qquad \Rightarrow \qquad dz \left( \frac{DZ}{\partial \bar z} \right) = \frac{ \partial \zeta}{\partial \bar z} = \frac{2}{\sigma ^2} \left\langle \frac{\partial f}{\partial \bar z} , \frac{DZ}{\partial \bar z} \right\rangle,$$
and comparison with (\ref{E:comparison}) yields a formula needed later:
\begin{equation} \int _\Sigma \langle \pi _V \circ (D_2F)_{(f,g)} (\dot g),Z \rangle dA = \int _\Sigma \frac{4}{\lambda ^2} \left\langle \frac{DZ}{\partial \bar z} , \frac{\phi}{\sigma ^2} \frac{\partial f}{\partial \bar z} \right\rangle dA.\label{E:desiredformula} \end{equation}

\vskip .1 in
\noindent
We return now to the proof of Proposition~3.1.  Since $f$ is an immersion, the real and imaginary parts of holomorphic sections of ${\bf L}$ are elements of ${\mathcal G}$.  It is now easy to see that
$$\pi _V\circ DF_{(f,g)} : T_f(L^2_k)'(\Sigma ,M) \times T_g\hbox{Met}^2_{\ell }(M) \longrightarrow L^2_{k-2}(\Sigma ,f^*TM)$$
maps onto a complement to ${\mathcal G}$, when $f$ is $\omega $-harmonic with respect to $g$ and $k$ is large.  Indeed, the Jacobi operator
$$\pi _V \circ D_1F_{(f,g)} : T_f(L^2_k)'(\Sigma ,M) \longrightarrow L^2_{k-2}(\Sigma ,f^*TM)$$
maps onto a complement to the Jacobi fields for $f$, a closed subspace of finite codimension and Lemma~3.1 implies that all Jacobi fields are covered by $\pi _V \circ D_2F_{(f,g)}$ except for the elements of ${\mathcal G}$.  It therefore follows from the implicit function theorem that the space ${\mathcal S}'_\omega$ defined by (\ref{E:somegaprime}) is indeed a submanifold as claimed.

\vskip .1 in
\noindent
{\bf Lemma~3.5.} {\sl The projection on the second factor $\pi : {\mathcal S}'_\omega \rightarrow \hbox{Met}(M)$ is a Fredholm map of Fredholm index $\dim G$, where $G$ is the group of symmetries.}

\vskip .1 in
\noindent
Proof:  A straightforward calculation yields the tangent space
\begin{multline} T_{(f,g)}{\mathcal S}'_\omega = \{ (X, \dot g) \in T_f(L^2_k)' \oplus T_g\hbox{Met}^2_{\ell}(M) : \\ \hbox{$L(X) + \pi _V \circ (D_2F)_{(f,g)} (\dot g) = 0$} \}, \label{eq:expfortang1} \end{multline}
where $L = L_{(f,g)}$ is the Jacobi operator at $(f,g)$, a Fredholm map of Fredholm index zero, with
$$\hbox{Kernel of $L$} = {\mathcal G} \oplus {\mathcal J} \cong (\hbox{complement to Range of $L$}),$$
where ${\mathcal J}$ is a finite-dimensional space of nontangential Jacobi fields.  For the differential of the projection,
$$d\pi_{(f,g)}: T_{(f,g)} {\mathcal S} \rightarrow T_g\hbox{Met}^2_{\ell}(M),$$
we find that
$$\hbox{Ker}(d\pi_{(f,g)}) = \{ (X, 0) \in T_f(L^2_k)'(\Sigma ,M) \oplus T_g\hbox{Met}^2_{\ell}(M) : \\ \hbox{$L(X) = 0$}\} \cong \hbox{Ker}(L),$$
while
\begin{multline*}\hbox{Range}(d\pi_{(f,g)}) = \{ \dot g \in T_g\hbox{Met}^2_{\ell}(M) : \pi _V \circ (D_2F)_{(f,g)} (\dot g) + L(X) = 0, \\ \hbox{ for some $X \in T_f(L^2_k)'(\Sigma ,M)$} \}.\end{multline*}
Thus the range of $d\pi_{(f ,g)}$ consists of the elements
\begin{multline*}\dot g \in T_g\hbox{Met}^2_{\ell}(M) \quad \hbox{such that} \quad T(\dot g) \in \hbox{Range}(L), \\ \hbox{where} \quad T(\dot g) = \pi _V \circ (D_2F)(\dot g),\end{multline*}
the preimage of a closed space of finite codimension, hence a closed subspace of finite codimension itself.  Thus $d\pi_{(f,g)}$ is indeed a Fredholm map, and the dimension of the cokernel of $d\pi_{(f,g)}$ is no larger than the dimension of the cokernel of $L$.  But Lemma~3.3 shows that $dF_{(f ,g)}$ maps surjectively onto a complement to ${\mathcal G}$, so the linear space ${\mathcal J}$ must be covered by $T$ and the cokernel of $d\pi_{(f ,g)}$ has the same dimension as ${\mathcal J}$.  Hence $d\pi_{(f,g)}$ is a Fredholm map of Fredholm index $\dim {\mathcal G}$, and Lemma~3.5 is proven.

\vskip .1 in
\noindent
Finally, we use the Sard-Smale Theorem \cite{Sm2} to conclude that a generic subset of $\hbox{Met}(M)_{\ell}$ consists of regular values for the projection $\pi $ on the second factor.  Any such metric $g_0$ has the property that all somewhere injective $\omega $-harmonic surfaces for $g_0$ will have no nonvanishing Jacobi fields, except for the infinitesimal symmetries ${\mathcal G}$ generated by the $G$-action.  Such $\omega $-harmonic surfaces must lie on nondegenerate critical submanifolds in $(L^2_k)'(\Sigma ,M)$ which also isolated $G$-orbits of critical points.  This finishes the proof of Proposition~3.1.

\vskip .1 in
\noindent
Since prime $\omega $-harmonic immersions are automatically somewhat injective, we obtain the following consequence of the proof of Proposition 3.1:

\vskip .1 in
\noindent
{\bf Corollary~3.6.} {\sl If $M$ is a compact connected manifold of dimension at least three and $\omega \in {\mathcal T}_g$, then
\begin{multline*} {\mathcal P}^\emptyset _\omega = \{(f,g) \in L^2_k(\Sigma ,M) \times \hbox{Met}_{\ell}(M) : \hbox { $f$ is a prime }\\ \hbox{ $\omega $-harmonic immersion for $g$} \} \end{multline*}
is a smooth submanifold when $\ell $ is sufficiently large, and
$$ \pi : {\mathcal P}^\emptyset _\omega \longrightarrow \hbox{Met}_{\ell }(M) $$
is a Fredholm map of Fredholm index $d_\Sigma $.}

\vskip .1 in
\noindent
The hypotheses of Proposition 3.1 can be weakened slightly:

\vskip .1 in
\noindent
{\bf Proposition~3.7.} {\sl If $M$ is a compact connected manifold of dimension at least three and $\omega \in {\mathcal T}_g$, then or a generic choice of Riemannian metric on $M$, every compact somewhere injective $\omega $-harmonic surface $f : \Sigma \rightarrow M$ such that the space ${\mathcal O}({\bf L})$ has the same real dimension as $G$ lies on a nondegenerate critical submanifold of dimension $d_\Sigma $, where $d_\Sigma $ is the real dimension of the group $G$ of symmetries.}

\vskip .1 in
\noindent
The proof is identical to the proof of Proposition~3.1, except that we replace ${\mathcal S}'_\omega$ with
\begin{multline*} \tilde {\mathcal S}_\omega = \{ (f,g) \in (L^2_k)_0(\Sigma ,M) \times \hbox{Met}(M)_{k-1} \hbox{ such that} \\
\hbox{$f$ is $\omega $-harmonic for $g$ and $\dim {\mathcal O}({\bf L}) = \dim G$ }. \}\end{multline*}
The assumption on ${\mathcal O}({\bf L})$) is used to insure that the space of tangential Jacobi fields (isomorphic to the space of holomorphic sections of ${\bf L}$) is exactly the space generated by the action of $G$:
\begin{enumerate}
\item If $\Sigma $ is a two-sphere, this is the same as the condition that the branch locus be empty.
\item If $\Sigma$ is a torus, $\hbox{dim}_{\mathbb C}{\mathcal O}({\bf L}) = 1$, the minimal value exactly when $f$ has total branching order zero or one.
\item If $\Sigma $ has genus at least two, ${\mathcal O}({\bf L}) = 0$ when the total branching order is $< 2g-2$, and if ${\bf L}$ is nontrivial as a holomorphic bundle when the total branching order is $2g-2$.
\end{enumerate}

\vskip .1 in
\noindent
There is also a corresponding weakening of the hypotheses of Corollary~3.6, which the reader can easily state.

\section{The two-variable energy}
\label{S:two-variableenergy}

Throughout this section, we assume that $\Sigma $ is an oriented surface of genus at least one.  We should regard the proof of Proposition~3.1 as really a proof schema which we need to generalize to the cases where the $\omega $-energy is replaced by the two-variable energy
\begin{equation}E : \hbox{Map}(\Sigma ,M) \times {\mathcal T}_g \longrightarrow {\mathbb R},\label{E:twovariableE}\end{equation}
whose critical points are exactly the conformal harmonic maps or minimal surfaces, and further generalized to the case where the branch locus is nontrivial.  Here ${\mathcal T}_g$ is the Teichm\"uller space for compact Riemann surfaces of a fixed genus $g$.  Actually this map is induced by a more basic energy function:  Let $\hbox{Met}(\Sigma )$ denote the space of smooth Riemannian metrics on the compact oriented surface $\Sigma $.  We can write a typical element $\eta \in \hbox{Met}(\Sigma )$ in terms of fixed coordinates $(x_1,x_2)$ on $\Sigma $ as
$$\eta = \sum _{a,b = 1}^2 \eta _{ab} dx_a dx_b, \quad \hbox{and let} \quad (\eta^{ab}) = (\eta_{ab})^{-1},$$
and define the energy function
$$E :  \hbox{Map}(\Sigma ,M) \times \hbox{Met}(\Sigma ) \longrightarrow {\mathbb R}$$
by means of the coordinate formula
\begin{equation} E(f,\eta ) = \frac{1}{2}\int _{\Sigma} \sum _{a,b} \eta ^{ab} \sqrt{\det(\eta _{ab})}  \left\langle\frac{\partial f}{\partial x_a} , \frac{\partial f}{\partial x_b} \right\rangle dx_1 dx_2, \label{E:energyrange}\end{equation}
where
$$\left\langle\frac{\partial f}{\partial x_a} , \frac{\partial f}{\partial x_b} \right\rangle = g\left(\frac{\partial f}{\partial x_a} , \frac{\partial f}{\partial x_b} \right),$$
in terms of the Riemannian metric on $M$.  One verifies directly that the energy $E$ does not change when the metric $\eta $ is multiplied by a conformal factor, that is $E$ is invariant under Weyl rescalings, $(\eta_{ab}) \mapsto (\lambda ^2 \eta_{ab})$ where $\lambda ^2$ is a smooth positive function, the integrand being independent of the choice of local coordinates.  Moreover, any Riemannian metric on $\Sigma $ is equivalent via Weyl rescaling to a metric in $\hbox{Met}_0(M)$, the space of Riemannian metrics on $\Sigma $ that have constant curvature and total area one, so the energy descends immediately to
$$E :  \hbox{Map}(\Sigma ,M) \times \hbox{Met}_0(\Sigma ) \longrightarrow {\mathbb R}.$$ 
Dividing by the action of the diffeomorphism group then yields (\ref{E:twovariableE}).

\subsection{Action of the diffeomorphism group}
\label{S:diffeogroup}

The group $\hbox{Diff}^+(\Sigma )$ of orientation-preserving diffeomorphisms of $\Sigma $ acts on $\hbox{Met}_0(\Sigma )$ by
$$(\eta ,\varphi ) \in \hbox{Met}_0(\Sigma ) \times \hbox{Diff}^+(\Sigma ) \quad \mapsto \quad \varphi ^* \eta \in \hbox{Met}_0(\Sigma ),$$
as does the normal subgroup $\hbox{Diff}_0(\Sigma )$ of diffeomorphisms homotopic to the identity.  The quotients
$${\mathcal M} = \frac{ \hbox{Met}_0(\Sigma )}{\hbox{Diff}^+(\Sigma )} \quad \hbox{and} \quad {\mathcal T} = \frac{ \hbox{Met}_0(\Sigma )}{\hbox{Diff}_0(\Sigma )}$$
are known as the {\em moduli space\/} and {\em Teichm\"uller space\/}, respectively, the quotient $\Gamma = \hbox{Diff}^+(\Sigma )/ \hbox{Diff}_0(\Sigma )$ being the mapping class group mentioned in the Introduction.

The group $\hbox{Diff}^+(\Sigma )$ also acts on $\hbox{Map}(\Sigma ,M)$ by
$$(f ,\varphi ) \in \hbox{Map}(\Sigma ,M) \times \hbox{Diff}^+(\Sigma ) \quad \mapsto \quad f \circ \varphi ,$$
with the energy $E$ being invariant under the corresponding product action of $ \hbox{Diff}^+(\Sigma )$ on $\hbox{Map}(\Sigma ,M) \times \hbox{Met}_0(\Sigma )$.  The function $E$ thus descends to a map on the quotient
$${\mathcal M}(\Sigma ,M) = \frac{\hbox{Map}(\Sigma ,M) \times\hbox{Met}_0(\Sigma )}{\hbox{Diff}^+(\Sigma )}.$$
As pointed out in the Introduction, this space is the most natural domain for developing a Morse theory for the energy $E$.  However, it is simpler to regard the energy as a function on the space
\begin{equation}\widetilde {\mathcal M}(\Sigma ,M) = \frac{\hbox{Map}(\Sigma ,M) \times\hbox{Met}_0(\Sigma )}{\hbox{Diff}_0(\Sigma )}.\label{E:msigmam}\end{equation}
In this case, the projection on the second factor
$$\pi : \hbox{Map}(\Sigma ,M) \times \hbox{Met}_0(\Sigma ) \rightarrow \hbox{Met}_0(\Sigma )$$
descends to a fiber bundle projection from $\widetilde {\mathcal M}(\Sigma ,M)$ to ${\mathcal T}$, the latter space being diffeomorphic to a ball by a well-known theorem from Teichm\"uller theory.  An explicit diffeomorphism is given in Chapter 2, \S 3 of \cite{SY2}, and using it, we can construct a trivialization of the bundle; see also \cite{EE}.  Thus we can also regard the energy $E$ as a function on the product space $\hbox{Map}(\Sigma ,M) \times {\mathcal T}$, in agreement with (\ref{E:twovariableE}), and if $[f, \eta ]$ denotes the equivalence class of $(f,\eta ) \in \hbox{Map}(\Sigma ,M) \times \hbox{Met}_0(\Sigma )$, we can identify
$$[f,\eta] \in \widetilde {\mathcal M}(\Sigma ,M) \qquad \hbox{with} \qquad (f,\omega ) \in \hbox{Map}(\Sigma ,M) \times {\mathcal T},$$
$\omega $ being the conformal class of the metric $\eta $.  

To calculate the first derivative of $E$ on $\widetilde {\mathcal M}(\Sigma ,M)$, we first calculate the derivative of
$$E : \hbox{Map}(\Sigma ,M) \times \hbox{Met}(\Sigma ) \rightarrow {\mathbb R},$$
then restrict to the tangent space to $\hbox{Map}(\Sigma ,M) \times \hbox{Met}_0(\Sigma )$, and finally restrict to the tangent space to a slice for the action of $\hbox{Diff}_0(\Sigma )$ which is transversal to the orbits.  The slice which passes through a representative
$$(f, \eta) \in \hbox{Map}(\Sigma ,M) \times\hbox{Met}_0(\Sigma ) \quad \hbox{for} \quad [f,\eta] \in \widetilde {\mathcal M}(\Sigma ,M)$$
allows us to make the identification
\begin{equation} T_{[f,\eta ]}\widetilde {\mathcal M}(\Sigma ,M) \cong T_f\hbox{Map}(\Sigma ,M) \oplus T_\omega {\mathcal T}. \label{E:identification} \end{equation}
In the metric direction, the slice is chosen to be tangent to the one-parameter families of metrics of the form $\eta_{ij} (t) = \eta _{ij}(0) + t \dot \eta _{ij}$, with $\eta _{ij}(0) = \lambda ^2 \delta _{ij}$ for some conformal factor $\lambda ^2$, the $(\dot \eta _{ij})$ being required to be $L^2$-perpendicular to the orbits of $\hbox{Diff}_0(\Sigma )$ at $(f,\eta )$.  When the genus is at least two, the argument in \cite{SY2} shows that for this choice of slice, $\dot \eta $ satisfies the zero-trace condition $\dot \eta _{11} + \dot \eta _{22} = 0$, and the condition that
\begin{equation} \phi dz^2 = \frac{1}{2}(\dot \eta _{11} - i \dot \eta _{12}) dz^2 \label{E:defofphi0} \end{equation}
be a holomorphic quadratic differential, the factor $(1/2)$ being included so that
$$\phi dz^2 + \bar \phi d \bar z^2 = \sum _{a,b = 1}^2 \dot \eta _{ab} dx_a dx_b.$$
If the genus is one, our normalization of the area forces the trace-free condition, and an argument similar to that in \cite{SY2} shows that $\phi dz^2$ is once again holomorphic.

Indeed, following \cite{DH}, we can write
\begin{equation} T_\eta (\hbox{Met}(\Sigma )) \otimes {\mathbb C} = \Gamma ({\bf K}^{(2,0)}) \oplus \Gamma ({\bf K}^{(1,1)}) \oplus \Gamma ({\bf K}^{(0,2)}), \label{E:decomposition1} \end{equation}
where ${\bf K}$ is the canonical bundle of $\Sigma $, ${\bf \bar K}$ its conjugate, and ${\bf K}^{(p,q)} = {\bf K}^p \otimes{\bf \bar K}^q$.  Sections of ${\bf K}^{(1,1)}$ are exactly the metric deformations which correspond to Weyl rescalings, while
\begin{equation} T_\eta (\hbox{Met}_0(\Sigma )) \otimes {\mathbb C} = \Gamma ({\bf K}^{(2,0)}) \oplus \Gamma ({\bf K}^{(0,2)}), \label{E:decomposition1-0} \end{equation}
where the section of $\Gamma ({\bf K}^{(2,0)})$ is given by (\ref{E:defofphi0}) and the section of $\Gamma ({\bf K}^{(0,2)})$ is its conjugate.  Covariant differentiation with respect to the Levi-Civita connection on $\Sigma $ yields operators
$$D' : \Gamma ({\bf K}^{(p,q)}) \rightarrow \Gamma ({\bf K}^{(p+1,q)}), \quad D'' : \Gamma ({\bf K}^{(p,q)}) \rightarrow \Gamma ({\bf K}^{(p,q+1)}),$$
including an operator
$$D' : \Gamma ({\bf K}^{(1,0)}) \rightarrow \Gamma ({\bf K}^{(2,0)})$$
which represents the infinitesimal action of vector fields on quadratic differentials induced by the action of $\hbox{Diff}_0(\Sigma )$ on $\hbox{Met}_0(\Sigma )$.  The real inner product on $T\Sigma $ extends to a complex bilinear inner product $\langle \cdot , \cdot \rangle$ on the sum of all ${\bf K}^{(p,q)}$'s, which pairs ${\bf K}^{(p,q)}$ with ${\bf K}^{(q,p)}$, or to a Hermitian inner product with respect to which all the ${\bf K}^{(p,q)}$'s are orthogonal, and the real structure on the complexification determines a conjugate linear conjugation map $C : {\bf K}^{(p,q)} \rightarrow {\bf K}^{(q,p)}$.

With these preparations out of the way, we seek a formula for
$$dE([f,\eta]) : T_{[f,\eta ]}\widetilde {\mathcal M}(\Sigma ,M) \longrightarrow {\mathbb R},$$
which is called the first variation of $E$ at $[f,\eta]$.  At a critical point, this differential must vanish when restricted to $T_f\hbox{Map}(\Sigma ,M) \oplus \{ 0 \}$, and hence $f$ must be a harmonic map.  To calculate the derivative in the other direction, we will use the symbol $\eta $ for both the metric on $\Sigma $ itself and for $\det (\eta _{ab})$, the context hopefully making clear which is meant.  We differentiate (\ref{E:energyrange}), taking a perturbation given by the formula
\begin{equation} \eta_{ab} (t) = \eta _{ab} + t \dot \eta _{ab}, \label{E:varinmetric} \end{equation}
where the variation $\dot \eta _{ab}$ in the metric is trace-free ($\dot \eta _{11} + \dot \eta _{22} = 0$) and the initial metric is given by the formula
$$\eta _{ab} = \delta _{ab} \lambda ^2, \qquad  \sqrt \eta =  \sqrt{\det(\eta _{ab})} = \lambda ^2,$$
for some conformal factor $\lambda ^2$.  Then the formulae
$$\frac{d}{dt} \eta_{ab} (t) = \dot \eta _{ab},\qquad \left. \frac{d}{dt} \eta (t) \right|_{t=0} = \lambda ^2 (\dot \eta _{11} + \dot \eta _{22}) = 0.$$
imply
$$\left. \frac{d}{dt}\begin{pmatrix} \sqrt \eta \eta ^{11} & \sqrt \eta \eta ^{12} \cr \sqrt \eta \eta ^{21} & \sqrt \eta \eta ^{22} \end{pmatrix} \right|_{t=0} =
\lambda ^{-2} \begin{pmatrix} \dot \eta _{22} & - \dot \eta _{12} \cr - \dot \eta _{21} & \dot \eta _{11} \end{pmatrix},$$
and thus we find that
\begin{multline} \left. \frac{d}{dt} E(f,\eta _{ab}(t))\right|_{t=0}  = \frac{1}{2} \int _{\Sigma} \sum _{a,b} \left. \frac{d}{dt} \sqrt \eta \eta ^{ab} \right|_{t=0} \left\langle \frac{\partial f}{\partial x_a}, \frac{\partial f}{\partial x_b} \right\rangle dx_1 dx_2 \\ = - \frac{1}{2} \int _{\Sigma} \left[ \frac{\dot \eta _{11}}{\lambda ^2} \left(\left\langle \frac{\partial f}{\partial x_1}, \frac{\partial f} {\partial x_1} \right\rangle - \left\langle \frac{\partial f}{\partial x_2}, \frac{\partial f}{\partial x_2} \right\rangle \right) + \frac{2 \dot \eta _{12}}{\lambda ^2} \left\langle \frac{\partial f}{\partial x_1}, \frac{\partial f}{\partial x_2} \right\rangle \right] dx_1 dx_2 \\
= - 2  \int _{\Sigma} \frac{1}{\lambda ^2} \hbox{Re} \left[ \left(\dot \eta _{11} + i\dot \eta _{12}\right) \left\langle \frac{\partial f}{\partial z} , \frac{\partial f}{\partial z} \right\rangle \right] dx_1dx_2. \label{E:innerprodbetwdiff}\end{multline}
We define a complex bilinear map
\begin{equation}\langle \cdot , \cdot \rangle : \Gamma ({\bf K}^{(2,0)}) \times \Gamma ({\bf K}^{(0,2)}) \rightarrow {\mathbb C} \quad \hbox{by} \quad \langle \phi dz^2 , \bar \psi d\bar z^2 \rangle = 4 \int _\Sigma \frac{1}{\lambda ^4} \phi \bar \psi dA,\label{E:hiponquad} \end{equation}
the real part being an $L^2$ inner product on the space of quadratic differentials.  Then we can rewrite (\ref{E:innerprodbetwdiff}) as
\begin{equation} \left. \frac{d}{dt} E(f,\eta _{ab}(t))\right|_{t=0}  = - \hbox{Re} \frac{1}{2} \langle \Omega _f , (\dot \eta _{11} - i\dot \eta _{12}) dz^2 \rangle = - \hbox{Re} \langle \Omega _f , \phi dz^2 \rangle ,\label{E:p} \end{equation}
where $\Omega _f$ is the Hopf differential, which can be defined even if $f$ is not harmonic.  Thus the vanishing of the restriction of $dE([f,\eta])$ to $\{ 0\} \oplus T_\omega {\mathcal T}_g$ is equivalent to the condition
\begin{equation} P(\Omega_f ) = 0,\label{projtoholoquaddiff}\end{equation}
where $P$ is the $L^2$-orthogonal projection to holomorphic quadratic differentials.  Conformal maps satisfy this condition even if $f$ is not harmonic.

At a critical point for $E$ itself, the Hopf differential $\Omega _f$, which is now holomorphic, must vanish.  We thus recover the familiar fact that critical points for $E$ are indeed conformal harmonic maps, possibly with branch points.  In summary:

\vskip .1 in
\noindent
{\bf Proposition~4.1.1.} {\sl The first variation for energy is
\begin{multline} dE(f,\omega)\left((X,\dot \omega )\right) = \int _\Sigma \left[ \left\langle \frac{\partial f} {\partial x}, \frac{DX}{\partial x}\right\rangle + \left\langle \frac{\partial f} {\partial y} , \frac{DX}{\partial y} \right\rangle \right] dx dy  - \hbox{Re} \langle \Omega _f, \bar \phi d\bar z^2 \rangle  \\ = - \int _\Sigma \left[ \left\langle \frac{D}{\partial x} \left( \frac{\partial f} {\partial x}\right) + \frac{D}{\partial y} \left( \frac{\partial f} {\partial y}\right) , X \right\rangle \right] dx dy - \hbox{Re} \langle \Omega _f, \bar \phi d\bar z^2 \rangle. \label{eq:firstvarofE} \end{multline}
Here $\bar \phi d\bar z^2 $ is an antiholomorphic quadratic differential.}

\vskip .1 in
\noindent
Projection to the real part enables us to identify the space of antiholomorphic quadratic differentials with the tangent space to Teichm\"uller space at $\omega $.

\subsection{Second variation of two-variable energy}
\label{S:secondvariationtwovar}

The second derivative of energy at a critical point $(f,\eta )$ for the energy
$$E :  \hbox{Map}(\Sigma ,M) \times \hbox{Met}(\Sigma ) \longrightarrow {\mathbb R}$$
is a sum of three terms.  If $X$ is a section of $f^*TM$ and $\dot \eta$ is a trace-free variation of the metric on $\Sigma $,
\begin{multline} d^2E(f,\eta)((X,\dot \eta),(X,\dot \eta)) = d^2E_\eta (f)(X,X) \\ + 2 \int _{\Sigma} \sum _{a,b} \left. \frac{d}{dt} \left(\eta ^{ab} \sqrt{\det(\eta _{ab})}\right)\right|_{t=0} \left\langle \frac{DX}{\partial x_a}, \frac{\partial f} {\partial x_b} \right\rangle dx_1 dx_2 \\ + \frac{1}{2} \int _{\Sigma} \sum _{a,b} \left. \frac{d^2}{dt^2} \left(\eta ^{ab} \sqrt{\det(\eta _{ab})}\right)\right|_{t=0} \left\langle \frac{\partial f} {\partial x_a}, \frac{\partial f} {\partial x_b} \right\rangle dx_1 dx_2, \label{E:fullsecondder}\end{multline}
where $E_\eta $ denotes energy when the metric $\eta $ on $\Sigma $ is fixed.  We can think of these three terms as second-order partial derivatives in the three directions:
$$\hbox{Map}(\Sigma ,M) \times \hbox{Map}(\Sigma ,M), \quad \hbox{Map}(\Sigma ,M) \times \hbox{Met}(\Sigma ), \quad \hbox{Met}(\Sigma ) \times \hbox{Met}(\Sigma ).$$
The first term in (\ref{E:fullsecondder}) is the familiar index formula for harmonic maps,
\begin{equation} d^2E _{\eta }(f)(X,X) = \int _\Sigma \left[ \| DX \|^2 - \langle {\mathcal R}(X\wedge df), X\wedge df \rangle \right] dA, \label{E:classicalsecvar} \end{equation}
where in terms of the complex parameter $z = x_1 + ix_2$ on $\Sigma $,
$$\| DX \|^2 = \frac{1}{\lambda ^2} \left[ \left| \frac{DX}{\partial x_1} \right|^2 + \left| \frac{DX}{\partial x_2} \right|^2 \right]$$
and
$$\langle {\mathcal R}(X\wedge df), X\wedge df \rangle = \frac{1}{\lambda ^2} \left[ \left\langle R \left(X, \frac{\partial f} {\partial x_1} \right)\frac{\partial f} {\partial x_1},X\right\rangle + \left\langle R \left( X, \frac{\partial f} {\partial x_2} \right)\frac{\partial f} {\partial x_2},X \right\rangle \right],$$
$R$ being the Riemann-Christoffel curvature tensor of $M$.  To evaluate the other two terms, we take a metric variation
$$\eta_{ab} (t) = \lambda ^2 \delta _{ab} + t \dot \eta _{ab},$$
where $\dot \eta _{11} + \dot \eta _{22} = 0$.  Differentiation then yields
\begin{equation} \left.  \frac{d}{dt} \eta (t) \right|_{t=0} = 0, \qquad \left. \frac{d^2}{dt^2} \eta(t) \right|_{t=0} = 2\dot \eta _{11} \dot \eta _{22} - 2\dot \eta _{12} \dot \eta _{12}, \label{E:directcalc}\end{equation}
where we abuse notation a bit, and use the symbol $\eta (t)$ to denote $\det (\eta _{ab}(t))$.  Using the formulae
$$\frac{d}{dt} (\eta ^{-1/2}) = - \frac{1}{2} \eta ^{-3/2} \frac{d\eta }{dt}, \qquad \frac{d^2}{dt^2} (\eta ^{-1/2}) = - \frac{1}{2} \eta ^{-3/2} \frac{d^2\eta }{dt^2} + \frac{3}{4} \eta ^{-5/2} \left(\frac{d\eta }{dt}\right)^2,$$
we can differentiate the components appearing in (\ref{E:fullsecondder}),
$$\begin{pmatrix} \sqrt \eta \eta ^{11} & \sqrt \eta \eta ^{12} \cr \sqrt \eta \eta ^{21} & \sqrt \eta \eta ^{22} \end{pmatrix} = \eta ^{-1/2} \begin{pmatrix} \eta _{22} & - \eta _{12} \cr - \eta _{21} & \eta _{11} \end{pmatrix},$$
obtaining
$$\frac{d}{dt} \begin{pmatrix} \sqrt \eta \eta ^{11} & \sqrt \eta \eta ^{12} \cr \sqrt \eta \eta ^{21} & \sqrt \eta \eta ^{22} \end{pmatrix} = \frac{d}{dt} (\eta ^{-1/2})  \begin{pmatrix} \eta _{22} & - \eta _{12} \cr - \eta _{21} & \eta _{11} \end{pmatrix} - \eta ^{-1/2} \begin{pmatrix} \dot \eta _{11} & \dot \eta _{12} \cr \dot \eta _{21} & \dot \eta _{22} \end{pmatrix},$$
with a similar result for the second derivative.  Evaluation at $t = 0$ and use of (\ref{E:directcalc}) yields the results
\begin{equation}\left. \frac{d}{dt}\left( \sqrt \eta \eta ^{ab} \right) \right|_{t=0} = - \lambda ^{-2} \dot \eta _{ab},\label{eq:first}\end{equation}
\begin{equation} \left. \frac{d^2}{dt^2} \left( \sqrt \eta \eta ^{ab} \right) \right|_{t=0} = \lambda ^{-4} (\dot \eta _{11}^2 + \dot \eta _{12}^2)\delta _{ab}. \label{E:second}\end{equation}
It now follows from (\ref{eq:first}) that the second term in (\ref{E:fullsecondder}) is
$$ - 2 \int _{\Sigma} \sum _{a,b} \frac{\dot \eta _{ab}}{\lambda ^2} \left\langle \frac{DX}{\partial x_a}, \frac{\partial f} {\partial x_b} \right\rangle dx_1 dx_2, $$
while from (\ref{E:second}) and the trace-free condition we conclude that the third term is
$$\frac{1}{2} \int _\Sigma \sum _{a,b}  \left| \frac{\dot \eta _{ab}}{\lambda ^2} \right|^2 \left\| \frac{\partial f} {\partial x_b} \right\| ^2dx_1dx_2.$$
Adding the three terms together yields our second variation formula
\begin{multline} d^2E(f,\eta)((X,\dot \eta),(X,\dot \eta)) = d^2E_\eta(f)(X,X) \\ - 2 \int _{\Sigma} \sum _{a,b} \frac{\dot \eta _{ab}}{\lambda ^2} \left\langle \frac{DX}{\partial x_a}, \frac{\partial f} {\partial x_b} \right\rangle  dx_1 dx_2 \\ + \int _\Sigma \frac{1}{\lambda ^4} (\dot \eta _{11} ^2 + \dot \eta _{12}^2 ) \sigma ^2 dx_1 dx_2, \label{E:fullsecondderreal}\end{multline}
where
$$\sigma ^2 = \left| \frac{\partial f}{\partial x_1} \right|^2 =  \left| \frac{\partial f}{\partial x_2} \right|^2 = 2  \left| \frac{\partial f}{\partial z} \right|^2.$$

For complex variations, we rewrite (\ref{E:secondvariationforharmonic}) as
\begin{equation} d^2E _{\eta }(f)(Z,\bar Z) = \int _\Sigma \left[ 2\| D''Z \|^2 - \langle {\mathcal R}(Z \wedge \partial f), \overline{Z \wedge \partial f} \rangle \right] dA,\label{E:cxsecondvariation2} \end{equation}
where
$$\| D''Z \|^2 = \frac{2}{\lambda ^2}\left| \frac{DZ}{\partial z} \right|^2 \quad \hbox{and} \quad \langle {\mathcal R}(Z \wedge \partial f), \overline{Z \wedge \partial f} \rangle = \frac{4}{\lambda ^2}\left\langle R \left(Z, \frac{\partial f} {\partial z} \right)\frac{\partial f} {\partial \bar z},\overline Z \right\rangle.$$
In this formula, $z = x+iy$ is a complex parameter on $\Sigma $, and $Z$ is a section of $f^*TM \otimes {\mathbb C}$, so
\begin{multline*} D''Z \in \Omega ^{0,1}(f^*TM \otimes {\mathbb C}), \\ Z \wedge \partial f \quad \hbox{and} \quad {\mathcal R}(Z \wedge \partial f) \in \Omega ^{(0,1)}(\Lambda^2f^*TM \otimes {\mathbb C}).\end{multline*}
Of course, we can also rewrite this as
$$d^2E _{\eta }(f)(Z,\bar Z) = \int _\Sigma \left[ 2\| D'Z \|^2 - \langle {\mathcal R}(Z \wedge \bar \partial f), \overline{Z \wedge \bar \partial f} \rangle \right] dA,$$
where
$$\| D'Z \|^2 = \frac{2}{\lambda ^2}\left\| \frac{DZ}{\partial \bar z} \right\|^2 \quad \hbox{and} \quad \langle {\mathcal R}(Z \wedge \bar \partial f), \overline{Z \wedge \bar \partial f} \rangle = \frac{4}{\lambda ^2}\left\langle R \left(Z, \frac{\partial f} {\partial \bar z} \right)\frac{\partial f} {\partial z},\overline Z \right\rangle.$$

Using the various powers ${\bf K}^{(p,q)} = {\bf K}^p \otimes{\bf \bar K}^q$ of the canonical bundle on $\Sigma $ (as described in \S \ref{S:diffeogroup}), and the locally defined sections
$$\frac{\partial f}{\partial z} \quad \hbox{of} \quad {\bf E} = f^*TM \otimes {\mathbb C},$$
which generate the holomorphic line bundle ${\bf L} \subseteq {\bf E}$, we now define linear maps
$$\iota : \Gamma ({\bf K}^{(2,0)}) \longrightarrow \Omega ^{1,0}(\overline{\bf L}) \subseteq \Omega ^{1,0}(\overline{\bf E}), \qquad \iota : \Gamma ({\bf K}^{(0,2)}) \longrightarrow \Omega ^{0,1}({\bf L}) \subseteq \Omega ^{0,1}({\bf E})$$
by
$$\iota \left( \phi dz^2 \right) = \frac{\phi dz^2}{\lambda ^2 dz d\bar z}\frac{\partial f}{\partial \bar z} d\bar z = \frac{\phi}{\lambda ^2}\frac{\partial f}{\partial \bar z}dz, \qquad \iota \left( \bar \psi d\bar z^2 \right) = \frac{\bar \psi}{\lambda ^2}\frac{\partial f}{\partial z}d\bar z.$$
With these preparations out of the way, we claim:

\vskip .1 in
\noindent
{\bf Proposition 4.2.1.} {\sl The second variation of the two-variable energy $E$ is
\begin{multline}  d^2E(f,\omega)\left( \left(Z,\phi dz^2, \bar \psi d\bar z^2 \right), \left(\bar Z ,\psi dz^2, \bar \phi d\bar z^2 \right) \right) \\ = \int _\Sigma \left[ 2\| D''Z - \iota \left( \bar \psi d\bar z^2 \right) \|^2 - \langle {\mathcal R}(Z \wedge \partial f), \overline{Z \wedge \partial f} \rangle \right] dA \\ 
- 4 \hbox{Re} \int_{\Sigma} \left\langle D'Z, \iota \left( \phi dz^2 \right) \right\rangle dA + \int _\Sigma 2 \left\| \iota \left( \phi dz^2 \right) \right\|^2 dA. \label{E:fullseconddercomplex}\end{multline}}

\vskip .1 in
\noindent
{\bf Note:} We can interpret this formula in two ways.  In both cases, $Z$ is a section of ${\bf E} = f^*TM \otimes {\mathbb C}$.  But on the one hand we can take
$$(\phi dz^2, \bar \psi d\bar z^2) \in T_\eta (\hbox{Met}_0(\Sigma )) \otimes {\mathbb C}$$
if we want the full second variation formula on $\hbox{Map}(\Sigma ,M) \times \hbox{Met}_0(\Sigma )$, or we can require that $\phi dz^2$ and $\psi dz^2$ be holomorphic quadratic differentials, so that
$$(\phi dz^2, \bar \psi d\bar z^2) \in T_{[\eta]}{\mathcal T}_g \otimes {\mathbb C},$$
if we want to divide out by the infinite-dimensional group $\hbox{Diff}_0(\Sigma )$ and consider second variation on $\hbox{Map}(\Sigma ,M) \times{\mathcal T}_g$.

\vskip .1 in
\noindent
We can expand the first term in this second variation formula to obtain an equivalent version of the formula,
\begin{multline} d^2E(f,\eta)\left( \left(Z,\phi dz^2, \bar \psi d\bar z^2 \right), \left(\bar Z ,\psi dz^2, \bar \phi d\bar z^2 \right) \right) = d^2E_\eta (f)(Z,\bar Z) \\
- 4 \hbox{Re} \int_{\Sigma} \left\langle D''Z, \iota \left( \bar \psi d\bar z^2 \right) \right\rangle dA - 4 \hbox{Re} \int_{\Sigma} \left\langle D'Z, \iota \left( \phi dz^2 \right) \right\rangle dA\\
+ 2 \int _\Sigma \left( \left\| \iota \left(\bar \psi d\bar z^2 \right) \right\|^2 + \left\| \iota \left( \phi dz^2 \right) \right\|^2 \right) dA. \label{E:fullsecondderreal3}\end{multline}
Our goal is to prove that this is the correct formula.  Of course, we could polarize and expand (\ref{E:fullsecondderreal3}) to a full Hermitian symmetric bilinear map.  This Hermitian symmetric bilinear map is uniquely determined by its restriction to the real subspace invariant under the conjugate linear conjugation map
$$C : \left(Z,\phi dz^2, \bar \psi d\bar z^2 \right) \mapsto \left(\bar Z ,\psi dz^2, \bar \phi d\bar z^2 \right),$$
so the Hermitian map $d^2E(f,\eta)$ is completely determined by
\begin{multline} d^2E(f,\eta)\left( \left(X,\phi dz^2, \bar \phi d\bar z^2 \right), \left(X ,\phi dz^2, \bar \phi d\bar z^2 \right) \right) = d^2E_\eta (f)(X,X) \\
- 4 \hbox{Re} \int_{\Sigma} \left\langle D''X, \iota \left( \bar \phi d\bar z^2 \right) \right\rangle dA - 4 \hbox{Re} \int_{\Sigma} \left\langle D'X, \iota \left( \phi dz^2 \right) \right\rangle dA\\
+ 4 \int _\Sigma  \left\| \iota \left( \phi dz^2 \right) \right\|^2 dA, \label{E:fullsecondderreal9}\end{multline}
when $X$ is a section of $f^*TM$.  Thus it suffices to show that (\ref{E:fullsecondderreal9}) reduces to the second variation formula (\ref{E:fullsecondderreal}) we obtained before.

To do this, we utilize formula (\ref{E:defofphi0}) for $\phi$, and take the complex bilinear inner product between the two expressions
$$DX = \sum \frac{DX}{\partial x_a}dx_a = \frac{DX}{\partial z}dz + \frac{DX}{\partial \bar z}d\bar z = D'X + D''X$$
and
$$2\sum \dot \eta _{ab}  \frac{\partial f}{\partial x_b}dx_a = 4\phi \frac{\partial f}{\partial \bar z}dz + 4\bar \phi \frac{\partial f}{\partial z}d\bar z = 4\left( \iota \left( \phi dz^2 \right) + \iota \left( \bar \phi d\bar z^2 \right) \right).$$
Using the fact that $\langle dx_a,dx_b\rangle = (1/\lambda ^2)\delta _{ab}$ , we see that
\begin{multline*}  - 2 \int _{\Sigma} \sum _{a,b} \frac{\dot \eta _{ab}}{\lambda ^2} \left\langle \frac{DX}{\partial x_a}, \frac{\partial f} {\partial x_b} \right\rangle  dx_1 dx_2 \\ = - 4 \hbox{Re} \int_{\Sigma} \left\langle D''X, \iota \left( \bar \phi d\bar z^2 \right) \right\rangle dA - 4 \hbox{Re} \int_{\Sigma} \left\langle D'X, \iota \left( \phi dz^2 \right) \right\rangle dA.\end{multline*}
On the other hand, a straightforward calculation verifies that
$$\int _\Sigma \frac{2}{\lambda ^4} (\dot \eta _{11} ^2 + \dot \eta _{12}^2 )\left| \frac{\partial f}{\partial z} \right|^2 dx_1 dx_2 = 4 \int _\Sigma  \left\| \iota \left( \phi dz^2 \right) \right\|^2 dA.$$
Thus (\ref{E:fullsecondderreal}) and (\ref{E:fullsecondderreal9}) do indeed agree, establishing what we wanted to prove.

\subsection{Tangential Jacobi fields}
\label{S:tangentialjacobi}

The second variation formula simplifies when $Z$ is a section of the line bundle ${\bf L}$ and becomes particularly interesting in the case where $f$ has branch points.  Suppose that $\bar \psi d\bar z^2$ is an antiholomorphic quadratic differential which satisfy the equation
\begin{equation} \frac{DZ}{\partial \bar z} = \frac{\bar \psi}{\lambda ^2} \frac{\partial f}{\partial z}.\label{E:tangentialJacobi1} \end{equation}
We can then set $\phi = 0$ and obtain what we claim is a Jacobi field $(Z, 0, \bar \psi d\bar z^2)$.  To prove this, we need to check that
$$d^2E(f,\omega)\left( \left(Z, 0, \bar \psi d\bar z^2 \right), \left(\bar Z_1, \psi_1 dz^2, \bar \phi_1d\bar z^2 \right) \right) = 0$$
for all $(\bar Z_1, \psi_1dz^2, \bar \phi_1d\bar z^2)$.  But the curvature term vanishes when $Z$ is a section of ${\bf L}$, no matter what $Z_1$ is, and (\ref{E:tangentialJacobi1}) implies that the only terms that could be nonzero are those which appear in
$$d^2E(f,\omega)\left( \left(Z, 0, \bar \psi d\bar z^2 \right), \left(0 , 0, \bar \phi_1d\bar z^2\right) \right).$$
Indeed, we see that the only term in this expression that does not immediately vanish is 
$$- 8 \hbox{Re} \int_{\Sigma} \left\langle \frac{DZ}{\partial z}, \frac{\bar \phi _1}{\lambda ^2} \frac{\partial f}{\partial z} \right\rangle dx_1 dx_2,$$
and this term vanishes because ${\bf L} \oplus {\bf \bar L}$ is a parallel decomposition of $(f^*TM)^\top \otimes {\mathbb C}$, the tangential part of $(f^*TM) \otimes {\mathbb C}$, and
$$\left\langle \frac{\partial f}{\partial z}, \frac{\partial f}{\partial z} \right\rangle = 0 \quad \hbox{or equivalently,} \quad \langle {\bf L}, {\bf L}\rangle = 0.$$
Conversely, it is not difficult to show that all Jacobi fields $(Z, \phi dz^2, \bar \psi d\bar z^2)$ in which $Z$ is a section of ${\bf L}$ are of the form $(Z, 0, \bar \psi d\bar z^2)$ where $\psi dz^2$ is a holomorphic quadratic differential which satisfies
\begin{equation} D''(Z) = \iota \left( \bar \psi d\bar z^2 \right). \label{E:Dprime} \end{equation}
We let ${\cal J}({\bf L})$ denote the complex linear space of pairs $(Z, \bar \psi d\bar z^2)$ satisfying this equation.

\vskip .1 in
\noindent
{\bf Definition 4.3.1.}  By a {\em tangential Jacobi field\/} we mean either an element $(Z, \bar \psi d\bar z^2)$ of ${\cal J}({\bf L})$ or a pair
$$(X, \dot \omega ) = \hbox{Re}(Z, \bar \psi d\bar z^2) \in \hbox{Re}\left({\cal J}({\bf L})\right),$$
the context making clear which is meant.  Note that when $\Sigma $ has genus zero or one, ${\cal J}({\bf L})$ contains a subspace of elements of the form $(Z,0)$, where $Z$ is a holomorphic section generated by the action of $G$.

\vskip .1 in
\noindent
Whenever a parametrized minimal surface $(f,\eta )$ has a nontrivial branch locus, the branch points give rise to nontrivial tangential Jacobi fields.  Recall that the first Chern class of the line bundle ${\bf L}$ depends on the total branching order $\nu $ 0f f; if $g$ is the genus of $\Sigma $ and $\nu $ is the total branching order,
$$c_1({\bf L})[\Sigma ] = 2-2g + \nu.$$

\noindent
{\bf Proposition 4.3.2.} {\sl The complex dimension of the space of tangential Jacobi fields to a critical point $(f,\eta )$ for the two-variable energy $E$ is given by the formula:
$$ \hbox{dim}_{\mathbb C} {\cal J}({\bf L}) =  \begin{cases} 3 + \nu & \hbox{if $g=0$}, \cr 1 + \nu & \hbox{if $g=1$}, \cr \nu & \hbox{if $g\geq 2$}. \cr\end{cases}$$
}

\vskip .1 in
\noindent
This proposition is an application of the Riemann-Roch theorem to the operator
$$D'' : \Gamma ({\bf L}) \rightarrow \Gamma ({\bf L} \otimes {\bf K}^{0,1}).$$
When the total branching order $\nu $ is zero, the cokernel of $D''$ is just the space of holomorphic quadratic differentials, and the only solutions to (\ref{E:tangentialJacobi1}) are the holomorphic sections ${\mathcal O}(\bf{K}^{-1})$ of the holomorphic tangent bundle $\bf{K}^{-1}$.  This space of holomorphic sections is just the tangent space to the identity component $G$ of the group of conformal automorphisms, which is $PSL(2,{\mathbb C})$ when $g = 0$, $S^1 \times S^1$ when $g = 1$ and trivial when $g \geq 2$.  Thus when $\nu = 0$, $\hbox{dim}_{\mathbb C} {\cal J}({\bf L})$ is three, one or zero when the genus $g$ is zero, one or $\geq 2$ respectively.

Suppose next that $(f,\eta)$ is a conformal harmonic map with a single branch point $p$ of multiplicity one.  Then ${\bf L} =  {\bf \bar K} \otimes \zeta _p$, where $\zeta _p$ is the point bundle at $p$, and according to the Riemann-Roch theorem, the analytic index of $D''$ is given by the formula,
$$\dim {\mathcal O}({\bf L}) - \dim {\mathcal O}({\bf L}^{-1} \otimes {\bf K}) = 1-g + c_1({\bf L})[\Sigma ].$$
Either the dimension of the kernel of $D''$ increases by one which yields a new element of ${\cal J}({\bf L})$, or the dimension of the image of
$$D'' : \Gamma ({\bf \bar K} \otimes \zeta _p) \rightarrow \Gamma ({\bf \bar K} \otimes {\bf \bar K} \otimes \zeta _p)$$
decreases by one.  In the latter case, at least one nonzero antiholomorphic quadratic differential is covered by $D''$, which also gives a new element of ${\cal J}({\bf L})$.  (The space of new elements of ${\cal J}({\bf L})$ obtained increasing $\nu$ by one is at most one-dimensional, because given two such solutions, a nontrivial linear combination would vanish at $p$ and therefore lie within ${\cal J}({\bf \bar K})$.)

The general case of the Lemma is proven by induction on $\nu$.  Each time $\nu $ is increased by one, either a new holomorphic section of ${\bf L}$ is produced or a new antiholomorphic quadratic differential is covered by $D''$, in either case leading to an increase of one in the dimension of ${\cal J}({\bf L})$.

\vskip .1 in
\noindent
{\bf Proposition 4.3.3.} {\sl If $f : \Sigma \rightarrow M$ is a weakly conformal harmonic map with divisor of branch points
$$D(f) = \nu _1 p_1 + \cdots + \nu _n p_n,$$
then there exist linearly independent elements
$$(Z_i, \bar \psi _i d\bar z^2) \in {\cal J}({\bf L}) \quad \hbox{for $1 \leq i \leq n$, such that} \quad Z_i(p_j) = \delta _{ij} {\bf e}_i,$$
where ${\bf e}_i$ is a unit-length element chosen from ${\bf L}_{p_i}$.}

\vskip .1 in
\noindent
This is an immediate consequence of a stronger assertion:  If $z_1, \ldots z_n$ are local complex coordinates on $\Sigma $ centered at $p_1, \ldots p_n$, then we can define a map
\begin{multline} \Psi : {\cal J}({\bf L}) \longrightarrow \sum _{i=1}^n \overbrace{{\bf L}_{p_i} \oplus \cdots \oplus {\bf L}_{p_i}}^{\nu _i} \\ \hbox{by} \quad \Psi \left(Z, \bar \psi d\bar z^2\right) = \left(\ldots , \left(Z(p_i), \frac{DZ}{\partial z_i}, \ldots ,\frac{D^{\nu_i-1}Z}{\partial z_i^{\nu_i-1}}\right) , \ldots \right). \label{E:psiisomorphism}\end{multline}

\vskip .1 in
\noindent
{\bf Proposition 4.3.4.} {\sl The linear map $\Psi $ defined by (\ref{E:psiisomorphism}) is surjective, and its kernel is the space of holomorphic sections generated by the holomorphic automorphisms of $\Sigma $.}

\vskip .1 in
\noindent
Propositions 4.3.2 and 4.3.3 can be proven first for simple branch points by induction, then by a second induction on the total branching order for a given choice of $n$ branch points.

From the point of view of the Riemann surface $\Sigma $, the section $Z$ appearing in $\left(Z, \bar \psi d\bar z^2\right)$ has poles at the branch points.  Roughly speaking, the fact that $\Psi $ is an isomorphism can be regarded as stating that any combination of principal parts can be realized by elements of ${\cal J}({\bf L})$ at the branch points.

\vskip .1 in
\noindent
{\bf Definition 4.3.5.}  By ${\cal J}_0({\bf L})$ we will denote the kernel of the map
\begin{multline} \Psi _0: {\cal J}({\bf L}) \longrightarrow \sum _{i=1}^n {\bf L}_{p_i} \\ \hbox{defined by} \quad \Psi \left(Z, \bar \psi d\bar z^2\right) = \left(Z(p_1), Z(p_2), \ldots , Z(p_n) \right). \label{E:psiisomorphism1}\end{multline}
It is a vector space of real dimension $2(\nu - n) + \dim G$.  We will sometimes let ${\cal J}_1({\bf L})$ denote a $2n$-dimensional complement of ${\cal J}_0({\bf L})$ within ${\cal J}({\bf L})$.

\subsection{Minimal surfaces without branch points}
\label{S:minimalnobranch}

Now that we have the second variation formula for parametrized minimal surfaces at our disposal, it is easy to modify the proof schema presented in \S\ref{S:harmonic} to give the corresponding result for minimal surfaces. 

\vskip .1 in
\noindent
{\bf Proposition~4.4.1.} {\sl Suppose that $M$ is a compact connected manifold of dimension at least three.  For a generic choice of Riemannian metric on $M$, every prime oriented compact parametrized minimal surface $f : \Sigma \rightarrow M$ with no branch points is as nondegenerate as allowed by its connected group $G$ of symmetries.  If $G$ is trivial, it is Morse nondegenerate for the two-variable energy $E$ in the usual sense.  If $G = S^1\times S^1$ or $PSL(2,{\mathbb C})$, then all such minimal surfaces lie on nondegenerate critical submanifolds for $E$ which are orbits for the $G$-action.}

\vskip .1 in
\noindent
The proof of this proposition follows the same pattern as that used for Proposition~3.1. except that we replace the Euler-Lagrange map $F$ utilized in \S\ref{S:harmonic} by the map
\begin{equation} F : \hbox{Map}'(\Sigma ,M) \times {\mathcal T}_g \times \hbox{Met}(M) \longrightarrow \hbox{Map}(\Sigma ,TM) \times T{\mathcal T}_g \label{E:defofFmin}\end{equation} 
obtained by first variation of $E$, as described in Proposition~4.1.1, $\hbox{Map}'(\Sigma ,M)$ being the space of somewhere injective immersions and ${\mathcal T}_g$ the Teichm\"uller space of the closed oriented surface $\Sigma $ of genus $g$.  We can divide $F$ into two components
$$F_0 : \hbox{Map}'(\Sigma ,M) \times {\mathcal T} \times \hbox{Met}(M) \longrightarrow \hbox{Map}(\Sigma ,TM),$$
$$F_1 : \hbox{Map}'(\Sigma ,M) \times {\mathcal T} \times \hbox{Met}(M) \longrightarrow T{\mathcal T},$$
the first of which is the negative of the tension
\begin{equation} F_0(f,\omega ,g) = - \frac{1}{\lambda ^2} \left[ \frac{D^g}{\partial x_1} \left(\frac{\partial f}{\partial x_1}\right) + \frac{D^g}{\partial x_2} \left(\frac{\partial f}{\partial x_2}\right) \right], \label{eq:familiarmap} \end{equation}
where $(x_1,x_2)$ are isothermal coordinates with respect to $\omega $.  In accordance with (\ref{E:hiponquad}) and (\ref{E:p}), the second component is simply the map
$$F_1(f, \omega ,g) = - P\left( \left\langle \frac{\partial f}{\partial z}, \frac{\partial f}{\partial z} \right\rangle dz^2\right) \in T_\omega{\mathcal T},$$
where $P$ is the $L^2$-orthogonal projection to the space of holomorphic quadratic differentials (\ref{projtoholoquaddiff}).  The critical points for the energy $E$ are the points $(f,\omega ,g)$ such that $F(f,\omega ,g)$ is a zero-section of the tangent bundle to $\hbox{Map}'(\Sigma ,M) \times {\mathcal T}$.

The main step should be to show that
\begin{multline} {\mathcal P}^\emptyset = \{ (f,\omega ,g) \in \hbox{Map}(\Sigma ,M) \times {\mathcal T} \times \hbox{Met}(M) : \hbox{$f$ is a prime conformal} \\ \hbox{$\omega $-harmonic immersion for $g$ with no branch points} \}, \label{E:mathcalSprime} \end{multline}
is a submanifold, where the superscript $\emptyset $ denotes that the branch locus is empty.  (Note that prime minimal surfaces are automatically somewhere injective, so elements of ${\mathcal P}^\emptyset$ are injective on an open dense subset of $\Sigma $.)  This time we take
\begin{multline*} {\mathcal E} = \{ (f,\omega ,g,X,\dot \omega) \in (L^2_1)'(\Sigma ,M) \times {\mathcal T} \times \hbox{Met}(M) \times \hbox{Map}(\Sigma ,TM) \times T{\mathcal T} \\ \hbox{such that $\pi (X) = f$ and $\pi (\dot \omega ) = \omega $ }  \},\end{multline*}
and regard the map $F$ of (\ref{E:defofFmin}) as defining a section of the bundle
$${\mathcal E} \longrightarrow \hbox{Map}'(\Sigma ,M) \times {\mathcal T} \times \hbox{Met}(M).$$
If $(f,\omega, g)$ is a zero for $F$, we let
\begin{multline} {\mathcal V} = \{ (X, \dot \omega ) \in T_f\hbox{Map}'(\Sigma ,M) \oplus T_\omega {\mathcal T} : \hbox{ $(X, \dot \omega )$ is perpendicular}\\
\hbox{ to  the image of $\pi _V \circ DF_{(f,\omega ,g)}$} \}, \label{E:defofVminimal} \end{multline}
where $\pi _V$ is projection on the fiber.  The following Lemma will show that ${\mathcal V} = {\mathcal G}(f,\omega ,g)$ where
$${\mathcal G}(f,\omega ,g) \subseteq T_f\hbox{Map}'(\Sigma ,M) \oplus T_\omega {\mathcal T} \subseteq {\mathcal E}_{(f,\omega ,g)}$$
is the space generated by the $G$-action, $G$ being the group of symmetries of $\Sigma $.  If we restrict the map $F$ of (\ref{E:defofFmin}) to
$$F : \hbox{Map}'(\Sigma ,M) \times {\mathcal T} \times \hbox{Met}(M) \longrightarrow \hbox{Map}(\Sigma ,TM) \times T{\mathcal T},$$
it will then follow that
$$(\hbox{image of $DF_{(f,\omega ,g)}$}) + {\mathcal G}(f,\omega ,g) + T_{(f,\omega ,g)}{\mathcal Z} = T_{(f,\omega, g,0,0)}{\mathcal E},$$
allowing us to conclude from the implicit function theorem that
\begin{multline}{\mathcal S} = \{ (f,\omega ,g) \in {\mathfrak M}(\Sigma ,M) \times {\mathcal T} \times \hbox{Met}(M) : \hbox{$f$ is a prime conformal} \\ \hbox{$\omega $-harmonic immersion for $g$ with no branch points} \} \label{E:scriSminnobranch}\end{multline}
is a submanifold.

\vskip .1 in
\noindent
{\bf Lemma~4.4.2.} {\sl If $(f, \omega )$ is a critical point for $E$ with $f$ somewhere injective, the only elements of $T_{(f,\omega )}( \hbox{Map}(\Sigma ,M) \times {\mathcal T})$ that can be perpendicular to the image of every variation in the metric are the tangential Jacobi fields $(X, \dot \omega)$, which correspond to elements of ${\mathcal J}({\bf L})$.}

\vskip .1 in
\noindent
We prove this lemma first under the assumption that $f$ has no branch points.  Recall that the derivative of $F$ at a critical point $(f, \omega ,g)$ in the direction of $\hbox{Map}(\Sigma ,M) \times {\mathcal T}$, which we denote by $D_1F$, can be determined from the second variation formula for $E$ of \S\ref{S:secondvariationtwovar}.  This determines the Jacobi operator,
$$L_{f, \omega ,g} = (\pi _V \circ D_1F)_{(f, \omega ,g)} : \Gamma (f^*TM) \oplus T_\omega {\mathcal T} \rightarrow \Gamma (f^*TM) \oplus T_\omega {\mathcal T},$$
which divides into two components $L_0$ and $L_1$, taking values in $\Gamma (f^*TM)$ and $T_\omega {\mathcal T}$ respectively.

To calculate the derivative $D_2F$ with respect to the metric, we need to divide into components $D_2F_0$ and $D_2F_1$.  The component $D_2F_0$ is treated just as in \S\ref{S:harmonic}, where we note explicitly that (\ref{E:conformalfactor}) applies since $f$ is conformal as well as harmonic.  In particular, since $f$ is somewhere injective, the first part of the argument for Lemma~3.3 shows that if $(X,\dot \omega )$ is an element of ${\mathcal V}$, then $X$ must be tangent to $\Sigma $.  On the other hand, if a tangential variation $\dot g$ in the metric on $M$ is represented by the (not necessarily holomorphic) quadratic differential
$$\phi dz^2 = \frac{1}{2}(\dot g _{11} - i \dot g_{12}) dz^2, \qquad \hbox{with} \qquad \dot g_{11} + \dot g_{22} = 0,$$
the second part of the argument (in particular (\ref{E:desiredformula})) shows that
\begin{equation} \int _\Sigma \langle \pi _V \circ (D_2F_0)_{(f,\omega ,g)} (\dot g),Z \rangle dA = \int _\Sigma \frac{4}{\sigma ^2} \left\langle \frac{DZ}{\partial \bar z} , \frac{\phi }{\lambda ^2} \frac{\partial f}{\partial \bar z} \right\rangle dA,\label{E:compare1} \end{equation}
under the assumption that $\sigma ^2$ is nonzero, which holds when there are no branch points.

To calculate $D_2F_1$, we let  
$$\bar \psi d\bar z^2 = \frac{1}{2}(\dot \eta _{11} + i \dot \eta _{12})dz^2$$
be an antiholomorphic quadratic differential and set $\dot \eta _{22} = - \dot \eta _{11}$.  We then consider a smooth family of metrics $t \mapsto g(t) = \langle \cdot , \cdot \rangle _t$, and the corresponding family of Hopf differentials
\begin{multline*} \Omega _f(t) = \frac{1}{4} \left[ \left\langle \frac{\partial f}{\partial x_1} , \frac{\partial f}{\partial x_1} \right\rangle_t - \left\langle \frac{\partial f}{\partial x_2} , \frac{\partial f}{\partial x_2} \right\rangle_t  - 2i \left\langle \frac{\partial f}{\partial x_1} , \frac{\partial f}{\partial x_2} \right\rangle _t \right] dz^2 \\ = \frac{1}{4}\left[ (g_{11}(t) - g_{22}(t)) - 2i g_{12}(t)\right]dz^2.\end{multline*}
Using the inner product (\ref{E:hiponquad}), we find that
\begin{multline*} \hbox{Re} \langle F_1(f, \omega ,g(t)), \bar \psi d\bar z^2 \rangle = - \hbox{Re} \langle \Omega _f(t) , \dot \omega \rangle \\ = - \int _\Sigma \frac{1}{\lambda ^2}[ \dot \eta _{11} ( g_{11}(t) - g_{22}(t)) + 2 \dot \eta _{12} g_{12}(t)]dx_1dx_2. \end{multline*}
We differentiate this with respect to $t$, letting $\dot g$ denote the derivative of $g(t)$ at $t = 0$, obtaining the result
\begin{multline} \hbox{Re} \langle D_2F_1(f, \omega ,g)(\dot g), \bar \psi d\bar z^2 \rangle = - \int _\Sigma \frac{1}{ \lambda ^2}[ \dot \eta _{11} ( \dot g_{11} - \dot g_{22}) + 2 \dot \eta _{12} \dot g_{12}]dx_1dx_2
\\  = -  \int _\Sigma \frac{2}{\lambda ^4}\left[ \dot \eta _{11} \dot g_{11} + \dot \eta _{12} \dot g_{12}  \right] dA. \label{E:d2f1simple}\end{multline}
Since (\ref{E:conformalfactor}) applies to $f$,
$$\left\langle \frac{\partial f}{\partial z}, \frac{\partial f}{\partial \bar z} \right\rangle = \frac{\sigma ^2}{2},$$
and we can therefore infer that
\begin{equation} \hbox{Re} \langle D_2F_1(f, \omega ,g)(\dot g), \bar \psi d\bar z^2 \rangle = - \hbox{Re} \int _\Sigma \frac{4}{\sigma ^2} \left\langle \frac{\bar \psi }{\lambda ^2} \frac{\partial f}{\partial z}, \frac{\phi }{\lambda ^2} \frac{\partial f}{\partial \bar z} \right\rangle dA. \label{E:compare2} \end{equation}

From (\ref{E:compare1}) and (\ref{E:compare2}) we conclude that if $Z$ is a section of ${\bf L}$ and $\bar \psi d\bar z^2$ is an antiholomorphic quadratic differential, then
$$ \langle D_2F(f, \omega ,g)(\dot g), (Z, \bar \psi d\bar z^2) \rangle = 4 \int _\Sigma \left\langle \frac{DZ}{\partial \bar z} - \frac{\bar \psi }{\lambda ^2} \frac{\partial f}{\partial z}, \frac{\phi}{\lambda ^2 \sigma ^2} \frac{\partial f}{\partial \bar z} \right\rangle dA.$$
Thus if $(X, \dot \omega) = \hbox{Re}(Z,\bar \psi d\bar z^2)$, then $(X, \dot \omega )$ is perpendicular to all variations in the metric if and only if
\begin{equation} \frac{DZ}{\partial \bar z} = \frac{\bar \psi}{\lambda ^2} \frac{\partial f}{\partial z}, \label{E:eqnforjacminimal}\end{equation}
which is the condition that $(X, \dot \omega )$ be an element of ${\mathcal J}((f^*TM)^\top)$, the space of tangential Jacobi fields.

To handle the case where $f$ has branch points, we must restrict the metric variations $\dot g$ under consideration.  To deal with the fact that $\sigma ^2$ vanishes exactly at the branch points we multiply $\phi $ by a cutoff function which vanishes at the branch points.  We then conclude that (\ref{E:eqnforjacminimal}) holds except possibly at the branch points, but an easy removeable singularity theorem from complex analysis then shows that it holds everywhere, and Lemma~4.4.2 is proven.

\vskip .1 in
\noindent
We now use the fact that when $f$ is an immersion, the only tangential Jacobi fields are those generated by the $G$-action.  It then follows from Lemma~4.2.2 that ${\mathcal V} = {\mathcal G}$, and hence the space ${\mathcal S}$ defined by (\ref{E:scriSminnobranch}) is indeed a submanifold.  To complete the argument for Proposition~4.4.1, we need only verify the following analog of Lemma~3.5.

\vskip .1 in
\noindent
{\bf Lemma~4.4.3.} {\sl The projection on the third factor $\pi : {\mathcal S} \rightarrow \hbox{Met}(M)$ is a Fredholm map of Fredholm index $d_\Sigma $, where $d_\Sigma $ is the dimension of $G$.}

\vskip .1 in
\noindent
The proof is a straightforward modification of the proof of Lemma~3.5.

\vskip .1 in
\noindent
Finally, to finish the proof of Proposition~4.4.1, we simply apply the Sard-Smale Theorem as before.

\section{Simple branch points}
\label{S:simplebranch}

This section is the core of the article.  The ideas behind the proof of the Main Theorem are most transparent in the case of simple branch points, and we now focus on that case.  Moreover, we will assume for simplicity in this section that the self-intersection set of each conformal harmonic map for a given metric consists of isolated points, a hypothesis later removed.

\subsection{Eliminating branch points}\label{S:eliminate}
\label{S:eliminatesimple}

We would like to extend the argument of \S\ref{S:minimalnobranch} from the space ${\mathcal P}^\emptyset$ of prime minimal surfaces without branch points to the space ${\mathcal P}$ of all prime minimal surfaces.  The most natural strategy for accomplishing this would be to show that the Euler-Lagrange map
$$F : \hbox{Map}(\Sigma ,M) \times {\mathcal T} \times \hbox{Met}(M) \longrightarrow T(\hbox{Map}(\Sigma ,M) \times {\mathcal T})$$
is as transversal to the zero section of the tangent bundle as allowed by the $G$-action at every point $(f,\omega,g)$ of ${\mathcal P}$.  By Lemma~4.4.1 and earlier results, it would suffice to find variations $\dot g$ in the metric such that the image of the map
\begin{equation} \dot g \quad \mapsto \quad \pi _V\circ D_2F(f,\omega,g)(\dot g) \label{eq:perturbaway} \end{equation}
covers the tangential Jacobi fields provided by the branch points, a complement to the space ${\mathcal G}$ of tangential Jacobi fields generated by the $G$-action.  In this case the intersection of the image of $F$ with the zero section would be a submanifold ${\mathcal S}$, and the differential of the projection from ${\mathcal S}$ to the space $\hbox{Met}(M)$ of metrics would fail to be surjective at minimal surfaces with branch points, implying that such minimal surfaces do not exist for generic metrics.

It is the goal of this section to construct a family of metric variations whose image under (\ref{eq:perturbaway}) covers a large portion of the space of tangential Jacobi fields (but not the entire space when $f$ has branch points of branching order $\geq 2$).  In a rough sense, we might think of these metric variations as \lq\lq killing" the Jacobi fields, or \lq\lq perturbing away" branch points on minimal surfaces.

Where might such metric variations come from?  Since the Jacobi fields we want to eliminate are tangential, we expect that the variations in the metric should also be tangent to $\Sigma $.  Moreover, experience with minimal surfaces suggests that the variations should have a holomorphic description, so we should look for holomorphic or meromorphic sections (suitably modified) for some holomorphic line bundle over $\Sigma $.

A Riemann surface is said to be {\em conformally finite\/} if it is conformally equivalent to a compact Riemann surface with finitely many points (called {\em punctures\/}) removed.  Suppose that $\Sigma _0$ is constructed from a compact oriented surface $\Sigma $ of genus $g$ by deleting $n$ points $\{ p_1, \ldots , p_n \}$, where $n > 0$ and $n \geq 3$ if $g=0$.  As described in \cite{RS}, we can regard the Teichm\"uller space of conformally finite conformal structures on $\Sigma _0$ as the quotient
$${\mathcal T}_{g,n} = \frac{\left(\overbrace{\Sigma \times \cdots \times \Sigma }^n - \Delta^{(n)} \right) \times \hbox{Met}_0(\Sigma)}{\hbox{Diff}_{0}(\Sigma)} \cong \frac{\hbox{Met}_0(\Sigma)}{\hbox{Diff}_{0,D}(\Sigma)},$$
where $\Delta^{(n)}$ is the fat diagonal consisting of $n$-tuples $(p_1, \ldots , p_n)$ for which any two of the $p_i$'s are equal, and
$$\hbox{Diff}_{0,D} = \{ \phi \in \hbox{Diff}_0 : \hbox{$\phi (p_i) = p_i$ for each $i$} \}.$$
Well-known arguments (see for example \cite{Abi}) show that ${\mathcal T}_{g,n}$ is diffeomorphic to an open ball in Euclidean space of dimension $6g - 6 + 2n$ so long as this number is positive.  (There is a single conformal structure on the sphere minus three points.)  The extra variations ($2n$ when $\Sigma $ has genus at least two) can be thought of as follows:  If we consider a fixed conformal structure on the compact surface $\Sigma $ and allow the punctures $\{ p_1, \ldots , p_n \}$ to vary, we generally obtain a family of varying conformal structures on the punctured surface.  This family can be realized by a family of diffeomorphisms of $\Sigma $, isotopic to the identity, which do not leave given punctures fixed.

Suppose now that $f : \Sigma \rightarrow M$ is a prime parametrized minimal surface with no singularities except branch points.  Suppose moreover, that there are diffeomorphisms $\varphi $ of $\Sigma $ and $\tilde \varphi $ of $M$ such that
$$\tilde \varphi ^{-1} \circ f = f \circ \varphi ^{-1}.$$
Then invariance of the map
$$E : \hbox{Map}(\Sigma ,M) \times \hbox{Met}(\Sigma ) \times \hbox{Met}(M) \longrightarrow {\mathbb R}$$
under diffeomorphims on both domain and range implies that
$$E(f, \varphi ^*\eta,g) = E(f\circ \varphi ^{-1},\eta,g) = E(\tilde \varphi ^{-1} \circ f,\eta,g) = E(f,\eta , \tilde \varphi ^*g).$$
Thus metric deformations induced by diffeomorphisms of $M$ must be tracked by the corresponding diffeomorphisms on the punctured domain or else the conformality of $f$ will be destroyed,  Put another way, metric deformations of the range not tracked by corresponding diffeomorphisms of the domain will result in a change in conformal structure on the punctured Riemann surface $\Sigma - \{ p_1, \ldots , p_n \}$, where $\{ p_1, \ldots , p_n \}$ are the branch points, and these must perturb away the minimal surface.  It is well-known (\cite{Abi}, Chapter~2, \S 1, or \cite{Ber}) that the cotangent space to ${\mathcal T}_{g,n}$ consists of meromorphic (or antimeromorphic) quadratic differentials on the compact surface $\Sigma $ which have the possibility of simple poles only at the punctures.  The preceding discussion suggests that the metric deformations on $M$ that kill the tangential Jacobi fields under (\ref{eq:perturbaway}) should be constructed from cutoffs of meromorphic differentials with simple poles.

The argument is simplest when $f$ has no singularities except for the branch points themselves and in particular there are no self-intersections.  At the end of the section, Lemma~5.1.2 will show that isolated points of self-intersection do not interfere with the argument and in later sections we will show that isolated self-intersection sets are the only complications that can occur for generic metrics.

It remains only to construct the explicit variations.  For the actual construction, we utilize a canonical form for a minimal surface in the neighborhood of a branch point.  Recall that in terms of a local holomorphic coordinate $z$ centered at a branch point $p$ of order $\nu $, we can write $(\partial f/\partial z)(z) = z^\nu g(z)$ where $g$ is a holomorphic section with $g(0) \neq 0$, and it follows directly from Taylor's theorem that in terms of normal coordinates $(u_1, \ldots ,u_n)$ centered at $f(p)$, we can write $f$ in terms of these coordinates as
\begin{equation} (u_1 + i u_2)(z) = c z^{\nu + 1} + o_1(|z|^{\nu + 1}), \quad u_r(z) = o_1(|z|^{\nu + 1}) \quad \hbox{for} \quad 3 \leq r \leq n, \label{E:canonicalformbr} \end{equation}
$c$ being a nonzero complex constant.  Here $o_1(|z|^{k+1})$ stands for a term which is $o(|z|^{k+1})$, and which has a derivative that is $o(z^{k})$.  Stronger results can be established with more work.  Indeed, it follows from Theorem~1.3 of \cite{MW} that if we set $U(z) = (u_1(z), \ldots ,u_n(z))$, then
$$U(z) = H_{\nu + 1}(z) + H_{\nu + 2}(z) + o_2(|z|^{\nu + 2}),$$
where $H_{\nu + 1}$ and $H_{\nu + 2}$ are harmonic polynomials, where $o_2(|z|^{k+2})$ stands for a term which is $o(|z|^{k+2})$, and which has a derivative that is $o_1(|z|^{k+1})$.  It follows from the first part of the proof of Theorem~1.4 of \cite{MW}) that
$$H_{\nu + k}(z) = \beta _{\nu + k} z^{\nu +k} + \bar \beta _{\nu + k} \bar z^{\nu +k}, \qquad \hbox{for $k = 1,2$,}$$
where $\beta _{\nu + k}$ is an element of ${\mathbb C}^n$.   Since the coordinates $(u_1, \ldots , u_n)$ are normal and hence the Christoffel symbols vanish at $p$,
$$\left\langle \frac{\partial f}{\partial z},  \frac{\partial f}{\partial z} \right\rangle = 0 \quad \Rightarrow \quad \frac{\partial U}{\partial z} \cdot \frac{\partial U}{\partial z} = o_1(|z|^{2\nu + 1}),$$
and hence
$$(\nu + 1)^2 \beta_{\nu + 1} \cdot \beta _{\nu + 1} z^{2\nu} + (\nu + 1)(\nu + 2) \beta_{\nu+1} \cdot \beta _{\nu + 2} z^{2\nu + 1} = o_1(|z|^{2\nu + 1}).$$
it follows that
$$\beta_{\nu + 1} \cdot \beta_{\nu + 1} = 0 = \beta_{\nu+1} \cdot \beta _{\nu + 2}.$$
The first of these conditions implies that we can improve the expressions for the first two of the normal coordinates to
$$(u_1 + i u_2)(z) = z^{\nu + 1} + c_1 z^{\nu + 2} + c_2 \bar z^{\nu + 2} + o_2(|z|^{\nu +2}),$$
while the second shows that we can write 
$$(u_1 + i u_2)(z) = z^{\nu + 1} + c z^{\nu + 2} + o_2(|z|^{\nu +2}), \quad u_r(z) = o_1(z^{\nu + 1}) \quad \hbox{for} \quad 3 \leq r \leq n,$$
where $c$ is a complex constant.

This implies that if we regard $w = u_1 + i u_2$ as a complex coordinate on $\Sigma $ with conjugate $\bar w$, then
$$dw = \left[ (\nu + 1) z^{\nu} + c (\nu + 2) z^{\nu + 1} + o_1(|z^{\nu +1}|) \right] dz + o_1(|z|^{\nu +1}) d\bar z,$$
from which we derive a formula for the corresponding vector field,
\begin{equation} \frac{\partial }{\partial w} = \left[ \frac{1}{(\nu + 1)z^\nu} + o_1(|z|^{-\nu}) \right] \frac{\partial }{\partial z} + o_1(|z|^{-\nu + 1}) \frac{\partial }{\partial \bar z}. \label{E:partialwitoz} \end{equation}

Recall that $F$ divides into two components, $F_0$ and $F_1$, with values in $T\hbox{Map}(\Sigma ,M)$ and $T{\mathcal T}$ respectively, $F_0$ being the familiar map (\ref{eq:familiarmap}) such that $F_0(f,\omega ,g) = 0$ is the equation that $f$ be harmonic with respect to $\omega $ and $g$.  To differentiate $F_0$, we utilize the formula
$$\pi _V \circ (D_2F_0)_{(f,\omega ,g)} (\dot g) = - \frac{1}{\lambda ^2}\sum _{i,j,k=1}^n g^{kl} \dot \Gamma _{l,ij} \left( \frac{\partial u_i} {\partial x_1} \frac{\partial u_j}{\partial x_1} + \frac{\partial u_i}{\partial x_2}\frac{\partial u_j}{\partial x_2} \right) \frac{\partial }{\partial x_k}.$$
from \S\ref{S:harmonic}, $z = x_1 + i x_2$ being the usual complex coordinate on $\Sigma $.  If $z$ is centered at the branch point $p$, the canonical form implies that
$$\left( \frac{\partial u_i} {\partial x_1} \frac{\partial u_j}{\partial x_1} + \frac{\partial u_i}{\partial x_2}\frac{\partial u_j}{\partial x_2} \right) = \begin{cases} (\nu + 1)^2 r^{2\nu }+ o(r^{2\nu }) & \hbox{if $i = j = 1$ or $i = j = 2$},\cr
o(r^{2\nu }) & \hbox{otherwise,} \cr\end{cases}$$
where $r = |z|$, and hence
\begin{multline} \left\langle \pi _V \circ (D_2F)_{(f,g)} ( \dot g), \sum f^i \frac{\partial }{\partial x_i} \right\rangle dA \\ = 
- \sum _{k=1}^2 \left( (\nu +1)^2 r^{2\nu } +  o(r^{2\nu })\right) (\dot \Gamma _{k,11} + \dot \Gamma _{k,22}) f^k dx_1dx_2. \label{E:derivofF}\end{multline}
It is important to notice that the Christoffel symbols are calculated with respect to the normal coordinates $(u_1,u_2)$ on the range, not the domain, so the correct analog of (\ref{E:varinmetric1}) for a variation $\dot g$ supported near the branch point $p$ is
\begin{multline} \int _\Sigma \left\langle \pi _V \circ (D_2F)_{(f,g)} (\dot g), X \right\rangle dA \\ = - \int _\Sigma \left( (\nu +1)^2 r^{2\nu } +  o(r^{2\nu })\right) \left[\left(\frac{\partial M}{\partial u_1} + \frac{\partial N}{\partial u_2}\right) f^1 \right.\\ \left. + \left(\frac{\partial N} {\partial u_1} - \frac{\partial M} {\partial u_2}\right) f^2 \right] dx_1dx_2. \label{E:varinmetricbranch} \end{multline}
when $M$ and $N$ are defined by (\ref {E:tangentialvar}).  To differentiate $F_1$, on the other hand, we would simply use (\ref{E:d2f1simple}), a calculation to which we return in the following subsection.  We can summarize the derivative of $F_1$ in a map,
$$B_0 : (\hbox{trace-free tangential variations in $T_g\hbox{Met}(M)$}) \times \Gamma({\bf L}) \longrightarrow {\mathbb C},$$
defined by
\begin{equation} B_0 \left(\phi dz^2, Z\right) = \int _\Sigma \left\langle \pi _V \circ (D_2F)_{(f,g)} (\dot g), Z \right\rangle dA, \label{EdefinitionB0} \end{equation}
where the trace-free tangential variation $\dot g$ is represented by a quadratic differential
$$\phi dz^2 = \frac{1}{2}(\dot g _{11} - i \dot g_{12}) dz^2, \qquad \hbox{with} \qquad \dot g_{22} = - \dot g_{11}.$$

We recall that when $f: \Sigma \rightarrow M$ is a weakly conformal harmonic map with divisor of branch points $D(f) = \nu _1 p_1 + \cdots + \nu _n p_n$, Proposition~4.3.3 provides us with a collection
$$(Z_i, \bar \psi _i d\bar z^2) \in {\cal J}({\bf L}) \quad \hbox{for $1 \leq i \leq n$, such that} \quad Z_i(p_j) = \delta _{ij} {\bf e}_i,$$
where ${\bf e_i}$ is a fixed unit-length element in the fiber over $p_i$.  If we choose the complex coordinate $z_i$ on $\Sigma $ near $p_i$ so that $(\partial f/\partial z_i) = z_i^{\nu _i}g$, where $g(p_i)$ has unit length, we can take ${\bf e_i} = g(p_i)$.  Then $Z_k$ has the expansion
\begin{equation} Z_k = h_{ki}(z_i) \frac{\partial f}{\partial z_i}, \quad \hbox{where} \quad h_{ki}(z_i) = \frac{\delta _{ki}}{z_i^{\nu}}(1 + o_1(|z_i|)). \label{E:expansionofZk}\end{equation}

\vskip .1 in
\noindent
{\bf Lemma~5.1.1.} {\sl Suppose that $f: \Sigma \rightarrow M$ is a weakly conformal harmonic map with divisor of branch points $D(f) = \nu _1 p_1 + \cdots + \nu _n p_n$, and that the only singularities of $f$ are its branch points.  Then there is a collection
$$\phi _1dz^2, \ldots , \phi _ndz^2$$
of quadratic differentials on $\Sigma $ such that
\begin{equation} B_0\left(\phi _jdz^2, Z_k \right) = \delta _{jk} = \begin{cases} 1, & \hbox{if $j=k$},\cr 0, & \hbox{if $j \neq k$}, \end{cases}\label{E:tobeproved}\end{equation}
where $B_0$ is defined by (\ref{EdefinitionB0}).}

\vskip .1 in
\noindent
We first prove Lemma~5.1.1 under the assumptions $\Sigma $ has genus at least two and $f$ has no points of self-intersection.

We begin by choosing one of the branch points $p_j$ and constructing a meromorphic quadratic differential $\phi '_jdz^2$ with a pole of order one at $p_j$ and no other poles.  (The existence of such a meromorphic quadratic differential, under the assumption that $\Sigma $ has genus at least two, follows from the Riemann-Roch theorem.)  We need to multiply this quadratic differential by a cutoff function near the poles so that it determines a well-defined deformation of the metric.  To do this, we first choose a local conformal coordinate $z = x_1 + i x_2$ centered at a given branch point $p_i$ and let $r = \sqrt{x_1^2 + x_2^2}$ so that we can define the $\epsilon $-ball $B_\epsilon (p_i)$ about $p_i$ as the set of points satisfying $r \leq \epsilon $, for a small $\epsilon $, $0 < \epsilon < 1$.  We construct a smooth cutoff function $\zeta : \Sigma \rightarrow [0,1]$ which is identically one outside the $\epsilon$-ball about each $p_i$, identically zero on an $\epsilon ^2$-ball about $p_i$, and depends only on the radial coordinate $r$ near $p_i$.  We then let $\phi _j = \zeta \phi'_j$, thereby obtaining a new quadratic differential $\phi _jdz^2$ which is holomorphic except in small neighborhoods of the branch points.  Since we are assuming that $f$ is an imbedding on the support of $\zeta $, the quadratic differential $\phi _jdz^2$ can indeed be realized by a deformation of the metric on the ambient manifold $M$.  Locally, such a deformation is defined by (\ref{E:tangentialvar}), when $\phi _jdz^2 = (1/2)(M-iN)dz^2$.

If $U$ is an open ball in $\Sigma $ disjoint from the $\epsilon $-balls about the branch points, we can introduce coordinates $(u_1,u_2)$ just like did in the proof of Lemma~3.2 and if follows from (\ref{E:holocomparison}) that
$$ \int _U \langle \pi _V \circ (D_2F_0)_{(f,\omega ,g)} (\phi_j) , Z_k \rangle dA = 0, $$
so the only contributions to the integral $B\left(\phi _jdz^2, Z_k \right)$ come from the small neighborhoods of the branch points on which the derivative of the cutoff function $\zeta$ is nonzero.

We focus on a neighborhood of $p_j$ and write $\phi'_j = (1/2)(M'-iN')$ as in \S\ref{S:harmonic}.  In terms of the coordinate $z$ centered at $p_j$ on $\Sigma $ and the corresponding coordinate $w = u_1 + i u_2$ on $M$ which approximates $z^{\nu + 1}$, we can apply the conjugate of (\ref{E:partialwitoz}).  Since $\phi'dz^2$ is a meromorphic differential, we find that
\begin{multline*}\frac{\partial }{\partial \bar w} (M'-iN')\\
= \left[\left(\frac{1}{(\nu +1) \bar z^\nu} + o(|z|^{-\nu}) \right)\frac{\partial }{\partial \bar z}+ o(|z|^{-\nu + 1})\frac{\partial }{\partial z}\right](M'-iN') \\
= o(|z|^{-\nu - 1}),\end{multline*}
where the extra factor $|z|^{-2}$ comes from the fact that $(M'-iN') = b/z + o(|z|^{-1})$ for some complex constant $b$, while its derivative with respect to $z$ is $-b/z^2 + o(|z|^{-2})$.  Hence
$$\frac{1}{2}\left[ \left(\frac{\partial M'}{\partial u_1} + \frac{\partial N'}{\partial u_2}\right) - i \left(\frac{\partial N'}{\partial u_1} - \frac{\partial M'}{\partial u_2} \right)\right] = o(|z|^{-\nu - 1}),$$
except at the origin where the differential has a singularity.

But we are interested in the variation defined by $\phi = \zeta \phi' = (1/2)(M-iN) = (1/2)(\zeta M' - i\zeta N')$, and in this case the calculation yields
\begin{multline*} \frac{\partial }{\partial \bar w} (M-iN) = (M'-iN') \frac{1}{(\nu +1) \bar z^\nu }\frac{\partial \zeta }{\partial \bar z} + \zeta (r) o(|z|^{-\nu - 1})\\
= (M'-iN') \zeta'(r) \frac{1}{(\nu +1) \bar z^\nu}\frac{\partial r}{\partial \bar z} + \zeta (r) o(|z|^{-\nu - 1}).\end{multline*}
In terms of polar coordinates $z = r e^{i\theta }$,
\begin{multline*}\frac{\partial r }{\partial \bar z} = \frac{1}{2} \left(\frac{\partial }{\partial x_1} + i \frac{\partial }{\partial x_2} \right) \left(\sqrt{x_1^2+x_2^2} \right) \\
= \frac{1}{2} \left(\frac{x}{r} + i \frac{y}{r} \right) = \frac{1}{2}( \cos \theta + i \sin \theta ) = \frac{1}{2}e^{i\theta },\end{multline*}
while
$$(M'-iN') = b\frac{1}{z} + O(1) = \frac{be^{-i\theta }}{r} + O(1),$$
where $b$ is a nonzero complex constant.  Hence
$$\frac{1}{2}\left[ \left(\frac{\partial M}{\partial u_1} + \frac{\partial N}{\partial u_2}\right) - i \left(\frac{\partial N}{\partial u_1} - \frac{\partial M}{\partial u_2} \right)\right] = \frac{b \zeta'(r)}{2r} \frac{1}{(\nu +1) \bar z^\nu} + \zeta (r) o(|z|^{-\nu - 1}),$$
or equivalently,
$$(\dot \Gamma _{1,11} + \dot \Gamma _{1,22}) - i (\dot \Gamma _{2,11} + \dot \Gamma _{2,22}) = \frac{b \zeta'(r)}{r}\frac{1}{(\nu +1) \bar z^\nu} + \zeta (r) o(|z|^{-\nu - 1}).$$
Using the expression (\ref{E:expansionofZk}) for $Z_k$, we obtain
\begin{multline*} \left[ \sum _{a=1}^2 (\dot \Gamma _{1,aa} - i \dot \Gamma _{2,aa}) (h_{ki}(z)) \right] \\
= \frac{b \zeta'(r)}{r}\frac{1}{(\nu +1) \bar z^\nu } h_{ki}(z)) +  \zeta (r) o(|z|^{-2\nu - 1}) \\
= \frac{b \zeta'(r)}{r} \frac{1}{(\nu +1) \bar z^\nu } h_{ki}(z)) +  \zeta (r) o(|z|^{-2\nu - 1}).\end{multline*}
The real part of the last expression is the integrand in (\ref{E:derivofF}), so substitution into the complex form of (\ref{E:derivofF}) yields
\begin{multline*} \int _ {B_\epsilon (p_j)} \left\langle \pi _V \circ (D_2F)_{(f,\omega ,g)} (\phi_j), Z_k \right\rangle dA \\ = - (\nu + 1)\int _ {B_\epsilon (p_j)}
\left[ (\delta _{jk}b) \frac{\zeta'(r)}{r}r + \zeta (r) o(1) \right] dr d\theta \\ = - (\nu + 1)\int _ {B_\epsilon (p_j)} (\delta _{jk}b)\zeta'(r) dr d\theta + \hbox{Error}(\epsilon ), \end{multline*}
an expression linear in both $\phi $ and $Z$, where $\hbox{Error}(\epsilon ) \rightarrow 0$ as $\epsilon \rightarrow 0$.  Thus
\begin{multline*} \int _0^\epsilon \zeta'(r)dr = \zeta(\epsilon ) - \zeta (0) = 1 \quad \Rightarrow \\ \lim _{\epsilon \rightarrow 0}\int _ {B_\epsilon (p_j)} \left\langle \pi _V \circ (D_2F_0)_{(f,\omega ,g)} (\phi _j), Z_k \right\rangle dA = - 2 \pi \delta _{jk}(\nu + 1) b.\end{multline*}
Finally, we choose $b$ so that $2\pi (\nu + 1) b = - 1$.  A similar calculation shows that
$$\lim _{\epsilon \rightarrow 0}\int _ {B_\epsilon (p_i)} \left\langle \pi _V \circ (D_2F_0)_{(f,\omega ,g)} (\phi _j), Z_k \right\rangle dA = 0, \qquad \hbox{when $i \neq j$.}$$

These calculations show that our collection $Z_1, \ldots, Z_n$ of holomorphic sections of ${\bf L}$ and our collection $\phi_1dz^2, \ldots , \phi _ndz^2$ of quadratic differentials approximate (\ref{E:tobeproved}) better and better as $\epsilon \rightarrow 0$.  Choosing $\epsilon > 0$ sufficiently small and replacing $Z_1, \ldots, Z_n$ by linear combinations of $Z_1, \ldots, Z_n$ enables us to establish (\ref{E:tobeproved}) on the nose.

\vskip .1 in
\noindent
In the case that $\Sigma $ has genus zero or one, the preceding proof becomes obstructed because we may not be able to construct holomorphic differentials with poles at prescribed points.  Thus if $\Sigma = S^2$, we need to introduce three additional points $ \{ q,r,s \}$ not in the branch locus and allow the meromorphic differential $\phi_j'dz^2$ to have simple poles at these points; it is then possible to prescribe for a meromorphic quadratic differential a single additional simple pole at any particular choice of the branch point of $f$.  To compensate for this, we must replace each $Z_k$ by $Z'_k = Z_k + W_k$, where $W_k$ is an element of space ${\mathcal G}$ of infinitesimal symmetries corresponding to the $G$-action, chosen so that $Z'_k(q)$, $Z'_k(r)$ and $Z'_k(s)$ all vanish.  We multiply by cutoff functions near each of the points $q$, $r$ and $s$ and carry through the calculation we presented above.  Assuming that the meromorphic differential $\phi_j'dz^2$ is not regular, but has a simple pole at $q$, the above calculation yields a zero contribution at $q$ since $Z'_k(q) = 0$.  Since $Z'_k(r)$ and $Z'_k(s)$ also vanish, we find that the contributions at $r$ and $s$ are also zero.

When $\Sigma = T^2$ the argument is similar, except we need only let the meromorphic differential $\phi _j'dz^2$ have an additional simple pole at a single predetermined point $q$.  Once again, we replace each $Z_k$ by $Z'_k = Z_k + W_k$, where $W_k$ is an infinitesimal symmetry so that $Z'_k(q) = 0$, and find that this forces the contribution at $q$ to be zero.

\vskip .1 in
\noindent
Finally, we show that the hypothesis of Lemma~5.1.1 can be weakened to allow finitely many points of self-intersection.  To do this it suffices to show that we can replace $\phi dz^2$ by a metric variation which vanishes at all of the points of self-intersection:

\vskip .1 in
\noindent
{\bf Lemma~5.1.2.} {\sl Given a finite subset $F = \{ q_1, \ldots , q_N \} \subseteq \Sigma$, disjoint from the branch locus and a trace-free tangential variation $\phi dz^2$ in $T_g\hbox{Met}(M)$ which is bounded at each point of $F$, we can construct a family of smoooth functions $\psi _\epsilon : M \rightarrow [0,1]$ such that each $\psi _\epsilon \phi dz^2$ vanishes in a neighborhood of every point of $F$, and
$$\lim _{\epsilon \rightarrow 0} B_0 \left(\psi _\epsilon \phi dz^2, Z \right) = B_0 \left(\phi dz^2, Z\right),$$
for all $Z \in \Gamma ({\bf L})$.}

\vskip .1 in
\noindent
There is a standard cutoff useful for such purposes:  We first define a piecewise smooth map $\psi : {\mathbb R} \rightarrow {\mathbb R}$ by
$$\psi (r) = \begin{cases} 0, & \hbox{if $r \leq \epsilon ^2$,} \cr (\log (\epsilon ^2) - \log r)/(\log \epsilon ), & \hbox{if $\epsilon ^2 \leq r \leq \epsilon $,} \cr 1, & \hbox{if $\epsilon  \leq r$,} \end{cases}$$
so that
$$\frac{d\psi}{dr} (r) = \begin{cases} 0, & \hbox{if $r \leq \epsilon ^2$,} \cr (-1)/(r\log \epsilon ), & \hbox{if $\epsilon ^2 \leq r \leq \epsilon $,} \cr 0, & \hbox{if $\epsilon  \leq r$,} \end{cases}$$
and
$$\int _0^{2\pi}\int _0^\epsilon \left| \frac{d\psi}{dr} (r) \right| rdrd\theta = \int _{\epsilon ^2}^\epsilon \frac{2\pi}{(\log \epsilon )}dr = \frac{2 \pi(\epsilon - \epsilon ^2)}{\log \epsilon },$$
which approaches zero as $\epsilon \rightarrow 0$.  Thus if we define $\psi _i : \Sigma \rightarrow {\mathbb R}$ so that it is one outside an $\epsilon$-neightborhood of the point $q_i \in F$ and in terms of polar coordinates $(r_i, \theta _i)$ about the point $q_i$ is obtained from $\psi _i = \psi \circ r_i$ by smoothing the corners, and let $\psi _ \epsilon = \psi _1 \cdots \psi _N$, we obtain a cutoff function which vanishes at every element of $F$.  The family $\{ \psi _\epsilon : 0 < \epsilon \leq \epsilon _0 \}$ then satisfies the above limit.

\subsection{Bumpy metrics for simple branch points}
\label{S:simple}

We now have the means to prove the Main~Theorem for prime minimal surfaces with at worst only simple branch points and with self-intersection set consisting of finitely many points.  This will follow from Proposition~4.4.1 and

\vskip .1 in
\noindent
{\bf Proposition~5.2.1.} {\sl Suppose that $M$ is a compact connected manifold of dimension at least four.  For a generic choice of Riemannian metric on $M$, there are no prime minimal surfaces with finite self-intersection set, that have nontrivial branch points, all of which are simple.}

\vskip .1 in
\noindent
To prove this, we replace the space ${\mathcal P}^\emptyset $ defined by (\ref{E:mathcalSprime}), and used in the proof of Proposition~4.4.1, by
\begin{multline}{\mathcal P}^s = \{ (f,\omega ,g) \in \hbox{Map}(\Sigma ,M) \times {\mathcal T} \times \hbox{Met}(M) : \\
\hbox{$f$ is a prime conformal $\omega $-harmonic map } \\
\hbox{for $g$ which has only simple branch points} \}, \label{eq:mathcalPprime} \end{multline}
an open subset of the space of all prime harmonic maps.  As in \S \ref{S:minimalnobranch}, we break the $G$-symmetry, and study ${\mathcal P}^s$ by means of the Euler-Lagrange map
$$F : \hbox{Map}(\Sigma ,M) \times {\mathcal T} \times \hbox{Met}(M) \longrightarrow \hbox{Map}(\Sigma ,TM) \times T{\mathcal T},$$
which we regard as a section of the vector bundle
\begin{multline*} {\mathcal E} = \{ (f,\omega ,g,X,\dot \omega) \in \hbox{Map}(\Sigma ,M) \times {\mathcal T} \times \hbox{Met}(M) \times \hbox{Map}(\Sigma ,TM) \times T{\mathcal T} \\ \hbox{such that $\pi (X) = f$ and $\pi (\dot \omega ) = \omega $ }  \}.\end{multline*}
The space ${\mathcal P}^s$ consists of critical points for the two-variable energy $E$, which are exactly the points $(f,\omega ,g)$ such that $F(f,\omega ,g) = 0$.

Once again, if $(f,\omega, g)$ is a zero for $F$, we let
\begin{multline} {\mathcal V} = \{ (X, \dot \omega ) \in T_f\hbox{Map}(\Sigma ,M) \oplus T_\omega {\mathcal T} : \hbox{ $(X, \dot \omega )$ is perpendicular}\\
\hbox{ to  the image of $\pi _V \circ DF_{(f,\omega ,g)}$} \}, \label{E:defofVminimal5} \end{multline}
where $\pi _V$ is projection on the fiber, and the main step in the argument is to show that ${\mathcal V} = {\mathcal G}(f,\omega ,g)$ where ${\mathcal G}(f,\omega ,g)$ is the space generated by the $G$-action:

\vskip .1 in
\noindent
{\bf Lemma~5.2.2.}  {\sl If $(f, \omega )$ is a critical point for $E$ such that $f$ is somewhere injective and has only simple branch points, the only elements of $T_{(f,\omega )}( \hbox{Map}(\Sigma ,M) \times {\mathcal T})$ that can be perpendicular to the image of every variation in the metric are the elements of ${\mathcal G}(f,\omega ,g)$.}

\vskip .1 in
\noindent
Once we have this Lemma, we can conclude that
$$(\hbox{image of $DF_{(f,\omega ,g)}$}) + {\mathcal G}(f,\omega ,g) + T_{(f,\omega ,g)}{\mathcal Z} = T_{(f,\omega, g,0,0)}{\mathcal E},$$
and it will follow from the implicit function theorem that
\begin{multline}{\mathcal S} = \{ (f,\omega ,g) \in {\mathfrak M}(\Sigma ,M) \times {\mathcal T} \times \hbox{Met}(M) : \hbox{$f$ is a prime conformal} \\ \hbox{$\omega $-harmonic immersion for $g$ with only simple branch points} \} \label{E:scriSminsimplebranch}\end{multline}
is a submanifold.

\vskip .1 in
\noindent
To prove Lemma~5.2.2, we first note that by Lemma~4.4.2, the only elements $(X,\dot \omega )$ of ${\mathcal V}$ are tangential Jacobi fields $(X , \dot \omega)$ which are the real parts of elements $(Z, \bar \psi d \bar z^2)$ of ${\mathcal J}({\bf L})$.  In the special case where $\Sigma $ has genus zero, Lemma~5.1.1 (and Lemma~5.1.2 to handle self-intersections) implies that all such Jacobi fields can be covered by perturbations in the metric, and Lemma~5.2.2 is therefore established in this case.

To treat the higher genus case, we must allow for deformations in the complex structure on $\Sigma $, and for this, we extend the functional $B_0$ introduced in \S \ref{S:minimalnobranch} to a two component functional
$$B : (\hbox{trace-free tangential variations in $T_g\hbox{Met}(M)$}) \times \Gamma({\bf L}) \times \overline {\mathcal Q}(\Sigma ) \longrightarrow {\mathbb C}.$$
where $\overline {\mathcal Q}(\Sigma )$ is the space of quadratic differentials $\bar \eta d \bar z^2$, defined by
\begin{equation} B \left(\phi dz^2, (Z, \bar \eta d \bar z^2) \right) = \int _\Sigma \left\langle \pi _V \circ (D_2F)_{(f,g)} (\dot g), Z \right\rangle dA + \frac{1}{8} \int_\Sigma \frac{1}{\lambda ^4} \phi \bar \eta dA. \label{EdefinitionB} \end{equation}
We claim that we can in fact, enhance the conclusion of Lemma~5.1.1 to the following:  there is a collection $\phi _1dz^2, \ldots , \phi _ndz^2$ of quadratic differentials on $\Sigma $, vanishing near the branch points, such that with this new definition of $B$, we have
\begin{equation} B\left(\phi _jdz^2, Z_k \right) = \delta _{jk} = \begin{cases} 1, & \hbox{if $j=k$},\cr 0, & \hbox{if $j \neq k$}. \end{cases}\label{E:tobeproved1}\end{equation}

Indeed, there was some freedom in constructing the quadratic differential $\phi _jdz^2$.  Recall that it was constructed from a meromorphic differential $\phi _j'dz^2$with a simple pole at $p_j$ and no other poles.  To $\phi _j'$ we could have added an arbitrary holomorphic quadratic differential without effecting the first term $B_0$.  Indeed, suppose that Teichm\"uller space ${\mathcal T}$ has dimension $N$ and that
$$\eta _1 dz^2 , \eta _2 dz^2, \ldots , \eta _Ndz^2$$
is a basis for the holomorphic quadratic differentials at a point $\omega \in {\mathcal T}$.  Then if $\phi _jdz^2$ is obtained from a meromorphic differential with a single simple pole at $p_j$ and no other singularities by the cutoff proces we described, there exist unique constants $c_{ji} \in {\mathbb C}$ such that
$$\left\langle \phi _jdz^2 - \sum _{i = 1}^N \eta _i dz^2, \bar \eta _k d \bar z^2 \right\rangle = 0, \quad \hbox{for $1 \leq k \leq N$,}$$
where $\langle \cdot , \cdot \rangle$ is the complex bilinear form defined by (\ref{E:hiponquad}), which also appears in the second term of (\ref{EdefinitionB}).  This enables us to achieve (\ref{E:tobeproved1}), which is sufficient for proving Lemma~5.2.2.  

\vskip .1 in
\noindent
Once we know that the subspace ${\mathcal S}$ defined by (\ref{E:scriSminsimplebranch}) is a submanifold, the argument proceeds exactly as before.  The next Lemma is proven in exactly the same way as Lemma~3.5 and Lemma~4.4.3:

\vskip .1 in
\noindent
{\bf Lemma~5.2.3.} {\sl The projection $\pi : {\mathcal S} \rightarrow \hbox{Met}(M)$ is a Fredholm map of Fredholm index $\dim G$.}

\vskip .1 in
\noindent
To finish the proof of Proposition~5.2.1, we note that if $g$ admits parametrized minimal surfaces with a nontrivial divisor of simple branch points, then $g$ cannot be a regular value of the projection.  Therefore such minimal surfaces cannot occur for generic choice of metric.

\section{Self-intersections}
\label{S:selfintersection}

As mentioned before, prime parametrized minimal surfaces are automatically somewhere injective.  However, we will need a much more explicit description of the set of self-intersections, and for this we need to understand the behavior of minimal surfaces near branch points.

An explicit normal form near a branch point is given by Theorem~1.4 of Micallef and White \cite{MW} which states that if $p$ is a branch point of a weakly conformal harmonic map $f : \Sigma \rightarrow M$, then there are normal coordinates $(u_1, \ldots , u_n)$ centered at $f(p)$ on the ambient Riemannian manifold $M$ and a $C^{2,\alpha }$ diffeomorphism $z \mapsto \xi(z)$ near $p$ on $\Sigma $, $z$ being a holomorphic coordinate centered at $p$, such that $\xi(z) = c z + O(|z|^2)$ for some nonzero complex constant $c$, and if we abbreviate the composition $u_i \circ f \circ \xi$ to $u_i$, then
\begin{equation} u_1(z ) + i u_2(z )= \xi (z)^{\nu + 1}, \quad u_r(z) = f_r(\xi(z)) \quad \hbox{for} \quad 3 \leq r \leq n, \label{eq:coordinates} \end{equation}
where each $f_r(\xi )$ is $o_1(|\xi |^{\nu +1})$.  Here $o_1(z^{k+1})$ stands for a term which is $o(z^{k+1})$, and which has a derivative that is $o(z^{k})$.  Thus in terms of the nonconformal coordinate $\xi $, the error term in $u_1 + i u_2$ vanishes.  If we use $\xi$ as the local (nonconformal) parameter on $\Sigma $, we can write (\ref{eq:coordinates}) as
\begin{equation} \xi \mapsto (\xi ^{\nu + 1}, \tilde f(\xi)), \qquad \hbox{where $\tilde f$ takes values in ${\mathbb R}^{n-2}$}. \label{eq:canonicalformbranch1} \end{equation}
Theorem~1.4 of \cite{MW} goes on to state that if $\mu ^{\nu +1} = 1 $, then either $\tilde f(\mu \xi) - \tilde f(\xi)$ is identically zero, or
\begin{equation}\tilde f(\mu \xi) - \tilde f(\xi) = a \xi ^\rho + \bar a \bar \xi ^\rho + o_1(|\xi |^\rho), \label{eq:canonicalformbranch2} \end{equation}
where $a$ is a nonzero element of ${\mathbb C}^{n-2}$ and $\rho > \nu + 1$.

We say that a branch point $p$ for a harmonic map $f$ is {\em ramified\/} of order $\sigma $ if $\sigma + 1$ is the largest integer such that there is an open neighborhood $U$ of $p$ and a smooth map $\psi : U \rightarrow D$, where $D$ is the unit disk centered at $0$ in the complex plane, such that
\begin{enumerate}
\item $\psi (p) = 0$,
\item $\psi |(U - \{ p \})$ is a $(\sigma + 1)$-sheeted cover of $D - \{ 0 \}$, and
\item if $q,r \in U - \{ p \}$, then $f(q) = f(r) \Leftrightarrow \psi (q) = \psi (r)$.
\end{enumerate}
If $\sigma = \nu$, the branching order of $p$, we say that $p$ is a {\em false branch point\/}.  At the other extreme, we say that the branch point is {\em primitive\/} or {\em unramified\/} if the largest such integer $\sigma $ is zero.

Note that a branch point is primitive precisely when $\tilde f(\mu \xi) - \tilde f(\xi)$ is not identically zero for any nontrivial $(\nu + 1)$-fold root of unity.

\vskip .1 in
\noindent
{\bf Lemma~6.1.} {\sl If $f:\Sigma \rightarrow M$ is a prime parametrized minimal surface, all its branch points are primitive.}

\vskip .1 in
\noindent
This follows quickly from the Aronszajn unique continuation theorem and the theory of branched immersions presented by Gulliver, Osserman and Royden in \S 3 of \cite{GOR}.  For the convenience of the reader we briefly sketch the idea behind the argument, referring to \cite{GOR} for more details.

We begin by defining an equivalence relation $\sim $ on points of $\Sigma $ by setting $p \sim q$ if there are open neighborhoods $U_p$ and $U_q$ of $p$ and $q$ respectively, and a conformal or anti-conformal diffeomorphism $\psi :U_p \rightarrow U_q$ such that $f \circ \psi = f$.  (We make a minor modification to the construction in \S 3 of \cite{GOR} by allowing $\psi $ to be orientation-reversing so that we can allow for branched covers of nonorientable surfaces.)  Using the argument in \cite{GOR}, which is based upon Aronsjazn's unique continuation theorem, one shows that $\sim $ is indeed an equivalence relation and that the quotient space $\Sigma _0$ is a topological manifold (possibly nonorientable) which is smooth except at the branch points.  We can define $f_0: \Sigma _0 \rightarrow M$ by $f_0([p]) = f(p)$, where $[p]$ denotes the equivalence class of $p$, so that if $\pi : \Sigma \rightarrow \Sigma _0$ is the quotient map, $f_0 \circ \pi = f$.  We note that any point equivalent to a branch point is itself a branch point.  The restriction of $f_0$ to $\Sigma _0$ minus the equivalence classes of the branch points is a harmonic map of finite energy.  It therefore follows from the removeable singularity theorem of Sacks and Uhlenbeck (Theorem~3.6 of \cite{SU1}) that the restriction of $f_0$ can be extended to the equivalence classes of the branch points so as to be a harmonic map, which is conformal so long as $f$ is conformal.

Thus $f$ is a branched cover of $f_0$.  The hypothesis that $f$ is prime implies that $\Sigma _0$ must equal $\Sigma $, $\pi $ is the identity, and every branch point must be primitive, finishing our sketch of the proof.

\vskip .1 in
\noindent
A description of the intersection of two minimal surfaces was an important ingredient in the highly successful applications of minimal surface theory to the topology of three-manifolds developed by Meeks and Yau and others (see \cite{MY2} and \cite{FHS} for example), and we claim that a similar description can be established in the case in which the ambient space $M$ has arbitrary dimension $\geq 3$, using Lemma~4.1, the canonical form in a neighborhood of a branch point, and the following Lemma, a special case of Lemma~2.4 of Cheng \cite{Ch}.

\vskip .1 in
\noindent
{\bf Lemma~6.2.} {\sl Suppose that $f$ is a smooth real-valued function on an open neighborhood of the origin in ${\mathbb R}^2$ and $p$ is a nonzero homogeneous harmonic polynomial of degree $\rho \geq 2$ on ${\mathbb R}^2$ such that
$$ f(\xi ) - p(\xi ) = o_1(\xi |^\rho ), \qquad \hbox{for some positive integer $\rho $.} $$
Then there is a local $C^1$ diffeomorphism $\Phi $ of ${\mathbb R}^2$, fixing the origin, such that $f(\xi ) = p(\Phi (\xi ))$.}

\vskip .1 in
\noindent
We refer to \cite{Ch}, pages 48 and 49, for an elegant proof of this lemma.

\vskip .1 in
\noindent
With Lemmas~6.1 and 6.2 at our disposal, we are ready to present the main result of this section.

\vskip .1 in
\noindent
{\bf Lemma~6.3.} {\sl Let $f :\Sigma \rightarrow M$ is a prime compact connected parametrized minimal surface in a compact Riemannian manifold of dimension at least three, with self-intersection set
$$K = \{ p \in \Sigma : \hbox{ $f^{-1}(f(p))$ contains more than one point } \},$$
a compact subset of $\Sigma $.  Then each point within $K$ belongs to one of two classes:
\begin{enumerate}
\item an isolated point of self-intersection (of which there can be only finitely many),
\item a point which lies in a neighborhood $U \subset \Sigma $ such that $K \cap U$ is a subset of a finite number of $C^1$ curve segments of finite length in $U$ intersecting at the point.
\end{enumerate}}
 
\vskip .1 in
\noindent
{\bf Remarks.}  Since $\Sigma $ is compact, the Lemma implies that the self-intersection set $K$ itself is contained in a finite number of points and a finite number of $C^1$ curves of finite length.  If $M$ has dimension three, Lemma~4.3 is a direct consequence of Lemma~2 in \cite{MY2} or Lemma~1.4 in \cite{FHS}, and in fact the argument we present for Lemma~4.3 is just a modification of the argument given in \cite{MY2}.
\vskip .1 in
\noindent
To begin the proof, we first note that if $p \in K$, there are only finitely many points in $f^{-1}(f(p))$, or in other words only finitely many sheets of $f(\Sigma )$ which pass through $f(p)$.  Indeed, if $f$ is minimal, the intersection of each sheet with a small geodesic ball of fixed radius around $f(p)$ has area bounded below by a fixed $\epsilon > 0$ by (for example) Lemma~1, page 445 of \cite{MY1}, implying there can be only finitely many such sheets.

There is a simple idea which can be used to construct such an estimate, particularly if $f$ is an immersion.  The key point is that (as mentioned in \S \ref{S:prelim}) it follows directly from the Gauss equation for $f$ that the curvature of $f(\Sigma )$ in the induced metric is $\leq K_0$, where $K_0$ is the maximum value of the sectional curvatures of $M$.  Moreover, a geodesic ball $B$ of radius $\delta > 0$ about $f(p)$ in $M$ contains a geodesic ball of radius $\delta $ about $f(p)$ in $f(\Sigma )$.  By using Jacobi fields along geodesics to compare with a space form of constant curvature $K_0$, one obtains a lower bound on the area of one sheet of $f(\Sigma )\cap B$ passing through $f(p)$.  The case where $f$ is not an immersion is only a little more complicated; in this case, it is still possible to define geodesics emanating from $p$ (even though the induced metric is degenerate at $p$) and carry out the Jacobi field estimates.  Thus a bound on the energy implies that there are only finitely many sheets passing through $f(p)$.

By an easy induction, we see that it suffices to consider the self-intersection set of two harmonic maps $f_1 : D_1 \rightarrow M$ and $f_2 : D_2 \rightarrow M$ from disks centered at $p_1$ and $p_2$ respectively such that $f_1(p_1) = f_2(p_2) = q$, as well as the self-intersection of a single such map near a branch point.

We start with the first case under the assumption that neither $p_1$ nor $p_2$ is a branch point.  If $K$ is the intersection set of $f_1$ and $f_2$, we need to show that any point in $K$ lies in one of the two classes in the conclusion of the lemma.  The tangent planes $\Pi _1$ and $\Pi _2$ at $q$ to $f_1(D_1)$ and $f_2(D_2)$, respectively, can intersect in a point or a line, or the two planes may coincide.  If the planes intersect in a point, then $p_1 $ and $p_2$ belong to the first class in the conclusion of the lemma.

If $\Pi _1$ and $\Pi _2$ intersect in a line, the planes make an angle $\theta $ with each other, where $0 < \theta < \pi $.  In this case there is a plane $\Pi$ which makes an acute angle with both $\Pi _1$ and $\Pi _2$.  On the other hand, if $\Pi _1 = \Pi _2$, we simply set $\Pi = \Pi_1 = \Pi_2$.  In either case, each of the harmonic surfaces can be regarded as a graph over $\Pi$ near $q$, in terms of normal coordinates on $M$ centered at $q$.

In more detail, we choose normal coordinates $(u_1, \ldots , u_n)$ on an open neighborhood of $q$ in $M$ with $u_i(q) = 0$, so that the plane spanned by the coordinate vectors $(\partial /\partial u_1)$ and $(\partial /\partial u_2)$ at $q$ is $\Pi$.  Let $B_\epsilon $ be a ball of radius $\epsilon $ about $0$ in the $(u_1,u_2)$-plane.  We can consider the two surfaces as graphs of functions
$$h_i : B_\epsilon \rightarrow {\mathbb R}^{n-2}, \qquad i = 1,2,$$
which satisfy an explicit nonlinear partial differential equation.  Our convention is to denote the Riemannian metric on $M$ and the metric on $B_\epsilon $ induced by the immersion $(u_1,u_2) \mapsto (u_1,u_2,h_i(u_1,u_2))$ by
$$\sum _{i,j=1}^n g_{ij} du_idu_j \qquad \hbox{and} \qquad \sum _{a,b=1}^2 \eta_{ab} du_adu_a,$$
respectively.  With this convention, it follows from a relatively straightforward calculation (presented in the Appendix to \cite{MW}) that each $h_i$ is an ${\mathbb R}^{n-2}$-valued solution to an elliptic system in \lq\lq nonparametric form,"
$$ \sum _{a,b=1}^2 \eta^{ab}(u_1,u_2,h,Dh) \frac{\partial ^2 h}{\partial u_a \partial u_b} + \phi (u_1,u_2,h,Dh) = 0,$$
where $(\eta^{ab})$ is the matrix inverse of $(\eta _{ab})$, the operator on the left once again having scalar symbol.

It is readily checked that the difference $h = h_2-h_1$ satisfies an equation of the same form,
$$ \sum _{a,b=1}^2 \eta^{ab}(u_1,u_2,h,Dh) \frac{\partial ^2 h}{\partial u_a \partial u_b} + \phi (u_1,u_2,h,Dh) = 0,$$
in which the coefficients depend upon $h_1$, which is now regarded as a known function of $u_1$ and $u_2$, with $h$ or equivalently $h_2$ being regarded as unknown.  Thus by Corollary~1 of Theorem~1.1 in \cite{MW} (a corollary of a theorem of Hartman and Wintner), we see that either $h \equiv 0$ or
\begin{multline*} |h(u_1,u_2)| = O(|(u_1,u_2)|^m) \quad \hbox{for some positive integer $m$, and}  \\
\Rightarrow \qquad h(u_1,u_2) = p(u_1,u_2) + o_1(|(u_1,u_2)|^m),\end{multline*}
where $p(u_1,u_2)$ is a nonzero homogeneous vector-valued harmonic polynomial of degree $m$.  We recall that stating that $h - p$ is $o_1(|z|^m)$ means that $h - p$ is $o(|z|^m)$ and $D(h - p)$ is $o(|z|^{m-1})$.  It cannot be the case that $h \equiv 0$ by Lemma~4.1, since $f$ is assumed to be prime.

Since the domain is only two-dimensional, such harmonic polynomials are easily described; indeed,
$$p(u_1,u_2) = \hbox{Re} (a(u_1 + i u_2)^m), \qquad \hbox{where $a$ is a nonzero element of ${\mathbb C}^{n-2}$.}$$
If $m = 1$, the two surfaces intersect in a smooth curve near $q$.  If $m  > 1$, this vector-valued harmonic polynomial vanishes at a single point or on a $C^1$ collection of curves.  In the former case, the two surfaces intersect in an isolated point and we are done.  In the latter case, we can choose a one-dimensional linear subspace $V$ of ${\mathbb R}^{n-2}$ such that if $\pi_V$ denotes the orthogonal projection into this subspace, $\pi_V \circ p$ is not identically zero.  The zero set of the ${\mathbb R}^{n-2}$-valued function $h$ is contained in the zero set of the scalar function $\pi _V \circ h$, so it suffices to show that the zero set of $\pi _V \circ h$ is a finite collection of $C^1$ curves passing through the origin.

But now, following \cite{MY2}, we can apply Lemma~4.2 to conclude that there is a $C^1$ local diffeomorphism $\Phi $ near the origin in ${\mathbb R}^2$ such that $\pi_V \circ h = \pi _V \circ p \circ \Phi$.  Thus $\pi _V \circ h$ vanishes on the intersection of a finite number of $C^1$ curves in terms of the parameter $\xi $, and by the canonical form described at the beginning of the section, also in terms of the conformal parameter $z$.  This establishes the desired conclusion in this case.

To finish the proof in the case that branch points are present, recall that we need to consider:
\begin{enumerate}
\item the self-intersection set of a single small minimal disk $f : D \rightarrow M$ in which $D$ is centered at a branch point $p$, and
\item the intersection of two minimal disks $f_1 : D_1 \rightarrow M$ and $f_2 : D_2 \rightarrow M$ centered at $p_1$ and $p_2$, with $f(p_1) = f(p_2) = q$, in which one or both of the points $p_1$ and $p_2$ are branch points.
\end{enumerate}

If we were to make the assumption that $f$ is not only minimal, but also locally area minimizing, we could apply Theorem~4.2 of \cite{MW} to show that the self-intersection set of the restriction of $f$ to a neighborhood of the branch point $p$ is empty when the neighborhood is sufficiently small, completing the first of these steps.

In general, the self-intersections near a branch point $p$ on a single minimal disk $f : D \rightarrow M$ can be treated by the canonical form (\ref{eq:coordinates}).  We note that for particular points $\xi_1$ and $\xi_2$ of $D$, $f(\xi _1)$ can equal $f(\xi _2)$ only if $\xi_2 = \mu \xi_1$, where $\mu $ is a $(\nu + 1)$-fold root of unity, and $\tilde f(\mu \xi_1) = \tilde f(\xi_2)$, where $\tilde f$ is defined by (\ref{eq:canonicalformbranch1}).  According to Lemma~4.1, since $f$ is prime, it is primitive, and hence $\tilde f(\mu \xi) - \tilde f(\xi)$ is not identically zero for any nontrivial $(\nu + 1)$-root of unity.  Thus we can apply (\ref{eq:canonicalformbranch2}) to conclude that
$$\pi _V \circ \tilde f(\mu \xi) - \pi _V \circ \tilde f(\xi) = \pi _V(a \xi ^\rho + \bar a \bar \xi ^\rho) + o_1(|\xi |^\rho),$$
where $\pi _V $ is the projection into a three-dimensional subspace of ${\mathbb R}^n$ containing the tangent plane at  $p$ such that $\pi _V(a + \bar a) \neq 0$.  Applying Lemma~4.2 once again shows that $\pi _V \circ \tilde f(\mu \xi) - \pi _V \circ \tilde f(\xi)$ vanishes on a finite number of $C^1$ curves (in terms of the nonconformal parameter $\xi $) passing through $p$.  Once again, since the transformation from the conformal parameter $z$ to $\xi $ is $C^{2,\alpha }$, the self-intersection set of $f|D$ lies within a finite collection of $C^1 $ curves of finite length, finishing the argument in this case.

We next consider the case in which two minimal disks $f_1 : D_1 \rightarrow M$ and $f_2 : D_2 \rightarrow M$ intersect at the centers $p_1$ and $p_2$ of the disks.  If the tangent planes $\Pi _1$ and $\Pi _2$ to $f_1(D_1)$ and $f_2(D_2)$ at $ q = f(p_1) = f(p_2)$ intersect at a point, then the two minimal disks also intersect at a point, and $p_1$ and $p_2$ are isolated points of self-intersection.

Suppose next that the two planes $\Pi _1$ and $\Pi _2$ coincide and $p_1$ is a branch point of order $\nu$, while $p_2$ is not a branch point.  In this case, we replace $f_2 : D_2 \rightarrow M$ by a new disk $f_2' : D_2' \rightarrow M$ centered at $p_2'$ such that $f'_2(z) = f_2(z^{\nu + 1})$.  Thus $f_2'$ is a branched cover of $f_2$, and $f_2'$  has a false branch point at its center $p_2'$ with exactly the same order as the branch point at $p_1$.  Both $f_1$ and $f_2' $ have nonconformal parametrizations
$$ \xi \mapsto (\xi ^{\nu + 1}, \tilde f_1(\xi)), \qquad  \xi \mapsto (\xi ^{\nu + 1}, \tilde f_2'(\xi)), $$
in which $\xi $ is related to conformal parameters $z_1$ on $D_1$ and $z_2'$ on $D_2'$ by $C^{2,\alpha }$ diffeomorphisms.  It cannot be the case that $\tilde f_1 - \tilde f_2'$ is identically zero, because this would contradict unique continuation, since $f$ is prime.  Therefore it follows from Remark~1.6 of \cite{MW} that if $\mu$ is any nontrivial $(\nu + 1)$-fold root of unity,
$$\tilde f_1(\mu \xi) - \tilde f_2'(\xi) = a \xi ^\rho + \bar a \bar \xi ^\rho + o_1(|\xi |^\rho), $$
where $a$ is a nonzero element of ${\mathbb C}^{n-2}$ and $\rho > \nu + 1$.  An application of Lemma~4.2 again shows that $f_1$ and $f_2'$ intersect in a finite number of $C^1$ curves of finite length intersecting at $p_1$ or $p_2'$, in terms of the parameter $\xi $.  This holds also in terms of conformal parameters $z_1$ for $f_1$ and $z_2'$ for $f_2'$, just as before.  The images of these curves in the disk $D_2$ are $C^1$ curves intersecting at $p_2$.

If $\Pi _1$ and $\Pi _2$ coincide and both $p_1$ and $p_2$ are branch points of branching order $\nu_1$ and $\nu_2$, respectively, we let $\nu + 1$ be the lowest common multiple of $\nu _1 + 1$ and $\nu _2 + 1$.  We choose positive integers $\sigma _1$ and $\sigma _2$ so that $\nu + 1 = \sigma_1 (\nu_1 + 1)$ and $\nu + 1 = \sigma_2 (\nu_2 + 1)$.  For $i = 1,2$, we then replace $f_i : D_i \rightarrow M$ by a new disk $f_i' : D_i' \rightarrow M$ such that $f_i'(z) = f_i(z^{\sigma _i})$.  The two disks $f_1'$ and $f_2'$ are no longer necessarily primitive, but we can still apply Remark~1.6 of \cite{MW} to show that if $\mu$ is any $(\nu + 1)$-root of unity, either
$$\tilde f_1'(\mu \xi) - \tilde f_2'(\xi) \equiv 0 \quad \hbox{or} \quad \tilde f_1'(\mu \xi) - \tilde f_2'(\xi) = a \xi ^\rho + \bar a \bar \xi ^\rho + o_1(|\xi |^\rho), $$
where $a$ is a nonzero element of ${\mathbb C}^{n-2}$ and $\rho > \nu + 1$.  Once again the first case would contradict unique continuation, since $f$ is prime, so we can apply Lemma~4.2 yet again to conclude that $f_1'$ and $f_2'$ intersect in a finite collection of $C^1$ curves of finite length, and once again we argue that $f_1$ and $f_2$ themselves intersect in a finite number of $C^1$ curves of finite length.

Finally, we need to consider the case in which the two planes $\Pi _1$ and $\Pi _2$ intersect in a line.  As before, if $p_1$ is a branch point of order $\nu$, while $p_2$ is not a branch point, we replace $f_2$ by a branched cover $f_2'$ so that $f_1'$ and $f_2'$ both have branch points of order $\nu$.  If $p_1$ and $p_2$ are branch points of order $\nu_1$ and $\nu_2$ we let $\nu +1$ be the lowest common multiple of $\nu_1 + 1$ and $\nu_2 + 1$, and replace both $f_1$ and $f_2$ by branched covers $f_1'$ and $f_2'$ so that $f_1'$ and $f_2'$ both have branching order $\nu $.

Each of the maps $f_1'$ and $f_2'$ has a canonical form (\ref{E:canonicalformbr}), but since the tangent planes make an angle with each other, the normal coordinates in $M$ that are used to realize these forms are also at an angle.  Thus we can choose two normal coordinate systems
$$(u_1, u_2, u_3, u_4, \ldots , u_n) \quad \hbox{and} \quad  (u_1, v_2, v_3, u_4, \ldots , u_n)$$
related to each other by a rotation,
\begin{equation} u_2 = a v_2 - b v_3, \quad u_3 = b v_2 + a v_3, \quad a^2 + b^2 = 1, \label{eq:transformationbetnorm}\end{equation}
where $b \neq 0$, so that by (\ref{eq:coordinates}),
\begin{multline*} u_1(\xi _1 ) + i u_2( \xi_i )= \xi _1^{\nu + 1}, \quad u_r(\xi_1 ) = o_1(|\xi _1|^{\nu + 1}) \quad \hbox{for} \quad 3 \leq r \leq n, \\
u_1(\xi _2 ) + i v_2( \xi_2 )= \xi _2^{\nu + 1}, \quad v_3(\xi_2 ) = o_1(|\xi _2|^{\nu + 1}), \quad u_s(\xi_2 ) = o_1(|\xi _2|^{\nu + 1}) \\ \hbox{for} \quad 4 \leq s \leq n, \end{multline*}
where $\xi _1 $ and $\xi _2$ are nonconformal parameters on $D_1$ and $D_2$ respectively.  Note that $f_1'$ and $f_2'$ can intersect only where $\hbox{Re}(\xi_1) = \hbox{Re}(\xi_2)$.  Since both $\xi _1 $ and $\xi _2$ give Euclidean parametrizations of $\Pi _1$ and $\Pi _2$ respectively, we can arrange that $\xi_1$ and $\xi_2$ are related by a local diffeomorphism from a neighborhood of $0$ in $D_1$ to a neighborhood of $0$ in $D_2$ under which $\xi_1 \mapsto \xi _2$ with $u_1(\xi _1) = u_1(\xi _2)$.  

It follows from (\ref{eq:transformationbetnorm}) that
$$u_3(\xi _2) = b \ \hbox{Im}(\xi _1^{\nu + 1}) +  o_1(|\xi _1|^{\nu + 1}) = b \ \hbox{Im}(\xi _2^{\nu + 1}) +  o_1(|\xi _2|^{\nu + 1}),$$
$$v_3(\xi _1) = - b \ \hbox{Im}(\xi _2^{\nu + 1}) +  o_1(|\xi _2|^{\nu + 1}) = - b \ \hbox{Im}(\xi _1^{\nu + 1}) +  o_1(|\xi _1|^{\nu + 1}).$$
According to Lemma~4.2, we can find local $C^1$ diffeomorphisms $\Phi _1$ and $\Phi _2$ of $D_1$ and $D_2$ so that if $\Phi _i (\eta _i) = \xi _i$,
$$v_3 \circ \Phi _1(\eta_1) = - b \ \hbox{Im}(\eta _1^{\nu + 1}) \quad u_3 \circ \Phi _2(\eta_2) = b\ \hbox{Im}(\eta _2^{\nu + 1}).$$
The set of points in $D_1$ taken by $f_1' : D_1 \rightarrow M$ to points of $f_2'(D_2)$ are points such that $\hbox{Im}(\eta _1^{\nu + 1}) = 0$, while the points of $D_2$ taken by $f_2' : D_2 \rightarrow M$ to points of $f_1'(D_1)$ are points such that $\hbox{Im}(\eta _2^{\nu + 1}) = 0$.  Thus we find that the points of $K$ contributed in a neighborhood of a self-intersection of two disks whose tangent planes intersect in a line lie on finitely many $C^1$ curves of finite length.

Thus we have covered all cases and Lemma~6.3 is proven.

\section{Higher order branch points}
\label{S:higherorderbranch}

In the case of higher order branch points, we only know how to perturb away part of the space of tangential Jacobi fields by variations in the metric.  To compensate for the additional tangential Jacobi fields, we cut down the space $\hbox{Map}(\Sigma ,M)$ by imposing conditions at the branch points which conformal harmonic maps of a given branch type will automatically satisfy.  In a rough sense, we want to project away the Jacobi fields we cannot perturb away, but then to restore the Fredholm index zero property we need to impose conditions on the maps at the branch points.

To impose these conditions at the branch points, it is convenient to modify the configuration space on which we work to
$${\mathcal A}^n(\Sigma ,M) = \tilde{\mathcal M}^n(\Sigma, M) \times \hbox{Met}(M),$$
where $n$ is the number of branch points considered, and
\begin{equation} \tilde{\mathcal M}^n(\Sigma, M) = \frac{\hbox{Map}(\Sigma, M) \times \left( \overbrace{\Sigma \times \cdots \times \Sigma }^n - \Delta^{(n)} \right) \times \hbox{Met}_0(\Sigma )}{\hbox{Diff}_0(\Sigma )}, \label{E:tildem}\end{equation}
the tilde indicating that we divided by the action of $\hbox{Diff}_0(\Sigma )$ rather than the full diffeomorphism group.  Here $\Delta^{(n)}$ once again denotes the fat diagonal consisting of $n$-tuples $(p_1, \ldots , p_n)$ for which no two of the elements are equal.  Of course, use of a section would enable us to write this space more simply as
$$\tilde{\mathcal M}^n(\Sigma, M) = \hbox{Map}(\Sigma, M) \times {\mathcal T}_{g,n},$$
but that would hide the explicit role of the branch points.  We are interested in investigating the subspace
\begin{multline}{\mathcal P}^\Lambda = \{ ([f,{\bf p}, \eta], g) \in {\mathcal A}^n(\Sigma ,M) \hbox{ such that $f$ is a prime } \\ \hbox{conformal harmonic map with respect to $\eta $ and $g$} \\ \hbox{of branch type $\Lambda $ with branch locus ${\bf p} = (p_1, \ldots ,p_n)$} \}, \label{E:plambda}\end{multline}
It would be nice if we could show that ${\mathcal P}^\Lambda $ were a submanifold of ${\mathcal A}^n(\Sigma ,M)$.  Instead, we will merely show that ${\mathcal P}^\Lambda $ lies in the union of a countable collection of submanifolds, each of which has a Fredholm projection to $\hbox{Met}(M)$ with Fredholm index zero.  We then use these submanifolds to show that for generic metrics on ambient spaces of dimension at least four, minimal surfaces within ${\mathcal P}^\Lambda$ have transverse points of self-intersection as their only singularities (except for the branch points themselves).  Once we know that the singularities are well-behaved, we can modify the definition of these submanifolds and apply Lemma~5.1.1 to show that the modified submanifolds have Fredholm projection with negative Fredholm index to the space of metrics, thereby proving the Main Theorem.  We proceed now to fill in the details of this outline.

\subsection{Stratification according to branch type}
\label{S:stratification}

By a {\em branch type\/} for total branching order $\nu $, we mean a sequence of positive integers $(\nu _1, \ldots , \nu _n)$, such that
$$\nu _1 \geq \nu _2 \geq \cdots \geq \nu _n \qquad \hbox{and} \qquad \sum_{i=1}^n \nu _i = \nu.$$
We sometimes denote the total branching order of $\Lambda $ by $|\Lambda |$ and note that there are only finitely many branch types $\Lambda $ of total branching order $|\Lambda | = \nu $.  We define a partial order on branch types by demanding that if $\Lambda = (\nu _1, \ldots , \nu _n)$ and $\Lambda '= (\nu _1', \ldots , \nu _{n'}')$, then
$$\Lambda  \leq \Lambda '\quad \Leftrightarrow \quad n \leq n' \quad \hbox{and} \quad \nu _i \leq \nu _i' \quad \hbox{for $1  \leq i \leq n$.} $$
We write $\Lambda  < \Lambda '$ if $\Lambda \leq \Lambda '$ and $\Lambda \neq \Lambda '$.  There is a unique branch type $\emptyset $ of total branching order $0$, and $\emptyset \leq \Lambda $ for every branch type $\Lambda $.  We say that a divisor $D$ on $\Sigma $ has branch type $\Lambda = (\nu _1, \ldots , \nu _n)$ if there is an ordered collection of distinct points $(p_1, \ldots , p_n)$ such that $D = \sum \nu _i p_i$, and that a harmonic map $f : \Sigma \rightarrow M$ has branch type $\Lambda$ if its divisor of branch points $D(f)$ has branch type $\Lambda $.

We would like to describe the space of maps within $\hbox{Map}(\Sigma, M)$ which have \lq\lq branch type" $\Lambda $, but it is easier to describe a corresponding space
$$\hbox{Map}^{\Lambda}(\Sigma, M) \subseteq \hbox{Map}(\Sigma, M) \times \left( \overbrace{\Sigma \times \cdots \times \Sigma }^n - \Delta^{(n)} \right).$$
As elements of $\hbox{Map}^{\Lambda}(\Sigma, M)$, we take
$$(f,(p_1, \ldots p_n)) \in \hbox{Map}(\Sigma, M) \times \left( \Sigma \times \cdots \times \Sigma - \Delta^{(n)} \right)$$
such that $f$ is an immersion except at each point $p_i$ where the $\nu _i$-jet of $f$ vanishes, while the $(\nu_i + 1)$-jet does not.  (By saying that the $\nu $-jet of $f$ vanishes at $p$, we mean that all covariant derivatives of $f$ of order between one and $\nu$ vanish at $p$.)  Using the Sobolev imbedding theorem and the smoothness of the $\nu _i$-jet evaluation map described in \cite{AR}, Theorem~10.4, for manifolds of $C^k$ maps, one checks that the $L^2_k$ completion of $\hbox{Map}^{\Lambda}(\Sigma, M)$ is a $C^{k-\nu _0 -2}$ submanifold of finite codimension in
$$L^2_k(\Sigma, M) \times \left( \Sigma \times \cdots \times \Sigma - \Delta^{(n)} \right),$$
where $\nu _0$ is the maximum of the branching orders $\nu _1, \ldots , \nu_n$ at the points.

Let us determine the tangent space to $\hbox{Map}^{\Lambda}(\Sigma, M)$ in the case where $\Lambda = (1,1, \ldots ,1)$ at a point $(f,(p_1, \ldots p_n))$.  Suppose that, for $t \in (-\epsilon, \epsilon)$,
$$t \mapsto \left(f(t), (\gamma _1(t), \ldots , \gamma _n(t)) \right)$$
is a variation in $\hbox{Map}^{\Lambda}(\Sigma, M)$ with
$$\left(f(0), (\gamma _1(0), \ldots , \gamma _n(0)) \right) = \left( f,(p_1, \ldots , p_n) \right).$$
About a given point $p_i$ we choose local coordinates $(x_1,x_2)$ such that $x_1(p_i) = x_2(p_i) = 0$.  If we differentiate the defining relation $df(t)(\gamma _i(t)) = 0$ with respect to $t$ and set $t = 0$, we obtain
$$\frac{\partial ^2f}{\partial x_a\partial t} (p) + \sum _{b=1}^2\frac{\partial ^2f}{\partial x_a\partial x_b} (p) x_b'(0) = 0.$$
Thus we find that the corresponding tangent vector is
$$(X, v_1, \ldots v_n) \in T_f\hbox{Map}(\Sigma, M) \oplus T_{p_1}\Sigma \oplus \cdots \oplus T_{p_n}\Sigma ,$$
where the components satisfy the conditions
$$\nabla X(p_i) = - d^2f(p_k)(v_i, \cdot ).$$
These conditions could also be obtained from the explicit derivatives of the one-jet evaluation map as described in Abraham and Robbin \cite{AR}, page 27.  Note that $\nabla X(p_i)$ and $d^2f(p_i)(v_i,\cdot)$ are both linear maps from $T_{p_i}\Sigma $ to $T_{f(p_i)}M$.  In the general case where $\Lambda = (\nu _1, \ldots , \nu _n)$, one can show by a similar argument that the tangent space consists of the elements
$$(X, v_1, \ldots v_n) \in T_f\hbox{Map}(\Sigma, M) \oplus T_{p_1}\Sigma \oplus \cdots \oplus T_{p_n}\Sigma $$
which satisfy the conditions
$$ \nabla^j X (p_i) = 0, \hbox{for $1 \leq j \leq \nu _i - 1$} \quad \hbox{and} \quad \nabla ^{\nu _i} X (p_i) = - (d^{\nu _{i+1}}f)(p_i)(v_i, \cdot , \ldots , \cdot ).$$

We have a fibration
$$\hbox{Map}^{\Lambda}(\Sigma, M) \longrightarrow \overbrace{\Sigma \times \cdots \times \Sigma }^n - \Delta^{(n)},$$
in which the fiber over $(p_1, \ldots ,p_n)$ is a submanifold
$$\hbox{Map}^{\Lambda}(\Sigma, M)(p_1, \ldots ,p_n)$$
with tangent space consisting of $X \in T_f\hbox{Map}(\Sigma, M)$ such that 
$$ \nabla^j X (p_i) = 0, \quad \hbox{for $1 \leq j \leq \nu _i$}.$$
Indeed, the manifold $\hbox{Map}^{\Lambda}(\Sigma, M)$ is foliated by these fibers, submanifolds of codimension $2n$.

We now construct a configuration space within ${\mathcal A}^n(\Sigma ,M)$ adapted to the branch type $\Lambda $ by imposing some further conditions on a point
$$((f, (p_1, \ldots ,p_n),\eta ),g)  \in \hbox{Map}^{\Lambda}(\Sigma, M) \times \hbox{Met}_0(\Sigma ) \times \hbox{Met}(M).$$
Suppose that $\Lambda = \nu_1p_1 + \cdots + \nu_n p_n$, focus on one of the points $p_i$, and abbreviate to $(p_i, \nu_i)$ to $(p,\nu)$.  We use the metric $g = \langle \cdot , \cdot \rangle $ in formulating:

\vskip .1in
\noindent
{\bf Lemma 7.1.1.} {\sl If $\sigma \geq \nu$, there exists a finite collection of linear conditions on the $(\sigma +1)$-jet of $f$ at $p$ such that the $(\sigma +1)$-jet of $f$ satisfies these conditions if and only if whenever $z$ is a complex coordinate on $\Sigma $, then
\begin{equation} j_\sigma \left( \frac{\partial f}{\partial z} \right) = j_\sigma \left(z^\nu h(z) \right), \quad \langle h(p), h(p) \rangle = 0, \quad h(p) \neq 0, \label{E:jetconditions}\end{equation}
where $j_\sigma $ denotes $\sigma $-jet and $h$ is a section of $f^*TM \otimes {\mathbb C}$ defined over some open neighborhood $U$ of $p$.}

\vskip .1in
\noindent
(We actually only need the case $\sigma = \nu$.)  It is convenient to express these conditions in terms of canonical isothermal coordinates $(x_1,x_2)$ centered at $p \in \Sigma $ for the metric on $\Sigma $ of constant curvature and total area one, and normal coordinates $(u_1, \ldots ,u_n)$ on $M$ determined by the Riemannian metric on $M$ and centered at $f(p) \in M$.  The jet conditions (\ref{E:jetconditions}) are then expressible in terms of the coefficients of the Taylor expansions.

When $\sigma = \nu$, we can consider the $(\nu +1)$-differential
$$d^{\nu +1}f : T_p \Sigma \times \cdots \times T_p\Sigma \longrightarrow T_{f(p)}M$$
as a polynomial approximation to $f$ near $p$, which we can write as
$$u = (u_1(x_1,x_2), u_2(x_1,x_2), \ldots ,u_n(x_1,x_2)).$$
If $x_1 + i x_2 = r \cos \theta $, we can also write this as
\begin{multline*} u_i(x_1,x_2) = r^{\nu +1} \left[ \sum_{j,k \geq 0, j+k = \nu + 1} a_{ijk} \cos ^j\theta \sin ^k \theta \right] \\ = r^{\nu +1} \left[  \sum \{ c_{ijk}\hbox{exp}\left( \frac{\sqrt{-1} (j-k)\theta }{2} \right) : j + k = \nu + 1, j \geq 0, k \geq 0 \right], \end{multline*}
where the $c_{ijk}$'s are complex constants satisfying the conditions $c_{ijk} = \bar c_{ikj}$.  The jet condition (\ref{E:jetconditions}) can now be expressed in terms of the vectors
$${\bf c}_j = {\bf \bar c}_{\nu - j} = ( c_{1j,\nu-j}, \ldots , c_{nj,\nu-j}), \quad 0 \leq j \leq \nu,$$
as
\begin{equation} {\bf c}_0 \cdot {\bf c}_0 = 0, \quad {\bf c}_0 = {\bf \bar c}_\nu \neq 0, \quad {\bf c}_j = 0 \quad \hbox{for $j \neq 0, \nu$.}\label{E:jetconditions1}\end{equation}
Note that these conditions are invariant under change of normal coordinates on $M$, and they include the condition that the Fourier expansion of the $r^{\nu + 1}$-term in $u_i(r,\theta)$ is a linear combination of $\cos ((\nu + 1) \theta)$ and $\sin ((\nu + 1) \theta)$.

Finally, all the conditions we have described are invariant under the action of the group $\hbox{Diff}_0(\Sigma )$ of diffeomorphisms isotopic to the identity, so we can pass to the quotient obtaining a submanifold
\begin{multline*} {\mathcal A}^\Lambda (\Sigma,M) = \{ ([f, (p_1, \ldots ,p_n),\eta ],g) \in {\mathcal A}^n(\Sigma, M) : \\ (f,(p_1, \ldots p_n)) \in \hbox{Map}^{\Lambda}(\Sigma, M) \hbox{ and}\\
\hbox{ conditions (\ref{E:jetconditions1}) hold at every $p_i$} \}.\end{multline*}
One easily checks that ${\mathcal A}^\Lambda (\Sigma,M)$ is a submanifold of finite codimension in ${\mathcal A}^n(\Sigma, M)$ which contains all conformal harmonic maps $f:\Sigma \rightarrow M$ of branch type $\Lambda $.  Moreover at each point of ${\mathcal A}^\Lambda (\Sigma,M)$ we have a well-defined family
\begin{equation} ([f, (p_1, \ldots ,p_n),\eta ],g) \quad \mapsto \quad (E_1(p_1), \ldots , E(p_n)) \label{E:familytwoplanes}\end{equation}
of two-planes (each $E_i(p_i)$ varying smoothly because the $(\nu _i+ 1)$-jet evaluation map is smooth), such that $E_i(p)$ is spanned by the image of the $(\nu_i+1)$-jet of $f$ at $p_i$, for $1 \leq i \leq n$.

The key feature of the submanifold ${\mathcal A}^\Lambda (\Sigma,M)$ is that it comes equipped with the family of two-planes (\ref{E:familytwoplanes}) which are automatically the images of $\hbox{Re}({\bf L})(p_i)$ whenever the element of ${\mathcal A}^\Lambda (\Sigma,M)$ is a conformal harmonic map.

\subsection{Submanifolds of ${\mathcal A}^n(\Sigma ,M)$ containing ${\mathcal P}^\Lambda $}
\label{S:submanifoldsoftildeM}

Given a branch type $\Lambda = (\nu _1, \ldots , \nu _n)$, we want to cover the space ${\mathcal P}^\Lambda$ of prime minimal surfaces of branch type $\Lambda $ by submanifolds lying in a countable open cover of ${\mathcal A}^n(\Sigma ,M)$, submanifolds which have Fredholm projection to the space of metrics on $M$.  If the genus of $\Sigma $ is at least two, we don't have to deal with a positive-dimensional group $G$ of symmetries, so we treat that simpler case first.  Then we will describe the minor modifications necessary for the case of genus zero or one.

\vskip .1in
\noindent
{\bf Proposition 7.2.1.} {\sl Suppose that $\Sigma $ has genus $\geq 2$.  If $\Lambda = (\nu_1, \ldots , \nu _n)$ is a given branch type of total branching order $\nu $, then there is a countable open covering $\{ {\mathcal W}_i : i \in {\mathbb N} \}$ of ${\mathcal P}^\Lambda$ by open sets in ${\mathcal A}^n(\Sigma ,M)$, and for each $i \in {\mathbb N}$ a submanifold
$${\mathcal Q}^\Lambda_i \subseteq {\mathcal W}_i \subseteq {\mathcal A}^n(\Sigma ,M)$$
such that ${\mathcal P}^\Lambda \cap {\mathcal W}_i \subseteq {\mathcal Q}^\Lambda_i$ and the projection $\pi :  {\mathcal Q}^\Lambda_i \rightarrow \hbox{Met}(M)$ is a proper Fredholm map of Fredholm index zero.}

\vskip .1in
\noindent
To prove this we start with a given conformal harmonic map
$$([f_0,{\bf p}_0, \eta_0], g_0) \in {\mathcal P}^\Lambda \subseteq {\mathcal A}^n(\Sigma ,M) = \tilde{\mathcal M}^n(\Sigma ,M) \times \hbox{Met}(M)$$
of branch type $\Lambda = (\nu _1, \ldots ,\nu _n)$.  Our goal is to construct an open neighborhood ${\mathcal U}$ of this element of ${\mathcal A}^n(\Sigma ,M)$ and a submanifold
$${\mathcal Q}^\Lambda \subseteq {\mathcal U} \subseteq {\mathcal A}^n(\Sigma ,M)$$
so that all conformal harmonic maps of branch type $\Lambda $ within ${\mathcal U}$ lie in ${\mathcal Q}^\Lambda$ and the projection from ${\mathcal Q}^\Lambda $ to $\hbox{Met}(M)$ is Fredholm of Fredholm index zero.  

Use of a local section allows us to simplify matters by replacing the space ${\mathcal A}^n(\Sigma ,M)$ by an open neighborhood in the product
$$\hbox{Map}(\Sigma, M) \times \left( \Sigma ^n - \Delta^{(n)} \right) \times {\mathcal T}_g \times \hbox{Met}(M).$$
Thus we seek to construct a neighborhood ${\mathcal U}$ of a given conformal harmonic map,
$$(f_0,{\bf p}_0, \omega_0, g_0) \in \hbox{Map}(\Sigma, M) \times \left( \Sigma ^n - \Delta^{(n)} \right) \times {\mathcal T}_g \times \hbox{Met}(M),$$
of branch type $\Lambda$, and a submanifold ${\mathcal Q}^\Lambda \subseteq {\mathcal U}$ so that all conformal harmonic maps of branch type $\Lambda $ within ${\mathcal U}$ lie in ${\mathcal Q}^\Lambda$.  We can assume that ${\mathcal U}$ is the domain for a submanifold chart for the submanifold ${\mathcal A}^\Lambda (\Sigma,M)$ constructed in the previous section.  Thus we can assume that the model Banach space for ${\mathcal A}^n(\Sigma ,M)$ is a direct sum $E \oplus F$, where $F$ is finite-dimensional and that we are given a submanifold chart
$$\phi : {\mathcal U} \longrightarrow E \oplus F, \quad \hbox{with} \quad \phi({\mathcal A}^\Lambda (\Sigma,M) \cap {\mathcal U}) = \phi ({\mathcal U}) \cap ( E \oplus \{ 0 \}),$$
the dimension of $F$ being just the codimension of ${\mathcal A}^\Lambda (\Sigma,M)$ in ${\mathcal A}^n(\Sigma ,M)$.  Moreover, we can assume that
$$\phi \left( f_0,{\bf p}_0, \omega_0, g_0 \right) = 0.$$

Within the open set ${\mathcal U}$, we want to obtain a submanifold of lower codimension than ${\mathcal A}^\Lambda (\Sigma,M) \cap {\mathcal U}$, in fact of finite codimension.  To do this, we first use the submanifold chart described above to extend the two-planes (\ref{E:familytwoplanes}) from the submanifold ${\mathcal A}^\Lambda (\Sigma,M) \cap {\mathcal U}$ to all of ${\mathcal U}$.  We then utilize complex coordinates $z_i$ on $\Sigma $ in a neighborhood of each $p_i$, and the Levi-Civita connection $D$ for the Riemannian metric $g$ on $M$.  If $\Pi_{E(p_i)}$ denotes the orthogonal projection onto the two-plane $E(p_i)$, we can then impose the conditions that
\begin{multline} \Pi_{E(p_i)}\frac{\partial f}{\partial z_i}(p_i) = \Pi_{E(p_i)}\frac{D^2 f}{\partial z_i^2}(p_i) = \cdots = \Pi_{E(p_i)}\frac{D^{\nu _i} f}{\partial z_i^{\nu _i}}(p_i) = 0, \\ \Pi_{E(p_i)}\frac{D^{\nu _i+1} f}{\partial z_i^{\nu _i+1}}(p_i) \neq 0, \label{E:vanishingofnuijet}\end{multline}
where $\nu _i$ is the branching order at $p_i$.  (The inequality in (\ref{E:vanishingofnuijet}) is included to eliminate branch types which are $> \Lambda $.)  Thus we impose only $2\nu_i$ real linear conditions at the branch point $p_i$.  Imposing the conditions at all of the branch points yields a new larger configuration space ${\mathcal B}^\Lambda (\Sigma ,M)$, a submanifold of codimension $2 \nu$ in ${\mathcal A}^n(\Sigma ,M) \cap {\mathcal U}$, where $\nu $ is the total branching order.  All conformal harmonic maps of branch type $\Lambda $ within ${\mathcal U}$ lie within ${\mathcal B}^\Lambda(\Sigma ,M)$ because the differentials of $f$ at each branch point $p_i$ must vanish up to order $\nu _i$.  Moreover, the codimension of ${\mathcal B}^\Lambda(\Sigma ,M)$ exactly equals the dimension of the space of tangential Jacobi fields that are introduced by the branch points.  By Proposition~5.3.3, the conditions (\ref{E:vanishingofnuijet}) prevent any elements from the space of tangential Jacobi fields $(X, \dot \omega ) = \hbox{Re}(Z, \bar \psi d\bar z^2)$ described in \S\ref{S:tangentialjacobi} from lying in the tangent space to ${\mathcal B}^\Lambda(\Sigma ,M)$.

As described before, the familiar Euler-Lagrange map
\begin{equation} F : \hbox{Map}(\Sigma ,M) \times {\mathcal T}_g \times \hbox{Met}(M) \longrightarrow \hbox{Map}(\Sigma ,TM) \oplus T{\mathcal T}_g \label{E:ELoriginal}\end{equation}
defines a section of a smooth bundle ${\mathcal E}$ over $\hbox{Map}(\Sigma ,M) \times {\mathcal T}_g \times \hbox{Met}(M)$ in such a way that conformal harmonic maps are just the elements of $F^{-1}({\mathcal Z})$, where ${\mathcal Z}$ is the image of the zero section in ${\mathcal E}$.   This induces an Euler-Lagrange map on the mapping space ${\mathcal A}^n(\Sigma ,M)$ defined by (\ref{E:tildem}),
\begin{multline} \hat F : {\mathcal A}^n(\Sigma ,M) \\ \longrightarrow T\tilde{\mathcal M}^n(\Sigma ,M) = T\hbox{Map}(\Sigma ,M) \oplus (T\Sigma \oplus \cdots \oplus T\Sigma ) \oplus T{\mathcal T}_g, \label{E:ELhat}\end{multline}
in which composition with the projection
$$ \pi : T\tilde{\mathcal M}^n(\Sigma ,M) \longrightarrow T\Sigma \oplus \cdots \oplus T\Sigma $$
is just the zero map.  We restrict this $\hat F$ to an Euler-Lagrange map on the submanifold ${\mathcal B}^\Lambda (\Sigma ,M) \subseteq {\mathcal U}$ that we have just constructed:
\begin{equation} F^\Lambda  : {\mathcal B}^\Lambda (\Sigma ,M) \longrightarrow T\hbox{Map}(\Sigma ,M) \oplus (T\Sigma \oplus \cdots \oplus T\Sigma ) \oplus T{\mathcal T}_g. \label{E:ELqlambdai}\end{equation}

Differentiation of the original Euler-Lagrange map (\ref{E:ELoriginal}) yields the Jacobi operator $L_{(f, \omega ,g)} = \pi _V \circ (D_1F)_{(f, \omega ,g)}$, which is Fredholm of index zero (since we are in the special case where $G$ is trivial).  Since $(f_0, \omega _0, g_0)$ has branch type $\Lambda $ of total branching order $\nu$,
\begin{equation} L_{(f_0, \omega _0,g_0)} : T_{f_0}\hbox{Map}(\Sigma ,M) \oplus T_{\omega _0}{\mathcal T}_{g} \longrightarrow T_{f_0}\hbox{Map}(\Sigma ,M) \oplus T_{\omega _0}{\mathcal T}_g \label{E:originaljacobioperator}\end{equation}
has a $(2\nu)$-dimensional space of tangential Jacobi fields,
\begin{equation} \hbox{Re} \left( {\cal J}({\bf L}) \right) = \hbox{Span} \{ (X_1 \dot \omega _1), \ldots, (X_{2\nu}, \dot \omega _{2\nu})\}, \label{E:tangentialjacobifieldsbranchtypelambda}\end{equation}
and Jacobi fields in this space are just as smooth as the conformal harmonic map $(f_0,\omega _0, g_0)$.  The $(2\nu )$-dimensional space spanned by these Jacobi fields can be translated about the coordinate neighborhood ${\mathcal U}$ thereby generating a subbundle ${\mathcal F} \subseteq {\mathcal E}$ of fiber dimension $2 \nu$.  We can then then define an $L^2$-orthogonal projection
$$p : {\mathcal E} \longrightarrow {\mathcal F}^\bot$$
to the $L^2$-orthogonal complement ${\mathcal F}^\bot$ to ${\mathcal F}$ in ${\mathcal E}$.

Let ${\mathcal Z}$ be the image of the zero-section in the restriction of ${\mathcal F}^\bot$ to ${\mathcal B}^\Lambda (\Sigma ,M)$.  We claim that $p \circ F^\Lambda $ is transverse to ${\mathcal Z}$ at the point $(f_0, \omega _0, g_0)$, and hence transverse to ${\mathcal Z}$ in some neighborhood ${\mathcal B}^\Lambda (\Sigma ,M) \cap {\mathcal W}$ of $(f_0, \omega _0, g_0)$, with ${\mathcal W} \subseteq {\mathcal U}$.  Then ${\mathcal Q} = (p \circ F^\Lambda )^{-1}({\mathcal Z}) \cap {\mathcal W}$ is a submanifold of ${\mathcal W}$ which contains all of the conformal harmonic maps of branch type $\Lambda $ within ${\mathcal W}$, and we will then show that ${\mathcal Q}$ has Fredholm projection to the space of metrics with Fredholm index zero.  Recall that Fredholm maps are locally proper by Theorem~1.6 of \cite{Sm2}.  Thus, since the Sobolev completions have countable bases, it will then follow that we can cover ${\mathcal P}^\Lambda$ by a countable collection of such open sets ${\mathcal W}_i$ satisfying all the conditions of Proposition~7.2.1, and thereby finishing its proof.

To prove transversality at $(f_0, \omega _0, g_0)$, we need to note that the original Jacobi operator (\ref{E:originaljacobioperator}) induces a corresponding operator
$$\hat L : T_{(f,(p_1, \ldots , p_n), \omega)}\tilde{\mathcal M}^n(\Sigma ,M) \longrightarrow T_{(f,(p_1, \ldots , p_n), \omega)}\tilde{\mathcal M}^n(\Sigma ,M)$$
which is also Fredholm of index zero and divides into a direct sum $\hat L = L \oplus I$, where
$$I : T_{p_1}\Sigma \oplus \cdots \oplus T_{p_n}\Sigma \longrightarrow T_{p_1}\Sigma \oplus \cdots \oplus T_{p_n}\Sigma$$
is just the identity map.  The identity summand in $\hat L$ shows that inclusion of the factor specifying the branch points has no real effect on the Euler-Lagrange operator, and indeed that factor was just included so we could properly define the subspace ${\mathcal B}^\Lambda (\Sigma ,M)$.  Finally, we restrict the operator $\hat L$ to our new configuration space to obtain
$$L^\Lambda : T_{(f,{\bf p}, \omega )}{\mathcal B}^\Lambda (\Sigma ,M) \longrightarrow T_{(f,{\bf p}, \omega)}\tilde{\mathcal M}^n(\Sigma ,M),$$
a Fredholm operator of Fredholm index $- 2 \nu$, because we have cut down the dimension of the domain with the $2 \nu$ real linear conditions described in (\ref{E:vanishingofnuijet}).  We then compose with the projection into ${\mathcal F}^\bot$ obtaining a map
\begin{equation} p \circ L^\Lambda : T_{(f,{\bf p}, \omega )}{\mathcal B}^\Lambda (\Sigma ,M) \longrightarrow {\mathcal F}^\bot_{(f,{\bf p}, \omega)}, \label{E:pLlambda}\end{equation}
which has Fredholm index zero at $(f_0, \omega _0, g_0)$ since we have eliminated spaces of equal dimension from the kernel and cokernel.

Since $(f_0, \omega _0, g_0)$ is prime, it is somewhere injective, hence injective on an open dense subset of $\Sigma $, so we can apply Lemma~3.3 to show that the derivative with respect to the metric
$$\pi _V \circ D_2F(f_0, \omega _0, g_0) : T_{g_0}\hbox{Met}(M) \longrightarrow T_{f_0}\hbox{Map}(\Sigma ,M) \oplus T_\omega{\mathcal T}_{g_0}$$
covers the complement of the space spanned by the tangential Jacobi fields (\ref{E:tangentialjacobifieldsbranchtypelambda}) in the space of all Jacobi fields.  On the other hand the tangential Jacobi fields (\ref{E:tangentialjacobifieldsbranchtypelambda}) get taken to zero under the projection $p$, so
$$p \circ \pi _V \circ D_2F^\Lambda(f_0, \omega _0, g_0) : T_{g_0}\hbox{Met}(M) \longrightarrow {\mathcal F}^\bot_{(f_0, \omega _0, g_0)}$$
covers a complement to the range of the projected Jacobi operator of (\ref{E:pLlambda}).  This implies transversality at the point $(f_0, \omega _0, g_0)$.  Transversality is an open condition, so $F^\Lambda |({\mathcal B}^\Lambda (\Sigma ,M) \cap {\mathcal V})$ must be transversal to ${\mathcal Z}$ when ${\mathcal W}$ is a sufficiently small neighborhood of $(f_0, \omega _0, g_0)$.  Thus
$${\mathcal Q} = (p \circ F^\Lambda )^{-1}({\mathcal Z}) \cap \left({\mathcal B}^\Lambda (\Sigma ,M) \cap {\mathcal W}\right)$$
is indeed a submanifold containing all points of ${\mathcal P}^\Lambda \cap {\mathcal W}$.

The proof that the projection is Fredholm of Fredholm index zero is a straightforward modification of the argument for Lemma~3.5.  Indeed,\begin{multline} T_{(f,\omega ,g)}{\mathcal Q} = \{ (X, \dot \omega, \dot g) \in T_{(f,\omega ,g)}{\mathcal B}^\Lambda (\Sigma ,M) : \\ \hbox{$p \circ L^\Lambda (X,\dot \omega ) + p \circ \pi _V \circ (D_2F)_{(f,\omega ,g)} (\dot g) = 0$} \}, \label{E:tangenttoQ} \end{multline}
where the Jacobi operator $p \circ L^\Lambda$ is of Fredholm index zero, with
$$\hbox{Kernel of $(p \circ L^\Lambda)$} = {\mathcal J} \cong (\hbox{complement to Range of $(p \circ L^\Lambda)$}),$$
where ${\mathcal J}$ is a finite-dimensional space of nontangential Jacobi fields.  For the differential of the projection,
$$d\pi_{(f,\omega ,g)}: T_{(f,\omega ,g)} {\mathcal Q} \rightarrow T_g\hbox{Met}(M),$$
we find that
$$\hbox{Ker}(d\pi_{(f,\omega ,g)}) \cong \hbox{Ker}(p \circ L^\Lambda),$$
while the range of $d\pi_{(f ,\omega ,g)}$ consists of the elements
\begin{multline*} \dot g \in T_g\hbox{Met}(M) \quad \hbox{such that} \quad T(\dot g) \in \hbox{Range}(p \circ L^\Lambda), \\ \hbox{where} \quad T(\dot g) = \pi _V \circ (D_2F)(\dot g).\end{multline*}
Thus the range of $d\pi_{(f ,\omega ,g)}$ is the preimage of a closed space of finite codimension, hence a closed subspace of finite codimension itself.  It follows that $d\pi_{(f,\omega ,g)}$ is indeed a Fredholm map, and the dimension of the cokernel of $d\pi_{(f,\omega ,g)}$ is no larger than the dimension of the cokernel of $L$.  But tranversality shows that $d(p \circ F^\Lambda )_{(f ,\omega ,g)}$ is surjective, so anything not in the image of $(p \circ L^\Lambda)$ must be covered by $T$ and the cokernels of $(p \circ L^\Lambda)$ and $d\pi $ must have equal dimension.  Hence $d\pi_{(f,\omega ,g)}$ is a Fredholm map of Fredholm index zero, and the required analog of Lemma~3.5 is proven.

\vskip .1in
\noindent
It remains to discuss the cases in which $\Sigma $ has genus zero or one, and so has a positive-dimensional group $G$ of symmetries.  In these cases, we need to modify the previous statement a little:

\vskip .1in
\noindent
{\bf Proposition 7.2.2.} {\sl Suppose that $\Sigma $ has genus zero or one, and $\Lambda = (\nu_1, \ldots , \nu _n)$ is a given branch type.  Then there is a countable collection $\{ {\mathcal W}_i : i \in {\mathbb N} \}$ of open sets in ${\mathcal A}^n(\Sigma ,M)$, and for each $i \in {\mathbb N}$ a submanifold
$${\mathcal Q}^\Lambda_i \subseteq {\mathcal W}_i \subseteq {\mathcal A}^n(\Sigma ,M) = \tilde{\mathcal M}^n(\Sigma ,M) \times \hbox{Met}(M)$$
such that every element of ${\mathcal P}^\Lambda $ lies on the $G$-orbit of some element in ${\mathcal Q}^\Lambda_i$, for some $i$, and the projection $\pi :  {\mathcal Q}^\Lambda_i \rightarrow \hbox{Met}(M)$ is a proper Fredholm map of Fredholm index zero.}

\vskip .1in
\noindent
To prove this, we need to modify the argument of Proposition~7.2.1 to allow for the $G$-action.  To destroy the $G$-symmetry in the two cases $\Sigma = S^2$ or $T^2$, we replace $\hbox{Map}(S^2,M)$ by the open subset ${\mathcal F}_j(S^2,M)$ of (\ref{E:defoffiS2}) or $\hbox{Map}(T^2,M)$ by the open subset ${\mathcal F}_j(T^2,M)$ of (\ref{E:defoffiT2}).  This has the effect of replacing every $G$-orbit by a finite number of points on that orbit.  We can then construct local configuration spaces
$${\mathcal A}^n_j(\Sigma ,M) = {\mathcal F}_j(\Sigma,M) \times \left( \Sigma ^n - \Delta^{(n)} \right) \times {\mathcal T}_g \times \hbox{Met}(M) $$
such that every conformal harmonic maps $f : \Sigma \rightarrow M$ of branch type $\Lambda $ has a $G$ orbit which intersects one of the spaces ${\mathcal A}^n_j(\Sigma ,M)$.  The constructions of \S \ref{S:stratification} can be applied to each ${\mathcal A}^n_j(\Sigma ,M)$ yielding a submanifold ${\mathcal A}^\Lambda_j(\Sigma ,M)$ of finite codimension in ${\mathcal A}^n_j(\Sigma ,M)$ as before.

Now we run through the proof of Proposition~7.2.1 starting with a conformal harmonic map
$$([f_0,{\bf p}_0, \eta_0], g_0) \in {\mathcal P}^\Lambda\cap{\mathcal A}^\Lambda_j(\Sigma ,M) \subseteq {\mathcal A}^n_j(\Sigma ,M).$$
Exactly the same proof leads to a countable covering $\{ {\mathcal W}_i : i \in {\mathbb N} \}$ of ${\mathcal P}^\Lambda \cap {\mathcal A}^n_j(\Sigma ,M)$ by open sets in ${\mathcal A}^n_j(\Sigma ,M)$, together with a countable collection of submanifolds,
$${\mathcal Q}^\Lambda_i \subseteq {\mathcal W}_i \subseteq {\mathcal A}^n_j(\Sigma ,M),$$
such that
$${\mathcal P}^\Lambda \cap {\mathcal W}_i \cap {\mathcal A}^n_j(\Sigma ,M) \subseteq {\mathcal Q}^\Lambda_i$$
and each projection ${\mathcal Q}^\Lambda_i \rightarrow \hbox{Met}(M)$ is a proper Fredholm map of Fredholm index zero.  This is exactly what is needed to establish Proposition~7.2.2.

\vskip .1in
\noindent
{\bf Corollary 7.2.3.} {\sl Suppose that $\Sigma $ is an oriented compact connected Riemann surface and $\Lambda $ is a given branch type.  For generic choice of metric on $M$, there are only countably many prime minimal surfaces $(f,\omega)$ from $\Sigma $ into $M$ of branch type $\Lambda $.}

\vskip .1in
\noindent
Of course, nongeneric metrics, for example the standard metrics on Hermitian symmetric spaces, can have uncountably many minimal immersions from spheres, for example.

\subsection{Generic self-intersection sets}
\label{S:genericselfintersection}

In the $m$-fold cartesian product $\Sigma ^m$, for $m$ a positive integer, we consider the subset
$$\Sigma ^m - \Delta ^{(m)} = \{ (q_1, \ldots q_m) \in \Sigma ^m : \hbox{$q_i \neq q_j$ when $i\neq j$} \}.$$
We also consider the multidiagonal in the $m$-fold cartesian product $M^m$,
$$\Delta _m = \{ (r_1, \ldots r_m) \in M^m : r_1 = r_2 = \cdots = r_m \}.$$
We say that an immersion $f : \Sigma \rightarrow M$ has {\em transversal crossings\/} if for every $m > 1$, the restricted product map,
\begin{equation} f^m = \overbrace{f \times \cdots \times f}^m : \Sigma ^m - \Delta ^{(m)} \longrightarrow M^m \quad \hbox{is transversal to $\Delta _m$.}\label{E:mtransversal}\end{equation} 
Thus if $\Sigma $ is a compact surface and $M$ is a smooth manifold of dimension at least five, an immersion with transversal crossings is a one-to-one immersion and hence an imbedding, while if $M$ has dimension four, such an immersion has only double points and the intersections at double points are transverse.

\vskip .1 in
\noindent
{\bf Proposition~7.3.1.} {\sl If the dimension of $M$ is at least four, then for generic choice of metric on $M$, the self-intersection set of the restriction of any conformal harmonic map $f : \Sigma \rightarrow M$ to $\Sigma - \{\hbox{branch points}\}$ is an immersion with transversal crossings.  Indeed, the space of metrics which admit conformal harmonic maps whose restriction to $\Sigma - \{\hbox{branch points}\}$ does not have transversal crossings has codimension at least two in the space of metrics, when $\dim M \neq 5$, and at least one when $\dim M = 5$.}

\vskip .1 in
\noindent
Proof:  This is a relatively straightforward consequence of the Transversal Density Theorem 19.1 of Abraham and Robbin \cite{AR}, but to apply this result, we need to construct a countable family of representations from open subsets of the space of metrics to families of maps from $\Sigma $ to $M$ such that the $G$-orbit of each prime minimal surface of a given branch type contains an element in one of the families.

However, that is exactly what is provided by Propositions~7.2.1 and 7.2.2.  Once we have fixed a branch type $\Lambda $ which has $n$ branch points, these propositions provide a countable collection $\{ {\mathcal W}_i : i \in {\mathbb N} \}$ of open sets in the configuration space ${\mathcal A}^n(\Sigma ,M)$, and for each ${\mathcal W}_i$ a submanifold
$${\mathcal Q}^\Lambda_i \subseteq {\mathcal W}_i \subseteq {\mathcal A}^n(\Sigma ,M) = \tilde{\mathcal M}^n(\Sigma ,M) \times \hbox{Met}(M)$$
such that every element of the space ${\mathcal P}^\Lambda $ of prime minimal surfaces of branch type $\Lambda $ lies on the $G$-orbit of some element in ${\mathcal Q}^\Lambda_i$, for some $i$.  Moreover, the projection $\pi :  {\mathcal Q}^\Lambda_i \rightarrow \hbox{Met}(M)$ is a proper Fredholm map of Fredholm index zero.  If
$$([f_i,{\bf p}_i, \eta_i], g_i) \in {\mathcal P}^\Lambda \subseteq {\mathcal Q}^\Lambda_i \subseteq {\mathcal A}^n(\Sigma ,M),$$
we can use the inverse function theorem to construct a local section
$$\sigma _i : {\mathcal V}_i \longrightarrow {\mathcal Q}^\Lambda_i \subseteq {\mathcal W}_i,$$
defined on some open neighborhood ${\mathcal V}_i$ of $g_i$, which contains the given element of ${\mathcal P}^\Lambda$.  Using the countable basis for the topology, we thus obtain a sequence of sections
$$\{ \sigma _i : {\mathcal V}_i \rightarrow {\mathcal W}_i \subseteq {\mathcal A}^n(\Sigma ,M) :  i \in {\mathbb N} \}$$
which contain all the elements of ${\mathcal P}^\Lambda $ in their images, and a corresponding sequence
$$\{ \rho_i = \pi \circ \sigma _i : {\mathcal V}_i \rightarrow \hbox{Map}(\Sigma ,M) : i \in {\mathbb N} \}$$
where $\pi $ denotes the projection into the manifold of maps.  After choosing sufficiently strong Sobolev norms, we can arrange that each $\rho _i$ defines a $C^r$ representation for a large positive integer $r$, in the terminology of Chapter~4 of Abraham and Robbin \cite{AR}.  This provides most of the hypotheses for application of the Transversal Density Theorem.

The one minor difficulty is that there is no hope of achieving transversality at the branch points.  So starting with a given element
$$([f_0,{\bf p}_0, \eta_0], g_0) \in {\mathcal P}^\Lambda \subseteq {\mathcal Q}^\Lambda \subseteq {\mathcal A}^n(\Sigma ,M),$$
where ${\mathcal Q}^\Lambda$ is one of the countable sequence of submanifolds described in the previous paragraph, with local section
$$\sigma : {\mathcal V} \longrightarrow {\mathcal Q}^\Lambda, \quad {\mathcal V} \subseteq \hbox{Met}(M),$$
we excise the branch points by the following construction:  Given $\epsilon _j = 1/2^j$, we let
$$D_{\epsilon _j}({\bf p}_0) = \{ q \in \Sigma : d(q,{\bf p}_0) < \epsilon _j \}, \quad K_j = \Sigma - D_{\epsilon _j}({\bf p}_0).$$
(We can assume that $D_{\epsilon _j}({\bf p}_0)$ consists of $m$ disjoint disks centered at the points of the branch locus ${\bf p}_0$.)  It suffices to show that for arbitrary $j \in {\mathbb N}$, and for a suitable contraction ${\mathcal V}_j$ of ${\mathcal V}$ to a possibly smaller neighborhood of $([f_0,{\bf p}_0, \eta_0], g_0)$ (chosen so that all branch points of elements of $\sigma ({\mathcal V}_j)$ lie within $D_{\epsilon _j}({\bf p}_0)$) a generic metric $g$ within ${\mathcal V}_j$ has the property that
$$\rho (g) = \pi \circ \sigma (g) \in  \hbox{Map}(\Sigma ,M)$$
is an immersion with transversal crossings on $K_j$, that is that
$$f^m = f \times \cdots \times f : \Sigma ^m - \Delta ^{(m)} \longrightarrow M^m$$
is transversal to $\Delta _m$ on $K_j^m - K_j^m\cap \Delta^{(m)}$.

Now starting with $\rho : {\mathcal V} \rightarrow \hbox{Map}(\Sigma,M)$ we can build the map
$$\rho ^m :  {\mathcal V} \rightarrow \hbox{Map}(\Sigma^m,M^m),$$
and from it the evaluation map
\begin{equation} \hbox{ev}_{\rho ^m} : {\mathcal V} \times \Sigma^m \longrightarrow M^m. \label{E:evmtransversality} \end{equation}
The Transversal Density Theorem of \cite{AR} states that if this map is transversal to the multidiagonal $\Delta_m$ on
$${\mathcal V}_j \times (K_j^m - K_j^m\cap \Delta^{(m)}),$$
then the metrics $g \in {\mathcal V}_j$ such that the map $\rho (g)$ satisfies (\ref{E:mtransversal}) form a residual set.  From this it will follow that the metrics $g \in {\mathcal V}_j$ such that $\rho (g)$ is an immersion with transversal crossings on $K_j$ are residual.

Thus we need to check transversality of the evaluation map (\ref{E:evmtransversality}), and this just requires constructing perturbations in the metric that will eliminate nontransversal intersections.  It is here that the arguments of \S~\ref{S:selfintersection} play a crucial role, because they eliminate a priori many possibilities of self-intersections, and indeed Lemma~6.3 shows that when $m = 2$, we need only check curves of self-intersection, and that isolated double points automatically have tangent planes which intersect at a single point.

Suppose that $(f,\omega )$ is a conformal harmonic map into the Riemannian manifold $(M,g)$ and let $U_1$ and $U_2$ be two disjoint open balls in $\Sigma $ such that $f(U_1)$ and $f(U_2)$ intersect along a $C^1$ curve $C$ within $M$, the intersection being transverse.  Suppose that $V$ is an open set in $M$ such that
$$f^{-1}(V) = U_1 \cup U_2 \cup (\hbox{open components disjoint from $U_1$ and $U_2$}).$$
In analogy with the proof of Lemma~3.2, we can suppose that $V$ is the domain of local coordinates $(u_1 ,\ldots , u_n)$ such that 
\begin{enumerate}
\item $C$ is described by the equations $u_2 = \cdots = u_n = 0$,
\item $f(U_1)$ is described by the equations $u_3 = \cdots = u_n = 0$,
\item $u_a \circ f = x_a$ on $U_1$, for $a = 1,2$, where $x_1 + i x_2$ is a conformal parameter on $U$,
\item $f(U_2)$ is described by the equations $u_2 = 0$ and $u_4 = \cdots = u_n = 0$,
\item $u_1 \circ f = y_1$ and $u_3 \circ f = y_2$ on $U_2$, where $(y_1,y_2)$ is a coordinate system on $U_2$, and
\item the Riemannian metric $g$ on $V$ takes the form $\sum g_{ij} du_idu_j$, such that when restricted to $f(\Sigma ) \cap V$, $g_{ir} = \delta _{ir}$, for $1 \leq i \leq n$ and $4 \leq r \leq n$.
\end{enumerate}
Just as in (\ref{E:conformalfactor}), we can write
$$g_{ab} = \sigma ^2 \eta_{ab}, \quad \eta_{ab} = \lambda ^2 \delta_{ab}, \quad \hbox{for $1 \leq a,b \leq 2$,}$$
although we make no restriction on $g_{i3}$ for $1 \leq i \leq 3$.

We can define a metric variation of the form 
$$\dot g_{22}(u_1, \ldots ,u_n) = \sum_{r=4}^n u_r \phi_r (u_1,\ldots, u_n), \quad \dot g_{ij} = 0 \quad \hbox{when $(ij) \neq (22)$},$$
where $\phi _r$ has compact support within $V$.  It then follows that
$$\dot \Gamma _{r,22} = \frac{1}{2} \frac{\partial g_{22}}{\partial u_r} = \frac{1}{2}\phi_r(u_1,u_2, 0, \ldots ,0), \quad \hbox{for $4 \leq r \leq n$,}$$
while if $1 \leq i,j \leq 3$ and $(i,j) \neq (2,2)$, then
$$\dot \Gamma _{r,ij} = 0, \quad \hbox{for $4 \leq r \leq n$.}$$
Thus if we differentiate the tension
$$\tau (f) = \frac{1}{\lambda ^2} \left(\frac{\partial ^2 u_k}{\partial x_1^2} + \frac{\partial ^2 u_k}{\partial x_2^2}\right) + \sum _{i,j} \frac{1}{\lambda ^2} \Gamma ^k_{ij} \left(\frac{\partial u_i}{\partial x_1} \frac{\partial u_j}{\partial x_1} + \frac{\partial u_i}{\partial x_2} \frac{\partial u_j}{\partial x_2}\right)$$
with respect to the metric variation, we obtain
$$\dot \tau (f)|U_1 = \frac{1}{2\lambda ^2}\sum_{r=4}^n\phi_r(u_1,u_2, 0, \ldots ,0)\frac{\partial }{\partial u_r} \quad \hbox{while} \quad \dot \tau (f)|U_2 = 0,$$
the latter equation holding because $u_2 \equiv 0$ on $f(U_2)$.

Suppose that $q_1$ and $q_2$ are elements of $U_1$ and $U_2$ respectively with $f(q_1) = f(q_2) \in C$.  The preceding argument shows we can perturb the metric fixing $f|U_2$ while moving $f(U_1)$ in any direction normal to the three-dimensional space spanned by
$$f_*(T_{q_1}\Sigma ) + f_*(T_{q_2}\Sigma ) \subseteq T_{f(q_1)}M.$$
We can clearly do the same construction, fixing instead $f|U_1$ while moving $f(U_2)$.  From this it follows that the evaluation map $\hbox{ev}_{\rho ^2}$ of (\ref{E:evmtransversality}) is transversal to the diagonal $\Delta \subseteq M^2$ at the point $(g,q_1,q_2)$.  Since we have two independent variations, it follows from Lemma~19.3 in \cite{AR} that metrics admitting nontransversal immersions have codimension two.

In more detail (following the notation of \cite{AR}), if we set
\begin{equation} X = (\Sigma -K_j) \times (\Sigma -K_j) - \Delta \quad \hbox{and} \quad W = \Delta _{M\times M}, \label{E:notation1}\end{equation}
and we let
\begin{equation} {\mathcal B} = \left(\hbox{ev}_{\rho ^2}\right)^{-1}(W) \subseteq {\mathcal V} \times X, \label{E:notation2}\end{equation}
then $\pi : {\mathcal B} \rightarrow {\mathcal V}$ is a Fredholm map of index $\leq -2$, and that the regular values for $\pi $ are exactly the maps for which $\rho _2$ is transverse to the diagonal.

To finish the proof of Proposition~7.3.1 we need to check transversality of (\ref{E:evmtransversality}) for $m \geq 3$, but the needed perturbations in these cases are much easier to construct, so we leave that verification to the reader.  Note that if $\dim M = 4$, double points are automatically transversal by Lemma~6.3 (so long as there are no curves of intersection).  When $\dim M = 5$, double points have codimension one, and it for this reason that we must single out the case of dimension five in the statement of the proposition.

\vskip .1 in
\noindent
We can say that $f : \Sigma \rightarrow M$ has {\em weakly transversal crossings\/} if the only singularities of $f$ are double points at each of which the two tangent planes intersect in a single point.  We can then remove the restriction $\dim M \neq 5$ in the statement of Proposition~7.3.1 by replacing \lq\lq transversal" by \lq\lq weakly trransversal."  Indeed, we could introduce notation similar to (\ref{E:notation1}) and (\ref{E:notation2}) for all the $\rho_m$'s, and restate what we have proven as follows:

\vskip .1 in
\noindent
{\bf Proposition~7.3.2.} {\sl There is an open cover $\{ {\mathcal V}_i : i \in {\mathbb N} \}$ of $\hbox{Met}(\Sigma )$ and a corresponding collection of Fredholm maps $\pi _i : {\mathcal B}_i \rightarrow {\mathcal V}_i$ of Fredholm index $\leq -2$, such that 
\begin{equation} g \in {\mathcal V}_i \subseteq \hbox{Met}(\Sigma ) \quad \Rightarrow \quad \hbox{$g$ is a regular value for $\pi _i$}, \label{E:good} \end{equation}
then every parametrized minimal surface $f : \Sigma \rightarrow M$ (for fixed oriented surface $\Sigma $) has weakly transversal crossings (in addition to possible branch points).}

\subsection{Proof of the Main Theorem}
\label{S:proof}

We now make a small adjustment to the definition (\ref{E:plambda}) of prime conformal harmonic maps of branch type $\Lambda$:
\begin{multline}{\mathcal P}^\Lambda_0 = \{ ([f,{\bf p}, \eta], g) \in {\mathcal P}^\Lambda \hbox{ such that the only singularities of } \\ \hbox{$f$ are the branch points and weakly transversal crossings} \}. \label{E:plambda0}\end{multline}
The Main~Theorem is a consequence of Propositions~7.3.1 and 7.3.2, the Sard-Smale Theorem and the following proposition:

\vskip .1 in
\noindent
{\bf Proposition~7.4.1.} {\sl Suppose that $\Lambda = (\nu _1, \ldots ,\nu _n)$ is a given branch type, $\Lambda \neq \emptyset$.  Then there is a countable open cover $\{ {\mathcal W}_i : i \in {\mathbb N} \}$ of ${\mathcal P}^\Lambda_0 \subseteq {\mathcal A}^n(\Sigma ,M)$, and for each $i \in {\mathbb N}$ a submanifold
$${\mathcal R}^\Lambda_i \subseteq {\mathcal W}_i \subseteq {\mathcal A}^n(\Sigma ,M)$$
such that every element of ${\mathcal P}^\Lambda _0$ lies on the $G$-orbit of some element in ${\mathcal R}^\Lambda_i$, for some $i$, and the projection $\pi :  {\mathcal R}^\Lambda_i \rightarrow \hbox{Met}(M)$ is a proper Fredholm map of Fredholm index $-2n$, where $n$ is the number of branch points.}

\vskip .1 in
\noindent
The proof of Proposition~7.4.1 is a very slight modification of the proofs for Propositions~7.2.1 and 7.2.2.

Assume first $\Sigma $ has genus at least two.  As before, we start with a given conformal harmonic map
$$([f_0,{\bf p}_0, \eta_0], g_0) \in {\mathcal P}^\Lambda_0 \subseteq {\mathcal A}^n(\Sigma ,M) = \tilde{\mathcal M}^n(\Sigma ,M) \times \hbox{Met}(M)$$
of branch type $\Lambda$, and construct an open neighborhood ${\mathcal U}$ of this element in ${\mathcal A}(\Sigma ,M)$, and construct a submanifold
$${\mathcal R}^\Lambda \subseteq {\mathcal Q}^\Lambda \subseteq{\mathcal U} \subseteq {\mathcal A}^n(\Sigma ,M),$$
where ${\mathcal Q}^\Lambda$ is the submanifold constructed in Proposition~7.2.1, so that all conformal harmonic maps of branch type $\Lambda $ within ${\mathcal U}$ lie in ${\mathcal R}^\Lambda$, and the projection from ${\mathcal R}^\Lambda$ to $\hbox{Met}(M)$ is Fredholm of Fredholm index $-2n$.

The construction proceeds at first exactly like the proof of Proposition~7.2.1, the first change being that we replace the space $\hbox{Re} \left( {\cal J}({\bf L}) \right)$ in (\ref{E:tangentialjacobifieldsbranchtypelambda}) by the smaller space $\hbox{Re} \left( {\cal J}_0({\bf L})\right)$, $\hbox{Re} \left( {\cal J}_0({\bf L})\right)$ being the kernel of the map $\Psi _0$ as explained at the end of \S\ref{S:tangentialjacobi}.  As before, this $2(\nu - n)$-dimensional space can be translated about the coordinate neighborhood ${\mathcal U}$, thereby generating a subbundle ${\mathcal F}_0 \subseteq {\mathcal E}$ of fiber dimension $2(\nu - n)$.  We can then then define an $L^2$-orthogonal projection
$$p_0 : {\mathcal E} \longrightarrow {\mathcal F}_0^\bot$$
to the $L^2$-orthogonal complement ${\mathcal F}_0^\bot$ to ${\mathcal F}_0$ in ${\mathcal E}$.  Let ${\mathcal Z}$ be the image of the zero-section in the restriction of ${\mathcal F}_0^\bot$ to the submanifold ${\mathcal B}^\Lambda (\Sigma ,M)$ constructed in the earlier proof.

The hypothesis that $([f_0,{\bf p}_0, \eta_0], g_0)$ has weakly transversal crossings implies that the hypotheses of Lemmas~5.1.1 and 5.1.2 are satsified.  Therefore metric the variations constructed in the proof of Lemmas~5.1.1 exist and cover a complement to the smaller space $\hbox{Re}\left( {\cal J}_0({\bf L}) \right)$ within $\hbox{Re}\left( {\cal J}({\bf L}) \right)$.  This allows us to show that $p_0 \circ F^\Lambda $ is transverse to ${\mathcal Z}$ at the point $(f_0, \omega _0, g_0)$, and hence transverse to ${\mathcal Z}$ in some neighborhood ${\mathcal B}^\Lambda (\Sigma ,M) \cap {\mathcal W}$ of $(f_0, \omega _0, g_0)$, with ${\mathcal W} \subseteq {\mathcal U}$.  Then ${\mathcal R} = (p_0 \circ F^\Lambda )^{-1}({\mathcal Z}) \cap {\mathcal W}$ is a submanifold of ${\mathcal W}$ which contains all of the conformal harmonic maps of branch type $\Lambda $ within ${\mathcal W}$.

But ${\mathcal F}_0$ has codimension $2n$ in ${\mathcal F}$, and instead of (\ref{E:pLlambda}), we now obtain a Jacobi operator
$$p_0 \circ L^\Lambda : T_{(f,{\bf p}, \omega )}{\mathcal B}^\Lambda (\Sigma ,M) \longrightarrow \left({\mathcal F}_0^\bot\right)_{(f,{\bf p}, \omega)}$$
which has Fredholm index $-2n$.  Carrying through the rest of the argument, we find that the submanifold ${\mathcal R}^\Lambda $ we obtain as the inverse image of the zero section now has Fredholm projection to $\hbox{Met}(M)$ with Fredholm index $-2n$ instead of zero.  This finishes the proof when $\Sigma $ has genus at least two.

The case where $\Sigma $ has genus zero or one is now treated exactly as in the proof of Proposition 7.2.2.  By replacing ${\mathcal A}^n(\Sigma ,M)$ with a sequence of submanifolds
$${\mathcal A}^n_j(\Sigma ,M) = {\mathcal F}_j(\Sigma,M) \times \left( \Sigma ^n - \Delta^{(n)} \right) \times {\mathcal T}_g \times \hbox{Met}(M) $$
of codimension $\dim G$ in ${\mathcal A}^n(\Sigma ,M)$, and once again the previous argument yields a sequence of submanifolds of the configuration space which contain representatives of every $G$-orbit of minimal surfaces of branch type $\Lambda $ with weakly transversal crossings, the submanifolds projecting to the space of metrics with Fredholm index $-2n$.  This finishes the proof of Proposition~7.4.1 and the Main Theorem.

\section{Nonorientable surfaces}
\label{S:nonorientable}

Any complete theory of minimal surfaces must take nonorientable surfaces into account.  Indeed, the first applications we have in mind are to minimal spheres and tori, and these can cover prime minimal projective planes and Klein bottles.  Therefore we give a brief description of how the preceding results can be modified to treat nonorientable surfaces, such as projective planes and Klein bottles.

First, we note that the lemmas of \S\ref{S:selfintersection} apply directly to nonorientable surfaces.  However, in the rest of the argument, much of which utilizes complex analysis, it is expedient to pass from the nonorientable surfaces themselves to their orientable double covers.  We describe briefly how this works, beginning with projective planes and Klein bottles.

A minimal projective plane can be regarded as a conformal harmonic map
$$f : S^2 \longrightarrow M \quad \hbox{such that} \quad f \circ A = f,$$
where $A : S^2 \rightarrow S^2$ is the antipodal map.  (From the point of view of conformal geometry, there are many antipodal maps, but once we have reduced the symmetry group to $SO(3)$ by requiring that the center of mass is zero, as described in the Introduction (\ref{eq:centerofmass}), there is a canonical one that is an isometry with respect to our canonical choice of Riemannian metric on $S^2$.)  To put this in variational form, we note that the energy is invariant under composition of $f$ with $A$, so to find minimal projective planes, we need only find critical points for the restriction of the energy to the space
$$\hbox{Map}_0({\mathbb R}P^2 ,M) = \{ f \in \hbox{Map}_0(S^2 ,M) : f \circ A = f \}$$
and restrict the energy to a function
\begin{equation} E : \hbox{Map}_0({\mathbb R}P^2 ,M)  \longrightarrow {\mathbb R}. \label{eq:e1} \end{equation}
The action of $A$ on $\hbox{Map}_0(S^2 ,M)$, defined by $A(f) = f \circ A$, commutes with the action of $SO(3)$ in this case, so the energy (\ref{eq:e1}) is $SO(3)$-invariant.

Similarly, a minimal Klein bottle will be covered by a minimal torus with a flat metric of area one that is invariant under an orientation-reversing deck transformation $A_s$, for $s \in S^1$, which consists of a translation composed with a reflection, expressed in terms of appropriate standard coordinates $(t_1,t_2)$ on the torus as
$$A_s(t_1,t_2) = (t_1 + 1/2, - t_2 - s),$$
which one easily checks satisfies the identity $A_s^2 = 1$ modulo ${\mathbb Z}^2$.  Recall that the Teichm\"uller space ${\mathcal T}_1$ for the torus is the upper half-plane
$${\mathbb H} = \{ \omega = u + iv \in {\mathbb C} : v > 0 \},$$
with the point $\omega = u + iv$ corresponding to the conformal class of the torus ${\mathbb C}/\Lambda $, where $\Lambda $ is the lattice in ${\mathbb C}$ generated by $1$ and $\omega $.  In the case of a double cover of a Klein bottle, we arrange that the differential of $A_s$ fixes the generator corresponding to $1$ in the fundamental parallelogram, and the orthogonal line in the complex plane must then be preserved under the reflection $u+iv \mapsto u-iv$, so the \lq\lq Teichm\"uller space" of flat Klein bottles with total area one consists of the positive real numbers $\omega = iv$ with $v > 0$, the fixed point set of the involution
$$A_\star : {\mathbb H} \longrightarrow {\mathbb H}, \qquad A_\star(u + iv) = -u + iv.$$
(See Wolf \cite{Wo}, Proposition~2.5.8, for example.)  If we let
$$\hbox{Map}(K^2 ,M) = \{ f \in \hbox{Map}(T^2 ,M) : f \circ A_s = f \hbox{ for some } s \in S^1 \},$$
minimal Klein bottles can then be regarded as critical points for the restricted two-variable function
\begin{equation}E : \hbox{Map}(K^2 ,M) \times \{ \omega \in {\mathcal T}_1 : A_\star(\omega ) = \omega \} \longrightarrow {\mathbb R}. \label{eq:e2} \end{equation}
The energy is invariant under the action of $S^1 \times S^1$ on $\hbox{Map}(K^2 ,M)$ defined by
$$f(t_1,t_2) \mapsto f(t_1+s_1,t_2+s_2), \quad \hbox{for} \quad (s_1,s_2) \in S^1 \times S^1.$$
In the subsequent discussion, we let $A = A_{s_0}$ for some choice of $s_0 \in S^1$, thereby breaking part of the $S^1 \times S^1$-symmetry.

More generally, if $f_0 : \Sigma _0 \rightarrow M$ is a nonorientable parametrized minimal surface with oriented double cover $f : \Sigma \rightarrow M$ of genus at least two, we let $A : \Sigma \rightarrow \Sigma$ be the sheet interchange map.  Suppose now that $f : \Sigma \rightarrow M$ is an oriented double cover of a nonorientable minimal surface $f_0 : \Sigma _0 \rightarrow M$.  The map $A$ induces an involution on $A_\star$ on $f^*TM$, as well as on the space of sections of $f^*TM$, and both of these actions extend complex linearly to the complexifications.  Moreover, these involutions preserve both the metric and the pullback of the Levi-Civita connection.  If ${\bf E} = f^*TM \otimes {\mathbb C}$, we can therefore construct a direct sum decomposition $\Gamma ({\bf E}) = \Gamma _+({\bf E}) \oplus \Gamma _-({\bf E})$, where
$$ \Gamma _+({\bf E}) = \{ X \in \Gamma ({\bf E}) : A_\star (X) = \bar X \}, \quad \Gamma _-({\bf E}) = \{ X \in \Gamma ({\bf E}) : A_\star (X) = - \bar X \}.$$
The sections of $\Gamma _+({\bf E})$ can be regarded as deformations of the nonorientable minimal surface $f_0 : \Sigma \rightarrow M$.  The argument presented in \S\ref{S:harmonic} shows that if $f_0 : \Sigma \rightarrow M$ is a somewhere injective conformal harmonic map from a nonorientable surface to $M$, nontangential elements of $\Gamma _+({\bf E})$ can be perturbed away by variations of the metric on $M$.  Moreover, the proof of Proposition~3.1 can be modified to apply to somewhere injective nonorientable harmonic maps, so long as the line bundle ${\bf L}$ on the orientable double cover has only the holomorphic sections demanded by the group of symmetries.

To carry out the second variation formula presented in \S\ref{S:two-variableenergy}, we need to define a Teichm\"uller space ${\mathcal T}_{\Sigma _0}$ for the nonorientable surface $\Sigma _0$.  Of course, for the projective plane the Teichm\"uller space is trivial, while for the Klein bottle, we have seen that the Teichm\"uller space is
$$ {\mathcal T}_{K^2} = \{ \omega \in {\mathcal T}_1 : A_\star(\omega ) = \omega \} .$$
For surfaces of higher genus, we can define the Teichm\"uller space of $\Sigma _0$ as the quotient
$${\mathcal T}_{\Sigma _0} = \frac{\hbox{Met}_0(\Sigma _0)}{\hbox{Diff}_0(\Sigma _0)},$$
where $\hbox{Met}_0(\Sigma _0)$ is the space of constant curvature metrics of total area one and $\hbox{Diff}_0(\Sigma _0)$ is the space of diffeomorphisms of $\Sigma _0$ which are isotopic to the identity, just as in the orientable case.  It is easily checked that this agrees with the definition we have given in the case of the Klein bottle.

When $\Sigma _0$ is a nonorientable surface whose orientable double cover $\Sigma $ has genus $g \geq 2$, we can also describe the Teichm\"uller space ${\mathcal T}_{\Sigma _0}$ of $\Sigma _0$ as a subspace of the Teichm\"uller space ${\mathcal T}_{\Sigma }$ of $\Sigma $.  To see this, note first that $A_\star$ induces a conjugate linear map on quadratic differentials on $\Sigma $, and the cotangent space to ${\mathcal T}_{\Sigma _0}$ should be the space of holomorphic quadratic differentials $hdz^2$  on $\Sigma $ such that $A_\star(hdz^2) = \bar h d\bar z^2$.  As described in Chapter~2, \S 3 of \cite{SY2}, once one picks an origin $\eta _0 \in {\mathcal T}_\Sigma$, there is a bijection ${\mathcal H}$ from ${\mathcal T}_\Sigma$ to the space of holomorphic quadratic differentials on $(\Sigma ,\omega _0)$, where $\omega _0$ is the conformal equivalence class of $\eta _0$ defined as follows:  If $\eta $ is a hyperbolic metric of area one on $\Sigma $, then ${\mathcal H}(\eta )$ is the Hopf differential of the harmonic map from $(\Sigma ,\eta _0)$ to $(\Sigma , \eta)$, which exists by a theorem of Eells and Sampson and is unique by a theorem of Hartman.  This isomorphism allows us to transfer $A_\star$ to an involution on the Teichm\"uller space ${\mathcal T}_\Sigma$, the fixed point set of the involution being exactly those hyperbolic metrics which are preserved by the deck transformation $A$ and thus descend to the nonorientable surface $\Sigma _0$.  Hence we can set
$${\mathcal T}_{\Sigma _0} = \{ \omega \in {\mathcal T}_\Sigma  : A_\star (\omega ) = \omega \}.$$

We can now calculate the second variation of the energy
$$E : \hbox{Map}(\Sigma _0, M) \times {\mathcal T}_{\Sigma _0} \longrightarrow {\mathbb R},$$
by applying the formula of \S\ref{S:secondvariationtwovar} to the oriented double cover.  In the second variation for $\Sigma _0$, the holomorphic differentials, or the corresponding Beltrami differentials, must be invariant under the conjugate linear involution $A_\star$.  Since $A$ is orientation reversing,
$$\frac{\partial f}{\partial z}dz \mapsto \frac{\partial f}{\partial \bar z}d\bar z,$$
under $A_\star $, so $A_\star $ exchanges the line bundle ${\bf L}$ over the orientable double cover of $\Sigma _0$ with its conjugate.  The tangential Jacobi fields are now isomorphic to the pairs $(Z, hdz^2)$, where $Z$ is a section of ${\bf L}$ and $h dz^2$ is a holomorphic quadratic differential, which satisfy (\ref{E:tangentialJacobi1}) as well as the identities
$$A_\star (Z) = \bar Z, \qquad A_\star (hdz^2) = \bar h d\bar z^2.$$   
We now let
\begin{multline*} {\mathcal O}_A({\bf L}) = \{ Z \in {\mathcal O}({\bf L}) : A_\star (Z) = \bar Z \}, \\ {\mathcal J}_A({\bf L}) = \{ (Z, hdz^2) \in {\mathcal J}({\bf L}) : A_\star ((Z, hdz^2)) = (\bar Z, \bar hd\bar z^2) \},\end{multline*}
real linear subspaces of ${\mathcal O}({\bf L})$ and ${\mathcal J}({\bf L})$, respectively.  The elements of ${\mathcal J}_A({\bf L})$ correspond to tangential Jacobi fields which generate deformations of the nonorientable minimal surface, as opposed to deformations which would separate the two sheets of the orientable double cover.  It is now straightforward to verify that Proposition~4.4.1 holds for nonorientable surfaces without branch points. 

Each branch point of the nonorientable minimal surface $f_0 : \Sigma_0 \rightarrow M$ gives rise to two branch points of the orientable double cover $f$ which are interchanged by $A$, so if $\nu $ is the total branching order of $f_0$, the total branching order of $f$ will be $2 \nu$.  One can therefore check that in the case of real projective planes, the real dimension of ${\mathcal O}_A({\bf L})$ will be $3 + 2\nu$, while in the case of Klein bottles, the real dimension of ${\mathcal J}_A({\bf L})$ will be $1 + 2 \nu$.  For orientable surfaces whose double has genus at least two, the dimension of ${\mathcal J}_A({\bf L})$ will be $2 \nu$.

The argument for the Main Theorem can now be extended so that it applies to nonorientable surfaces.  In the proof, we need a Teichm\"uller space for the punctured nonorientable surface
$$\Sigma _0 - \{ \hbox{branch points $p_1, \ldots , p_n$} \}.$$
Each branch point $p_i$ lifts to two branch points $\tilde p_i$ and $\hat p_i$ in the oriented double cover which are interchanged by $A$.  If $\Sigma $ has genus at least two, the cotangent space to the Teichm\"uller space for $\Sigma _0$ consists of the meromorphic quadratic differentials $hdz^2$ on the oriented double cover with at most simple poles at the branch points $\tilde p_i$ and $\hat p_i$ such that
$$A_\star(hdz^2) = \bar h d \bar z^2.$$
In the case where $\Sigma = T^2$, we excise one additional pair of points from $T^2$ which are interchanged by $A$, while if $\Sigma = S^2$, we excise three pairs of points interchanged by $A$.  The cotangent space to Teichm\"uller space for the new punctured surface consists of the real linear subspace of the space of meromorphic quadratic differentials on the oriented double cover which are invariant under $A$ and have at most simple poles at the branch points and excised points.

Once these changes are made, it is relatively easy to check that the previous arguments yield a version of the Main Theorem for nonorientable minimal surfaces.  For generic choice of metric, prime minimal projective planes and Klein bottles are free of branch points and lie on nondegenerate critical submanifolds with dimensions dictated by the group actions.  For generic choice of metric, prime nonorientable minimal surfaces whose orientable double covers have genus at least two are Morse nondegenerate in the usual sense.

We can now give a rough description of the full collection of parametrized minimal spheres and tori in a Riemannian manifold with generic metric.  For minimal two-spheres, we first have the constant maps which make up a nondegenerate critical submanifold at energy level zero.  Then we have the prime minimal two-spheres themselves.  Next we have the branched covers of minimal two-spheres, which can be thought of as more complicated parametrizations of the prime minimal two-spheres.  Finally, we have prime minimal projective planes as well as covers of these, including two-to-one covers of projective planes in which the parameter space is the two-sphere and other branched covers of prime minimal projective planes.

For the tori, we have the constant maps and the prime minimal tori once again.  This time, however, there are two types of covers, ordinary covers of prime minimal tori and Klein bottles without branch points, and branched covers of prime minimal two-spheres and projective planes.

When the genus is at least two, a complete critical point theory for minimal surfaces of genus $g$ would have to consider prime surfaces of genus $g$, branched covers of prime minimal surfaces of genus $\leq g$ and branched covers of prime nonorientable surfaces.

\section{Families of metrics}
\label{S:families}

The argument we presented actually establishes an important extension of the Main Theorem.  Recall that for each branch type $\Lambda $ we have a subspace ${\mathcal P}^\Lambda$ of ${\mathcal A}^n(\Sigma ,M)$ consisting of prime parametrized minimal surfaces of branch type $\Lambda $, where $n$ is the number of branch points in $\lambda $, and we have a natural projection
$$\pi : {\mathcal A}^n(\Sigma ,M) \longrightarrow \hbox{Met}(\Sigma ).$$
But Propositions 7.3.2 and 7.4.1 show that for $g$ belonging to a subspace of $\hbox{Met}(\Sigma )$ of codimension $\geq 2$, $\pi ^{-1}(g)$ is empty.  Thus if we avoid a subspace of metrics of codimension two, we can work with the much simpler space of minimal immersions with no branch points,
\begin{multline*}{\mathcal P}^\emptyset = \{ (f,\omega ,g) \in \hbox{Map}(\Sigma ,M) \times {\mathcal T} \times \hbox{Met}(M) : \\
\hbox{ $f$ is a prime minimal immersion, conformal with respect to $\omega $ and $g$ } \}.\end{multline*}
This is a submanifold of $ \hbox{Map}(\Sigma ,M) \times {\mathcal T} \times \hbox{Met}(M)$ with Fredholm projection
\begin{equation} \pi : {\mathcal P}^\emptyset \longrightarrow \hbox{Met}(M) \label{E:piforimmersions}\end{equation}
of Fredholm index $d_\Sigma $, where $d_\Sigma $ is the real dimension of the group of symmetries.

Suppose now that $g_0$ and $g_1$ are generic metrics and $\gamma : [0,1] \rightarrow \hbox{Met}(M)$ is a smooth path with $\gamma (0) = g_0$ and $\gamma (1) = g_1$.  We now modify the argument for Theorem~3.1 of \cite{Sm2} to construct a smooth approximation $\tilde \gamma $ to $\gamma $ with $\tilde \gamma (0) = g_0$ and $\tilde \gamma (1) = g_1$ such that $\tilde \gamma $ is transverse to the projection $\pi$ of (\ref{E:piforimmersions}) as well as finitely many of the Fredholm projections of Theorems 7.3.2 and 7.4.1.  (Compactness of the interval implies that only finitely many of the Fredholm projections are needed.)

\vskip .1in
\noindent
{\bf Theorem 9.1.}  {\sl Let $M$ be a compact connected smooth manifold of dimension at least four.  A generic path $\gamma : [0,1] \rightarrow \hbox{Met}(M)$ between generic metrics on $M$ will have the property that the only prime minimal surfaces $f : \Sigma \rightarrow M$ for metrics in $\gamma ([0,1])$ have no branch points.  Moreover, $\pi ^{-1}(\gamma ([0,1]))$ will be a smooth submanifold of dimension $d_\Sigma + 1$ of ${\mathcal P}^\emptyset$, where $d_\Sigma $ is the dimension of the group $G$ of symmetries of the energy $E$.  Finally, we can arrange that each element of $\pi ^{-1}(\gamma ([0,1]))$ be an immersion with transversal crossing if $\dim M \neq 5$, or an immersion with weakly transversal crossings if $\dim M = 5$.}

\vskip .1in
\noindent
The last statement follows from Propositions~7.3.1 and 7.3.2.

\end{document}